\documentclass[11pt]{article}
\usepackage{latexsym,bm}
\usepackage{color}
\usepackage{amssymb, graphicx, amsmath}  

\newcommand\qed{\qquad $\square$}

\def \[{\begin{equation}}
\def \]{\end{equation}}

\newtheorem{theorem}{Theorem}[section]

\newtheorem{definition}{Definition}[section]

\newtheorem{lemma}{Lemma}[section]

\newtheorem{corollary}{Corollary }[section]

\begin{document}

\begin{center}
{\bf The existence for the classical solution of the Navier-Stokes equations}

\medskip

  { Jianfeng Wang  }\\
Department of Mathematics, Hohai University, Nanjing, 210098,
  P.R. China.\\Email: wjf19702014@163.com, Acadmic email: 20020001@hhu.edu.cn.

\end{center}
\bigskip
{\narrower \noindent  {\bf Abstract.} In this paper we will discuss the existence for the classical solution of the Navier-Stokes equations. First, we transform it into generalized integral equations. Next, we discuss the existence of the classical solution by Leray-Schauder degree and Sobolev space\ $H^{-m_{1}}(\Omega_{1})$. 
    \par }
\vskip 1.0 true cm
 \section{Introduction}
We consider the dynamical equations for a viscous and
incompressible fluid as follows, {
\begin{eqnarray}
 &&\cfrac{\partial u_{1}}{\partial x}+\cfrac{\partial u_{2}}{\partial y}+\cfrac{\partial u_{3}}{\partial z}=0, \\&&\cfrac{\partial u}{\partial
 t}-\mu \Delta u+u_{1}\cfrac{\partial u}{\partial
 x}+u_{2}\cfrac{\partial u}{\partial
 y}+u_{3}\cfrac{\partial u}{\partial
 z}+\cfrac{1}{\rho}\
 grad p=F,\end{eqnarray}}where\ $u=(u_{1},\ u_{2},\ u_{3})^{T}$\ is the velocity vector,\ $\mu$\ is the
dynamic viscosity coefficient,\ $\rho$\
is the density of the material particle,\ $p$\ is intensity of the pressure,\ $F=(F_{1},\ F_{2},\ F_{3})^{T}$\ is the density of the body
force,\ $u=(u_{1},\ u_{2},\ u_{3})^{T},\ p$\ and\ $F$\ are all continuous functions with the variables of
the time\ $t$\ and position\ $(x,\ y,\ z)^{T}$. These equations are called the Navier-Stokes
equations. \\We assume\ $u=(u_{1},\ u_{2},\ u_{3})^{T}\in C^{2}(\overline{\Omega})\cap C^{1}[0,\ T]$,\ $p\in C^{1}(\overline{\Omega})\cap C[0,\ T]$,\ $\ t\in [0,\ T],\ (x,\ y,\ z)^{T}\in \overline{\Omega}\subset R^{3},\ \overline{\Omega}=\Omega\cup\partial\Omega$,\ $\Omega$\ is a bounded domain, and\ $\partial\Omega\in C^{1,\ \alpha},\ 0<\alpha\leq 1$, moreover,\ $\overline{\Omega}$\ is convex, the initial conditions should be known as\[u|_{t=0}\in C^{2}(\overline{\Omega}),\]and at least one of three boundary conditions as follows,\\(1)Dirichlet problem,\[u|_{\partial\Omega\times[0,\ T]}\in C^{2}(\partial\Omega)\cap C^{1}[0,\ T],\](2)Neumann problem,\[\cfrac{\partial u}{\partial n}|_{\partial\Omega\times[0,\ T]}\in C^{1}(\partial\Omega\times[0,\ T]),\] where\ $n=(n_{1},\ n_{2},\ n_{3})^{T}$\ is the exterior normal vector to\ $\partial\Omega$,\\(3)Robin problem,\[(\cfrac{\partial u}{\partial n}+\sigma u)|_{\partial\Omega\times[0,\ T]}\in C(\partial\Omega\times[0,\ T]),\] where\ $\sigma=\sigma(x,\ y,\ z,\ t)$\ is continuous, and\ $\sigma(x,\ y,\ z,\ t)>0,\ \forall (x,\ y,\ z,\ t)^{T}\in\partial\Omega\times[0,\ T].$\\Our goal is absolutely the classical solution of Eqs(1.1) and (1.2).\\
In Section 2, we will prove that the boundary in\ $C^{1,\ \alpha},\ 0<\alpha\leq 1$, is not only the Sobolev's imbedding surface, but also the Lyapunov's surface, the Hopf's surface for the strong maximum principle.\\
In section 3, we will transform Eqs(1.1) and (1.2) into the equivalent generalized integral equations as follows,\[ Z_{1}=T_{0}(Z_{1})=w_{1}(x,\ y,\ z,\ t)+w_{2}(x,\ y,\ z,\ t).\ast(\psi(Z_{1})),\]
where\ $Z_{1}=(u,\ p,\ \partial u\setminus\ u_{1x},\ \partial^{2} u, grad p)^{T}$,\ $\partial u=(u_{x},\ u_{y},\ u_{z})^{T},\ \partial u\setminus\ u_{1x}=(u_{2x},\ u_{3x},\ u_{y},\ u_{z})^{T}$,\ $u_{jx}=\partial u_{j}/\partial x,\ j=1,\ 2,\ 3$,\ $\partial^{2}u=(u_{xx},\ u_{xy},\ u_{xz},\ u_{yy},\ u_{yz},\ u_{zz})^{T}$,\ $.\ast$\ means the matrix convolution as follows,
\[ w_{2}.\ast(\psi(Z_{1}))=\int_{R^{4}}w_{2}(x-x_{1},\ y-y_{1},\ z-z_{1},\ t-\tau)(\psi(Z_{1}(x_{1},\ y_{1},\ z_{1},\ \tau)))dx_{1}dy_{1}dz_{1}d\tau. \]
Here generalized comes from\ $w_{1}$\ and\ $w_{2}$\ being related to the Dirac function. The term equivalent is defined in the following.
\begin{definition}\label{definition}The equations\ $f_{1}(x)=0$\ and\ $f_{2}(y)=0$\ are equivalent, if and only if there exist continuous mappings\ $T_{1},\ T_{2}$, such that\ $\forall x,\ y$\ are respectively the solutions of\\ $f_{1}(x)=0,\ f_{2}(y)=0$,\ we have\ $f_{1}(T_{2}(y))=0,\ f_{2}(T_{1}(x))=0,$\ and\ $T_{2}(T_{1}(x))=x$\ ,\ $T_{1}(T_{2}(y))=y$.\end{definition}
Hence, we can get a necessary and sufficient condition for there exist\ $u\in C^{2}(\overline{\Omega})\cap C^{1}[0,\ T]$,\ $p\in C^{1}(\overline{\Omega})\cap C[0,\ T]$\ satisfy Eqs(1.1) and (1.2) is that there exists\ $Z_{1}\in C(\overline{\Omega}\times[0,\ T])$\ satisfies\ $Z_{1}=T_{0}(Z_{1})$.\\
In Section 4, we will discuss the existence of the classical solution of\ $Z_{1}=T_{0}(Z_{1})$. We will use the theory on\ $H^{-m_{1}}(\Omega_{1})$, where\ $\Omega_{1}=\Omega\times(0,\ T)$, which is defined on page 130 in [10], and a primary theory on the Leray-Schauder degree.\\
At first, we construct a norm\ $\|\cdot\|_{-m_{1}}$\ for\ $T_{0}(Z_{1})=(T_{0,\ i}(Z_{1}))_{33\times 1}$\ as the following,\[ \|T_{0}(Z_{1})\|_{-m_{1}}=\max_{1\leq i\leq 33}\sup_{\varphi\in C_{0}^{\infty}(\Omega_{1})}\cfrac{|<T_{0,\ i}(Z_{1}),\ \varphi>|}{\|\varphi\|_{m_{1}}}.\]
Next, we make approximate ordinary integral equations\ $Z_{1}=T_{0\epsilon}(Z_{1}),\ \forall \epsilon>0$, where
\[T_{0\epsilon}(Z_{1})=\delta_{\epsilon}.\ast T_{0}(Z_{1}),\ \delta_{\epsilon}=\cfrac{1}{(\sqrt{\pi\epsilon})^{4}}e^{-|X|^{2}/\epsilon},\ X=(x,\ y,\ z,\ t)^{T}.\]
At last, we assume the following,
\[ \exists M>0,\ \exists \delta>0,\ \exists \delta^{\prime}>0,\ \forall \epsilon\in (0,\ \delta],\ \forall T\in (0,\ \delta^{\prime}],\ \mbox{we have}\ \tau(M,\ \epsilon,\ T)>0, \]
\[ \exists M>0,\ \exists \delta>0,\ \forall \epsilon\in (0,\ \delta],\ \tau(M,\ \epsilon)>0,\ \mbox{and}
\ \exists \epsilon_{0}\in (0,\ \delta],\ deg(Z_{1}-T_{0\epsilon_{0}}(Z_{1}),\ \Omega_{M},\ 0)\neq0,\]\[ \exists\ \epsilon_{k}>0,\ Z_{1\epsilon_{k}}\in \Omega_{M},\ Z_{1\epsilon_{k}}=T_{0\epsilon_{k}}(Z_{1\epsilon_{k}}),\ k\geq1,\ \lim_{k\rightarrow +\infty}\epsilon_{k}=0,\ \mbox{and}\ \sup_{k,\ l,\ i}S(\partial\Omega_{k,\ l,\ i}^{+})< +\infty,\]where$$\tau(M,\ \epsilon,\ T)=\inf_{\|Z_{1}\|_{\infty}=M}\|Z_{1}-T_{0\epsilon}(Z_{1})\|_{\infty},$$
if time\ $T$\ is fixed, then we denote\ $\tau(M,\ \epsilon,\ T)$\ into\ $\tau(M,\ \epsilon)$,\ $\Omega_{M}=\{Z_{1}\in C(\overline{\Omega_{1}}): \|Z_{1}\|_{\infty}<M\}$,  $$S(\partial\Omega_{k,\ l,\ i}^{+})=\int_{\{X\in \overline{\Omega}\times[0,\ T]:\ Z_{1,\ i,\ \epsilon_{k}}-Z_{1,\ i,\ \epsilon_{l}}= 0\}}dS,\ k\neq l,\ 1\leq i\leq 33. $$We will prove that the strong solution of\ $Z_{1}=T_{0}(Z_{1})$\ will exist locally under the condition (1.11) and exist globally under the condition (1.12), where the strong solution\ $Z_{1}^{\ast}$\ is required that there exist series\ $\epsilon_{k}\rightarrow0,\ k\rightarrow+\infty$, such that\[ \lim_{k\rightarrow+\infty}\|Z_{1\epsilon_{k}}-Z_{1}^{\ast}\|_{-m_{1}}=0,\ \mbox{moreover}\ \lim_{k\rightarrow+\infty}\|Z_{1\epsilon_{k}}-T_{0}(Z_{1\epsilon_{k}})\|_{-m_{1}}=0.\]
Moreover, the\ $L^{\infty}$\ solution of\ $Z_{1}=T_{0}(Z_{1})$\ will exist under the condition (1.13).\\
If the strong solution is locally integrable, then it is in\ $L^{\infty}$. If (1.11) or (1.12) or (1.13) does not hold, then the blow-up will happen.\\ The\ $L^{\infty}$\ solution of\ $Z_{1}=T_{0}(Z_{1})$\ will always exist and be unique except the blow-up.\\
We obtain\ $u^{\ast}\in W^{2,\ +\infty}(\overline{\Omega}),\ p^{\ast}\in W^{1,\ +\infty}(\overline{\Omega})$, if\ $Z_{1}^{\ast}=T_{0}(Z_{1}^{\ast})$,\ $Z_{1}^{\ast}\in L^{\infty}(\overline{\Omega}\times[0,\ T])$, where\ $Z_{1}^{\ast}=(u^{\ast},\ p^{\ast},\ \partial u^{\ast}\setminus\ u^{\ast}_{1x},\ \partial^{2} u^{\ast}, grad p^{\ast})^{T}$. Here\ $W^{1,\ +\infty}(\overline{\Omega}),\ W^{2,\ +\infty}(\overline{\Omega})$\ are Sobolev spaces defined on page 153 in [2]. From the condition that domain\ $\Omega$\ satisfies a uniform exterior and interior cone, if\ $\Omega$\ is bounded,\ $\partial\Omega\in C^{1,\ \alpha},\ 0<\alpha\leq 1$, we can get that\ $u^{\ast}\in C^{1,\ 1}(\overline{\Omega}),\ p^{\ast}\in C^{0,\ 1}(\overline{\Omega})$\ by imbedding. By using Morrey's inequality defined on page 163 in [2], we get\ $u^{\ast}$\ is twice classically differentiable and\ $p^{\ast}$\ is classically differentiable almost everywhere in\ $\overline{\Omega}$.\\If\ $F(T_{0}(Z_{1}^{\ast}))$\ is analytical, then\ $u^{\ast}$\ and\ $p^{\ast}$\ satisfy Eqs(1.1) and (1.2) almost everywhere in\ $\overline{\Omega}\times[0,\ T]$, where\ $F(T_{0}(Z_{1}^{\ast}))$\ is the Fourier transform of\ $T_{0}(Z_{1}^{\ast})$. That is near our goal.\\
Since we haven't got any similar paper, we have to put some books we have learnt in the reference.\\
\section{Boundary in\ $C^{1,\ \alpha}$}\setcounter{equation}{0}
In this section, we will explain why we chose\ $\partial\Omega\in C^{1,\ \alpha},\ 0<\alpha\leq 1,$\ instead of\ $\partial\Omega\in C^{2,\ \alpha},\ 0<\alpha\leq 1,$\ or the Sobolev's imbedding surface, the Lyapunov's surface, the Hopf's surface for the strong maximum principle.\\First of all, we know\ $C^{1,\ \alpha}\supset C^{2,\ \alpha},\ 0<\alpha\leq 1$.\ Simplicity and generality are our eternal pursuit. Secondly,\ $\partial\Omega\in C^{1,\ \alpha},\ 0<\alpha\leq 1$\ is just the Sobolev's imbedding surface from the following.
\begin{theorem} \label{Theorem2-1} If\ $\Omega$\ is bounded,\ $\partial\Omega\in C^{1,\ \alpha},\ 0<\alpha\leq 1,$\ then domain\ $\Omega$\ satisfies a uniform exterior cone condition, that is, there exists a fixed finite right circular cone\ $K$,\ such that each\ $P=(x,\ y,\ z)^{T}\in \partial\Omega$\ is the vertex of a cone\ $K(P)$,\ $\overline{K(P)}\cap\overline{\Omega}=P$, and\ $K(P)$\ is congruent to\ $K$. \end{theorem}
{\it Proof of theorem 2.1}. By using the equivalent definition of\ $C^{1,\ \alpha},\ 0<\alpha\leq 1,$\ on page 94 in [2], if\ $\partial\Omega\in C^{1,\ \alpha}$,\ $0<\alpha\leq 1,$\ then each\ $P_{0}=(x_{0},\ y_{0},\ z_{0})^{T}\in \partial\Omega$,\ there exists a neighborhood\ $U(P_{0},\ \delta_{0}(P_{0})),\ \delta_{0}(P_{0})>0$, where
 $$U(P_{0},\ \delta_{0}(P_{0}))=\{P:\ |P-P_{0}|<\delta_{0}(P_{0})\},\ \overline{U}(P_{0},\ \delta_{0}(P_{0}))=\{P:\ |P-P_{0}|\leq\delta_{0}(P_{0})\},$$
$P=(x,\ y,\ z)^{T}\in R^{3},\ \mid P-P_{0}\mid=\sqrt{(x-x_{0})^{2}+(y-y_{0})^{2}+(z-z_{0})^{2}}$, and\ $U(P_{0},\ \delta_{0}(P_{0}))\cap\partial\Omega$\ is a graph of a\ $C^{1,\ \alpha},\ 0<\alpha\leq 1,$\ function of two of the coordinates\ $x,\ y,\ z$.\\ Without loss of the generality, we assume such a function is\ $f_{0}(x,\ y)\in C^{1,\ \alpha},\ 0<\alpha\leq 1$, and we have the following,\ $z-f_{0}(x,\ y)=0$,\ if\ $P=(x,\ y,\ z)^{T}\in U(P_{0},\ \delta_{0}(P_{0}))\cap\partial\Omega$,\ $z-f_{0}(x,\ y)>0$,\ if\ $P=(x,\ y,\ z)^{T}\in U(P_{0},\ \delta_{0}(P_{0}))\setminus\overline{\Omega}$, and\ $z-f_{0}(x,\ y)<0$,\ if\ $P=(x,\ y,\ z)^{T}\in U(P_{0},\ \delta_{0}(P_{0}))\cap\Omega$.\\ We will get the same results if\ $U(P_{0},\ \delta_{0}(P_{0}))\cap\partial\Omega$\ is a graph of other functions.\\ We introduce the following lemma.
\begin{lemma} \label{lemma1} If\ $\Omega$\ is bounded,\ $\partial\Omega\in C^{1,\ \alpha},\ 0<\alpha\leq 1,$\ then for each\ $P_{0}=(x_{0},\ y_{0},\ z_{0})^{T}\in \partial\Omega$,\ there exists a\ $C(P_{0})>0$, related to\ $f_{0}(x,\ y)$, such that\[\mid {\bf r}_{P_{0}P}\cdot n_{p_{0}}\mid\leq C(P_{0})\mid P-P_{0}\mid^{1+\alpha},\ \forall P=(x,\ y,\ z)^{T}\in U(P_{0},\ \delta_{0}(P_{0}))\cap\partial\Omega,\]where\ ${\bf r}_{P_{0}P}=\overrightarrow{P_{0}P},\ n_{p_{0}}$\ is exterior normal vector to\ $\partial\Omega$\ at point\ $P_{0}$.\end{lemma}
{\it Proof of lemma 2.1}. From\ $n_{p_{0}}=\cfrac{1}{\sqrt{f_{0x}^{2}(x_{0},\ y_{0})+f_{0y}^{2}(x_{0},\ y_{0})+1}}\ (-f_{0x}(x_{0},\ y_{0}),\ -f_{0y}(x_{0},\ y_{0}),\ 1)^{T}$, where\ $f_{0x},\ f_{0y}$\ are partial derivatives of\ $f_{0}$, we get\[{\bf r}_{P_{0}P}\cdot n_{p_{0}}=\cfrac{z-z_{0}-(x-x_{0})f_{0x}(x_{0},\ y_{0})-(y-y_{0})f_{0y}(x_{0},\ y_{0})}{\sqrt{f_{0x}^{2}(x_{0},\ y_{0})+f_{0y}^{2}(x_{0},\ y_{0})+1}}.\]From\ $z-z_{0}=f_{0}(x,\ y)-f_{0}(x_{0},\ y_{0})$, and\ $f_{0}(x,\ y)\in C^{1,\ \alpha}$, we have\[z-z_{0}=(x-x_{0})f_{0x}(tx_{0}+(1-t)x,\ ty_{0}+(1-t)y)+(y-y_{0})f_{0y}(tx_{0}+(1-t)x,\ ty_{0}+(1-t)y),\] where\ $0\leq t\leq 1$. From\ $f_{0}(x,\ y)\in C^{1,\ \alpha}$, we know there exists\ $C_{\alpha}(P_{0})>0$, related to\ $f_{0}(x,\ y)$, such that\ $\forall P_{1}=(x_{1},\ y_{1},\ z_{1})^{T},\ P_{2}=(x_{2},\ y_{2},\ z_{2})^{T}\in \overline{U}(P_{0},\ \delta_{0}(P_{0}))\cap\partial\Omega$,\begin{eqnarray*}&&\mid f_{0x}(x_{1},\ y_{1})-f_{0x}(x_{2},\ y_{2})\mid \leq C_{\alpha}(P_{0})\mid (P_{1}-P_{2})_{1}\mid^{\alpha},\\ &&\mid f_{0y}(x_{1},\ y_{1})-f_{0y}(x_{2},\ y_{2})\mid \leq C_{\alpha}(P_{0})\mid (P_{1}-P_{2})_{1}\mid^{\alpha},\end{eqnarray*} where\
$(P_{1}-P_{2})_{1}=(x_{1}-x_{2},\ y_{1}-y_{2},\ 0)^{T}$.\\
Hence we obtain the following,\begin{eqnarray}&&\mid f_{0x}(tx_{0}+(1-t)x,\ ty_{0}+(1-t)y)-f_{0x}(x_{0},\ y_{0})\mid \leq C_{\alpha}(P_{0})\mid (P-P_{0})_{1}\mid^{\alpha},\\ &&\mid f_{0y}(tx_{0}+(1-t)x,\ ty_{0}+(1-t)y)-f_{0y}(x_{0},\ y_{0})\mid \leq C_{\alpha}(P_{0})\mid (P-P_{0})_{1}\mid^{\alpha},\end{eqnarray} where\
$(P-P_{0})_{1}=(x-x_{0},\ y-y_{0},\ 0)^{T}$.\\ From\ $\sqrt{f_{0x}^{2}(x_{0},\ y_{0})+f_{0y}^{2}(x_{0},\ y_{0})+1}\geq 1$, we obtain,
\[\mid{\bf r}_{P_{0}P}\cdot n_{p_{0}}\mid\leq 2C_{\alpha}(P_{0})\mid (P-P_{0})_{1}\mid^{1+\alpha}\leq C(P_{0})\mid P-P_{0}\mid^{1+\alpha},\] where\ $C(P_{0})=2C_{\alpha}(P_{0})$.
\qed\\
Now we can make an exterior finite right circular cone\ $K(P_{0})$\ at point\ $P_{0}$. We choose a\ $\delta_{1}(P_{0})<\delta_{0}(P_{0})$\ that is small enough such that\ $C(P_{0})(2\delta_{1}(P_{0}))^{\alpha}<1$.\ We let\ $P_{0}$\ be the vertex of a cone\ $K(P_{0})$\ and\ $n_{p_{0}}$\ be the symmetry axis. The polar angle is\ $\theta(P_{0})=\arccos(C(P_{0})(2\delta_{1}(P_{0}))^{\alpha})\in(0,\ \pi/2)$,\ and the length of generatrix is\ $\delta_{1}(P_{0})/3$. \\We will prove that\ $\overline{K(P_{0})}\cap\overline{\Omega}=P_{0}$. \\ From\ $\mid P-P_{0}\mid\leq \delta_{1}(P_{0})/3,\ \forall P\in \overline{K(P_{0})}$, we can obtain\ $\overline{K(P_{0})}\subset U(P_{0},\ \delta_{1}(P_{0}))$.\\ And from\[{\bf r}_{P_{0}P}\cdot n_{p_{0}}\geq\mid P-P_{0}\mid C(P_{0})(2\delta_{1}(P_{0}))^{\alpha}> C(P_{0})\mid P-P_{0}\mid^{1+\alpha},\ \forall P\in \overline{K(P_{0})},\ P\neq P_{0}, \]we can obtain\ $\overline{K(P_{0})}\cap \partial\Omega=P_{0}$.\\ Finally, if there exists\ $P\in \overline{K(P_{0})}\cap\Omega$,\ then we can obtain\ $z<f_{0}(x,\ y)$, where\ $P=(x,\ y,\ z)^{T}$. From (2.2), we can get\ ${\bf r}_{P_{0}P}\cdot n_{p_{0}}< C(P_{0})\mid P-P_{0}\mid^{1+\alpha}$. This contradicts (2.7). Hence\ $\overline{K(P_{0})}\cap\overline{\Omega}=P_{0}$.\\
Then, we can make a uniform exterior finite right circular cone\ $K$\ for each\ $P\in U(P_{0},\ \delta_{1}(P_{0})/3)\cap\partial\Omega$.\ $K$\ is congruent to\ $K(P_{0})$, if we let\ $P$\ be the vertex of a cone\ $K(P)$,\ $n_{p}$\ the symmetry axis, a polar angle of\ $\theta(P)=\theta(P_{0})=\arccos(C(P_{0})(2\delta_{1}(P_{0}))^{\alpha})$,\ and the length of generatrix\ $\delta_{1}(P_{0})/3$. From lemma 2.1, we obtain the following,\ $\forall P_{1}=(x_{1},\ y_{1},\ z_{1})^{T},\ P_{2}=(x_{2},\ y_{2},\ z_{2})^{T}\in \overline{U}(P_{0},\ \delta_{0}(P_{0}))\cap\partial\Omega$,
\[\mid{\bf r}_{P_{2}P_{1}}\cdot n_{p_{2}}\mid\leq C(P_{0})\mid P_{1}-P_{2}\mid^{1+\alpha}.\]
If there exists\ $P^{\prime}\neq P\in K(P)\cap \overline{\Omega}$, then it will contradict the following,
\[{\bf r}_{PP^{\prime}}\cdot n_{p}\geq\mid P^{\prime}-P\mid C(P_{0})(2\delta_{1}(P_{0}))^{\alpha}> C(P_{0})\mid P^{\prime}-P\mid^{1+\alpha},\ \forall P^{\prime}\in \overline{K(P)},\ P^{\prime}\neq P. \]
From arbitrary\ $P_{0}$, we obtain that\[\bigcup_{P_{0}\in \partial\Omega}[U(P_{0},\ \delta_{1}(P_{0})/3)\cap\partial\Omega]\] is an open cover for\ $\partial\Omega$. From the Heine-Borel theorem, we see that there exists a finite sub-cover for\ $\partial\Omega$ as follows,\ $\exists N>0,\ \exists P_{k}\in\partial\Omega,\ \exists \delta_{1}(P_{k})>0,\ 1\leq k\leq N$, such that
\[ \bigcup_{k=1}^{N}[U(P_{k},\ \delta_{1}(P_{k})/3)\cap\partial\Omega]\supset\partial\Omega,\] and the definitions of\ $\delta_{0}(P_{k}),\ C_{\alpha}(P_{k}),\ C(P_{k}),\ \delta_{1}(P_{k}),\ \theta(P_{k})$\ are the same as\ $\delta_{0}(P_{0}),\ C_{\alpha}(P_{0}),\ C(P_{0})$,\ $\delta_{1}(P_{0}),\ \theta(P_{0})$,\ $1\leq k\leq N$.\\ So we can make a uniform exterior finite right circular cone\ $K$\ for each\ $P\in \partial\Omega$.\ All the\ $K(P)$\ are congruent, if we let\ $P$\ be the vertex of a cone\ $K(P)$,\ $n_{p}$\ be the symmetry axis, have a polar angle of\ $\theta^{\ast}$,\ and length of generatrix of\ $\delta_{1}^{\ast}$, where
 $$ \theta^{\ast}=\min_{1\leq k\leq N}\theta(P_{k}),\ \delta_{1}^{\ast}=\min_{1\leq k\leq N}\delta_{1}(P_{k})/3.$$ Hence the statement stands. \qed
\\If we choose\ $-n_{p_{0}}$\ as the symmetry axis of the cone, then we can transform the uniform exterior cone into a uniform interior cone. Hence, domain\ $\Omega$\ satisfies a uniform exterior and interior cone condition, if\ $\Omega$\ is bounded,\ $\partial\Omega\in C^{1,\ \alpha},\ 0<\alpha\leq 1.$\ It looks like\ $C^{1,\ \alpha},\ 0<\alpha\leq 1,$\ is better than\ $C^{0,\ 1}$\ for the Sobolev's imbedding. However, we do not discuss the weak solution of Eqs(1.1) and (1.2) directly here. We will be discussing the classical solution directly. We want to obtain the equivalent equations that the classical solution should satisfy.\\ Thirdly, the following two theorems demonstrate that\ $\partial\Omega\in C^{1,\ \alpha},\ 0<\alpha\leq 1,$\ may be taken as the Lyapunov's surface. We will use them in section 3. We require two lemmas as follows.\begin{lemma} \label{lemma2.2} If\ $\Omega$\ is bounded,\ $\partial\Omega\in C^{1,\ \alpha},\ 0<\alpha\leq 1$,\ $f(M,\ P),\ M\neq P$\ is continuous, and\ $\forall M_{0}\in \partial\Omega,\ \forall \epsilon>0,\ \exists \delta>0$,\ such that\ $\forall M\in U(M_{0},\ \delta)$, we have
$$|\int_{U(M_{0},\ \delta)\cap\partial\Omega}f(M,\ P)dS_{P}|\leq \epsilon, $$ then \[\omega(M)=\int_{\partial\Omega}f(M,\ P)dS_{P},\ \forall M=(x,\ y,\ z)^{T}\in R^{3},\] will be continuous. In particular, if\ $\forall M_{0}\in \partial\Omega$,\ there exists a neighbourhood\ $U(M_{0})$,\ and\ $\delta_{1}\in(0,\ 1],\ C>0$, such that\[ |f(M,\ P)|\leq \cfrac{C}{r_{MP}^{2-\delta_{1}}}, \ \forall M\in U(M_{0}),\ \forall P\in  U(M_{0})\cap\partial\Omega ,\]where\ $r_{MP}=|M-P|$, then\ $\omega(M)$\ will also be continuous.\end{lemma}The proof is available on pages 178 to 180 of [6].
\begin{lemma} \label{lemma2.3} If\ $\Omega$\ is bounded,\ $\partial\Omega\in C^{1,\ \alpha},\ 0<\alpha\leq 1,$\ then there exist\ $\delta_{0}>0,\ C_{0}>0$, such that for each\ $P_{0}=(x_{0},\ y_{0},\ z_{0})^{T}\in \partial\Omega$,\ we have the following,\[\mid n_{p}-n_{p_{0}}\mid\leq C_{0}r_{PP_{0}}^{\alpha},\ \forall P=(x,\ y,\ z)^{T}\in U(P_{0},\ \delta_{0})\cap\partial\Omega.\]\end{lemma}
{\it Proof of lemma 2.3}. If\ $\Omega$\ is bounded,\ $\partial\Omega\in C^{1,\ \alpha},\ 0<\alpha\leq 1,$\ then from Theorem 2.1, for each\ $P_{0}=(x_{0},\ y_{0},\ z_{0})^{T}\in \partial\Omega$,\ there exist\ $P_{k}\in \partial\Omega,\ 1\leq k\leq N$, such that\ $P_{0}\in U(P_{k},\ \delta_{1}(P_{k})/3)\cap\partial\Omega$. And\ $U(P_{k},\ \delta_{1}(P_{k}))\cap\partial\Omega$\ is a graph of a\ $C^{1,\ \alpha},\ 0<\alpha\leq 1,$\ function of two of the coordinates\ $x,\ y,\ z$. \\Without loss of the generality, we assume such a function is\ $f_{k}(x,\ y)\in C^{1,\ \alpha},\ 0<\alpha\leq 1$, and we have the following,\ $z-f_{k}(x,\ y)=0$,\ if\ $P=(x,\ y,\ z)^{T}\in U(P_{k},\ \delta_{1}(P_{k}))\cap\partial\Omega$,\ $z-f_{k}(x,\ y)>0$,\ if\ $P=(x,\ y,\ z)^{T}\in U(P_{k},\ \delta_{1}(P_{k}))\setminus\overline{\Omega}$, and\ $z-f_{k}(x,\ y)<0$,\ if\ $P=(x,\ y,\ z)^{T}\in U(P_{k},\ \delta_{1}(P_{k}))\cap\Omega$.\\ We obtain the following results if\ $U(P_{k},\ \delta_{1}(P_{k}))\cap\partial\Omega$\ is a graph of other functions.\\ If we let\ $\delta_{0}=\delta_{1}^{\ast}$, then from\ $\delta_{1}^{\ast}\leq \delta_{1}(P_{k})/3$, we obtain the following,
$$U(P_{0},\ \delta_{0})\cap\partial\Omega\subset U(P_{k},\ \delta_{1}(P_{k}))\cap\partial\Omega. $$
Then\ $\forall P=(x,\ y,\ z)^{T}\in U(P_{0},\ \delta_{0})\cap\partial\Omega$, we obtain the following,
\[n_{p}=\cfrac{1}{\varphi(x,\ y)}\ (-f_{kx}(x,\ y),\ -f_{ky}(x,\ y),\ 1)^{T},\ n_{p_{0}}=\cfrac{1}{\varphi(x_{0},\ y_{0})}\ (-f_{kx}(x_{0},\ y_{0}),\ -f_{ky}(x_{0},\ y_{0}),\ 1)^{T},\] where\ $\varphi(x,\ y)=\sqrt{f_{kx}^{2}(x,\ y)+f_{ky}^{2}(x,\ y)+1}$,\ $f_{kx},\ f_{ky}$\ are partial derivatives of\ $f_{k}$.\\ We can deduce the following,
\begin{eqnarray*}\cfrac{f_{kx}(x,\ y)}{\varphi(x,\ y)}-\cfrac{f_{kx}(x_{0},\ y_{0})}{\varphi(x_{0},\ y_{0})}&=&\cfrac{\varphi(x_{0},\ y_{0})f_{kx}(x,\ y)-\varphi(x,\ y)f_{kx}(x_{0},\ y_{0})}{\varphi(x,\ y)\varphi(x_{0},\ y_{0})}\\&=&\cfrac{(\varphi(x_{0},\ y_{0})-\varphi(x,\ y))f_{kx}(x,\ y)+\varphi(x,\ y)(f_{kx}(x,\ y)-f_{kx}(x_{0},\ y_{0}))}{\varphi(x,\ y)\varphi(x_{0},\ y_{0})},\end{eqnarray*}
\begin{eqnarray*}\cfrac{f_{ky}(x,\ y)}{\varphi(x,\ y)}-\cfrac{f_{ky}(x_{0},\ y_{0})}{\varphi(x_{0},\ y_{0})}&=&\cfrac{\varphi(x_{0},\ y_{0})f_{ky}(x,\ y)-\varphi(x,\ y)f_{ky}(x_{0},\ y_{0})}{\varphi(x,\ y)\varphi(x_{0},\ y_{0})}\\&=&\cfrac{(\varphi(x_{0},\ y_{0})-\varphi(x,\ y))f_{ky}(x,\ y)+\varphi(x,\ y)(f_{ky}(x,\ y)-f_{ky}(x_{0},\ y_{0}))}{\varphi(x,\ y)\varphi(x_{0},\ y_{0})},\end{eqnarray*}
\begin{eqnarray*} \varphi(x_{0},\ y_{0})-\varphi(x,\ y)&=&\cfrac{\varphi^{2}(x_{0},\ y_{0})-\varphi^{2}(x,\ y)}{\varphi(x_{0},\ y_{0})+\varphi(x,\ y)}\\ \varphi^{2}(x_{0},\ y_{0})-\varphi^{2}(x,\ y)&=&f^{2}_{kx}(x_{0},\ y_{0})-f_{kx}^{2}(x,\ y)+f^{2}_{ky}(x_{0},\ y_{0})-f_{ky}^{2}(x,\ y).\end{eqnarray*}
If we assume,\[ M_{k}=\max_{P=(x,\ y,\ z)^{T}\in \overline{U}(P_{k},\ \delta_{1}(P_{k}))\cap\partial\Omega}(|f_{kx}(x,\ y)|,\ |f_{ky}(x,\ y)|),\ 1\leq k\leq N,\]then from\ $\varphi(x,\ y)\geq 1$, we can obtain,\[\mid n_{p}-n_{p_{0}}\mid\leq (4M_{k}^{2}+4M_{k}+1)C_{\alpha}(P_{k})r_{PP_{0}}^{\alpha},\ \forall P=(x,\ y,\ z)^{T}\in U(P_{0},\ \delta_{0})\cap\partial\Omega,\] where\ $C_{\alpha}(P_{k})$\ is defined in the same way as in Theorem 2.1.\\
If we denote,\[ M_{0}=\max_{1\leq k\leq N}M_{k},\ C_{\alpha}=\max_{1\leq k\leq N}C_{\alpha}(P_{k}),\ C_{0}=(4M_{0}^{2}+4M_{0}+1)C_{\alpha},\]
then we know (2.14) holds.\qed\\
\begin{theorem} \label{Theorem2-2} If\ $\Omega$\ is bounded,\ $\partial\Omega\in C^{1,\ \alpha},\ 0<\alpha\leq 1,$\ then we have the following for an ababsolute solid angle that is defined in Lyapunov's potential theory on page 182 of [6],
\[ \max_{M\in R^{3}}\int_{\partial\Omega}\cfrac{|{\bf r}_{MP}\cdot n_{p}|}{r_{MP}^{3}}dS_{P}<+\infty.\]\end{theorem}
{\it Proof of theorem 2.2}. If\ $M=(x,\ y,\ z)^{T}\in\partial\Omega$, then from lemma 2.1, we know there exists a neighborhood\ $U(M,\ \delta(M)),\ \delta(M)>0$,\ and\ $C_{\alpha}>0$, such that\ $|{\bf r}_{MP}\cdot n_{p}|\leq 2C_{\alpha} r_{MP}^{(1+\alpha)}$,\ $\forall P\in U(M,\ \delta(M))\cap \partial\Omega$, where\ $C_{\alpha}$\ is defined in the same way as in lemma 2.3. \\We obtain the following,$$\cfrac{|{\bf r}_{MP}\cdot n_{p}|}{r_{MP}^{3}}\leq \cfrac{2C_{\alpha}}{r_{MP}^{2-\alpha}},\ 0<\alpha\leq 1,\ \forall P\in U(M,\ \delta(M))\cap \partial\Omega.$$
From lemma 2.2, we can see that\[ \int_{\partial\Omega}\cfrac{|{\bf r}_{MP}\cdot n_{p}|}{r_{MP}^{3}}dS_{P},\] is continuous on\ $\partial\Omega$. Hence we get the following,\[ \max_{M\in \partial\Omega}\int_{\partial\Omega}\cfrac{|{\bf r}_{MP}\cdot n_{p}|}{r_{MP}^{3}}dS_{P}<+\infty.\]
Next we assume,
\[\Omega_{1}=\{ M\in R^{3}\setminus\overline{\Omega}:\ dist(M,\ \partial\Omega)\leq \delta_{2}^{\ast}\},\ \Omega_{2}=\{ M\in \Omega:\ dist(M,\ \partial\Omega)\leq \delta_{2}^{\ast}\},\]where$$\delta_{2}^{\ast}=[\delta_{1}^{\ast}\cos(\theta^{\ast})]/2,\ dist(M,\ \partial\Omega)=\min_{P\in\partial\Omega}|M-P|,$$ $\delta_{1}^{\ast},\ \theta^{\ast}$\ are defined in the same way as in Theorem 2.1.
\\If\ $M=(x,\ y,\ z)^{T}\in R^{3}\setminus(\Omega_{1}\cup\Omega_{2})$, then we obtain the following,
$$\int_{\partial\Omega}\cfrac{|{\bf r}_{MP}\cdot n_{p}|}{r_{MP}^{3}}dS_{P}\leq \cfrac{S(\partial\Omega)}{(\delta_{2}^{\ast})^{2}}, $$ where
$$S(\partial\Omega)=\int_{\partial\Omega}dS_{P}.$$Hence we get as follows,\[ \max_{M\in R^{3}\setminus(\Omega_{1}\cup\Omega_{2})}\int_{\partial\Omega}\cfrac{|{\bf r}_{MP}\cdot n_{p}|}{r_{MP}^{3}}dS_{P}<+\infty.\]
\\If\ $M=(x,\ y,\ z)^{T}\in \Omega_{1}$, then there exists\ $P_{0}\in\partial\Omega$, such that\ $|M-P_{0}|=\ dist(M,\ \partial\Omega)$. We discuss a smooth curve on\ $\partial\Omega$\ that passes through\ $P_{0}$. The parameter coordinates of the point\ $P$\ on the curve are\ $(x(\theta),\ y(\theta),\ z(\theta))$.The parameter coordinates of \ $P_{0}$\ are\ $(x(\theta_{0}),\ y(\theta_{0}),\ z(\theta_{0}))$. We assume the tangent vector at\ $P_{0}$\ is as follows,
$$ s_{0}=(x^{\prime}(\theta_{0}),\ y^{\prime}(\theta_{0}),\ z^{\prime}(\theta_{0}))^{T}.$$ We denote\ $f(\theta)$\ as follows,
$$f(\theta)=r_{MP}^{2}=(x-x(\theta))^{2}+(y-y(\theta))^{2}+(z-z(\theta))^{2}.$$ Then\ $f(\theta)$\ attains the minimum at\ $\theta_{0}$. Since\ $f(\theta)$\ is smooth, we see that$$ f^{\prime}(\theta)|_{\theta=\theta_{0}}=-2(M-P_{0})\cdot s_{0}=0.$$ So\ ${\bf r}_{MP_{0}}$\ is perpendicular to\ $s_{0}$. From the arbitrary of tangent vector\ $s_{0}$, we obtain that\ ${\bf r}_{MP_{0}}$\ is parallel to\ $n_{p_{0}}$. \\
From Theorem 2.1, we know that there exists a uniform exterior finite right circular cone\ $K$\ for\ $P_{0}$,\ $P_{0}$\ is the vertex of a cone\ $K(P_{0})$,\ $n_{p_{0}}$\ is the symmetry axis, the polar angle is\ $\theta^{\ast}$,\ and the length of generatrix is\ $\delta_{1}^{\ast}$. From\ $|M-P_{0}|\leq \delta_{2}^{\ast}=[\delta_{1}^{\ast}\cos(\theta^{\ast})]/2$,\ and\ $M\in\Omega_{1},\ {\bf r}_{MP_{0}}$\ is parallel to\ $n_{p_{0}}$, and we can obtain\ $M\in K(P_{0})$, and\ $M$\ is on the symmetry axis\ $n_{p_{0}}$.\\
Again from Theorem 2.1, we know there exists\ $P_{k}\in \partial\Omega,\ 1\leq k\leq N$, such that\ $P_{0}\in U(P_{k},\ \delta_{1}(P_{k})/3)\cap\partial\Omega$. And\ $U(P_{k},\ \delta_{1}(P_{k}))\cap\partial\Omega$\ is a graph of a\ $C^{1,\ \alpha},\ 0<\alpha\leq 1,$\ function of two of the coordinates\ $x,\ y,\ z$. \\From\ $\delta_{1}^{\ast}\leq \delta_{1}(P_{k})/3$, we get\ $U(P_{0},\ \delta_{1}^{\ast})\cap\partial\Omega\subset U(P_{k},\ \delta_{1}(P_{k}))\cap\partial\Omega.$\\
If we denote\ $\theta_{1}=({\bf r}_{P_{0}M},\ {\bf r}_{P_{0}P}),\ \forall P\in U(P_{0},\ \delta_{1}^{\ast})\cap\partial\Omega$, where\ $({\bf r}_{P_{0}M},\ {\bf r}_{P_{0}P})$\ is the angle between\ ${\bf r}_{P_{0}M}$\ and\ ${\bf r}_{P_{0}P}$, then we obtain\ $\theta_{1}> \theta^{\ast}$. We denote\ $\theta_{2}=({\bf r}_{P_{0}M},\ {\bf r}_{PM})$.\\ If\ $\theta_{1}<\pi/2$, then we obtain
$$r_{P_{0}P}=\cfrac{ r_{MP}\sin\theta_{2}}{\sin\theta_{1}}\leq \cfrac{1}{\sin\theta^{\ast}} r_{MP}.$$ If\ $\theta_{1}\geq\pi/2$, then we obtain
$$r_{P_{0}P}\leq  r_{MP}\leq\cfrac{1}{\sin\theta^{\ast}} r_{MP},\ \mbox{where}\ \theta^{\ast}\in(0,\ \cfrac{\pi}{2}).$$ Hence we obtain
\[r_{P_{0}P}\leq\cfrac{1}{\sin\theta^{\ast}} r_{MP},\ \forall P\in U(P_{0},\ \delta_{1}^{\ast})\cap\partial\Omega.\]
From lemma 2.3, we can obtain\ $\forall P\in U(P_{0},\ \delta_{1}^{\ast})\cap\partial\Omega,$
\begin{eqnarray*}|\cfrac{{\bf r}_{MP}\cdot n_{p}}{r_{MP}^{3}}-\cfrac{{\bf r}_{MP}\cdot n_{p_{0}}}{r_{MP}^{3}}|&=&|\cfrac{{\bf r}_{MP}\cdot (n_{p}-n_{p_{0}})}{r_{MP}^{3}}|\\&\leq& C_{0}\cfrac{r_{P_{0}P}^{\alpha}}{r_{MP}^{2}}\leq \cfrac{C_{0}}{(\sin\theta^{\ast})^{\alpha}}\cfrac{1}{r_{MP}^{2-\alpha}}.\end{eqnarray*}
And from lemma 2.1, we can get\ $\forall P\in U(P_{0},\ \delta_{1}^{\ast})\cap\partial\Omega,$\begin{eqnarray*}|\cfrac{{\bf r}_{MP}\cdot n_{p_{0}}}{r_{MP}^{3}}-\cfrac{{\bf r}_{MP_{0}}\cdot n_{p_{0}}}{r_{MP}^{3}}|&=&|\cfrac{{\bf r}_{P_{0}P}\cdot n_{p_{0}}}{r_{MP}^{3}}|\\&\leq& 2C_{\alpha}\cfrac{r_{P_{0}P}^{1+\alpha}}{r_{MP}^{3}}\leq \cfrac{2C_{\alpha}}{(\sin\theta^{\ast})^{1+\alpha}}\cfrac{1}{r_{MP}^{2-\alpha}}.\end{eqnarray*}
Hence we can see the following,\begin{eqnarray*}\int_{U(P_{0},\ \delta_{1}^{\ast})\cap\partial\Omega}\cfrac{|{\bf r}_{MP}\cdot n_{p}|}{r_{MP}^{3}}dS_{P}&\leq&\int_{U(P_{0},\ \delta_{1}^{\ast})\cap\partial\Omega}\cfrac{|{\bf r}_{MP_{0}}\cdot n_{p_{0}}|}{r_{MP}^{3}}dS_{P}+\\&&[\cfrac{C_{0}}{(\sin\theta^{\ast})^{\alpha}}+\cfrac{2C_{\alpha}}{(\sin\theta^{\ast})^{1+\alpha}}]\int_{U(P_{0},\ \delta_{1}^{\ast})\cap\partial\Omega}\cfrac{1}{r_{MP}^{2-\alpha}}dS_{P}.\end{eqnarray*}From lemma 2.2, we know$$\int_{\partial\Omega}\cfrac{1}{r_{MP}^{2-\alpha}}dS_{P}$$ is continuous on\ $M$. So we can assume,\[ C_{1}=\max_{M\in \Omega_{1}\cup\Omega_{2}\cup\partial\Omega}\int_{\partial\Omega}\cfrac{1}{r_{MP}^{2-\alpha}}dS_{P}.\]
Hence we can obtain\[\int_{U(P_{0},\ \delta_{1}^{\ast})\cap\partial\Omega}\cfrac{|{\bf r}_{MP}\cdot n_{p}|}{r_{MP}^{3}}dS_{P}\leq\int_{U(P_{0},\ \delta_{1}^{\ast})\cap\partial\Omega}\cfrac{|{\bf r}_{MP_{0}}\cdot n_{p_{0}}|}{r_{MP}^{3}}dS_{P}+[\cfrac{C_{0}}{(\sin\theta^{\ast})^{\alpha}}+\cfrac{2C_{\alpha}}{(\sin\theta^{\ast})^{1+\alpha}}]C_{1}.\]
From\ $M$\ is on\ $n_{p_{0}}$, we can get that\ $|{\bf r}_{MP_{0}}\cdot n_{p_{0}}|=r_{MP_{0}}$. From the cosine law, we can obtain
\begin{eqnarray*}r_{MP}^{2}&=&r_{MP_{0}}^{2}+r_{P_{0}P}^{2}-2r_{MP_{0}}r_{P_{0}P}\cos\theta_{1}\\
&\geq&r_{MP_{0}}^{2}+r_{P_{0}P}^{2}-2r_{MP_{0}}r_{P_{0}P}\cos\theta^{\ast}\\&\geq&(1-\cos\theta^{\ast})(r_{MP_{0}}^{2}+r_{P_{0}P}^{2}).\end{eqnarray*}
We assume\ $U(P_{j},\ \delta_{1}(P_{j}))\cap\partial\Omega$\ is a graph of a\ $C^{1,\ \alpha},\ 0<\alpha\leq 1,$\ function\ $f_{j}(x_{j1},\ x_{j2})$, where\ $x_{j1},\ x_{j2}$\ are two of the coordinates\ $x,\ y,\ z,\ 1\leq j\leq N$. We assume
$$ C_{2}=\max_{1\leq j\leq N}C_{2}(P_{j}),\ C_{2}(P_{j})=\max_{P=(x,\ y,\ z)^{T}\in \overline{U}(P_{j},\ \delta_{1}(P_{j}))}\sqrt{f^{2}_{j1}(P(x_{j1}),\ P(x_{j2}))+f^{2}_{j2}(P(x_{j1}),\ P(x_{j2}))+1},$$where\ $f_{j1},\ f_{j2}$\ are partial derivatives of\ $f_{j},\ P(x_{jl})$\ is the value of coordinate\ $x_{jl}$\ at point\ $P$,\ $l=1,\ 2,\ 1\leq j\leq N$.\\Now we can obtain the following,
\begin{eqnarray*} \int_{U(P_{0},\ \delta_{1}^{\ast})\cap\partial\Omega}\cfrac{r_{MP_{0}}}{r_{MP}^{3}}dS_{P}&\leq &\cfrac{1}{(1-\cos\theta^{\ast})^{3/2}} \int_{U(P_{0},\ \delta_{1}^{\ast})\cap\partial\Omega}\cfrac{r_{MP_{0}}}{(r_{MP_{0}}^{2}+r_{P_{0}P}^{2})^{3/2}}dS_{P} \\(a=r_{MP_{0}})
&\leq& \cfrac{2\pi C_{2}}{(1-\cos\theta^{\ast})^{3/2}}\int_{0}^{\delta_{1}^{\ast}}\cfrac{ardr}{(r^{2}+a^{2})^{3/2}}.\end{eqnarray*}
And we can work out the following,\ $\forall a>0$,
$$ \int_{0}^{\delta_{1}^{\ast}}\cfrac{ardr}{(r^{2}+a^{2})^{3/2}}= \cfrac{-a}{\sqrt{r^{2}+a^{2}}}|_{0}^{\delta_{1}^{\ast}}
= 1-\cfrac{a}{\sqrt{(\delta_{1}^{\ast})^{2}+a^{2}}}\leq1.$$
From (2.26), we can obtain,\[\int_{U(P_{0},\ \delta_{1}^{\ast})\cap\partial\Omega}\cfrac{|{\bf r}_{MP}\cdot n_{p}|}{r_{MP}^{3}}dS_{P}\leq\cfrac{2\pi C_{2}}{(1-\cos\theta^{\ast})^{3/2}}+[\cfrac{C_{0}}{(\sin\theta^{\ast})^{\alpha}}+\cfrac{2C_{\alpha}}{(\sin\theta^{\ast})^{1+\alpha}}]C_{1}.\]
And from\ $M\in\Omega_{1}$, we can get\ $r_{MP_{0}}\leq \delta_{2}^{\ast}$. Hence we can obtain$$r_{MP}\geq r_{P_{0}P}-r_{MP_{0}}\geq \delta_{1}^{\ast}-\delta_{2}^{\ast},\ \forall P\in\partial\Omega\setminus U(P_{0},\ \delta_{1}^{\ast}).$$ So we get as follows,
 \[\int_{\partial\Omega\setminus U(P_{0},\ \delta_{1}^{\ast})}\cfrac{|{\bf r}_{MP}\cdot n_{p}|}{r_{MP}^{3}}dS_{P}\leq  \cfrac{S(\partial\Omega)}{(\delta_{1}^{\ast}-\delta_{2}^{\ast})^{2}}, \]
where\ $$S(\partial\Omega)=\int_{\partial\Omega}dS_{P}.$$Hence we get the following,\ $\forall M\in\Omega_{1}$,
\[\int_{\partial\Omega}\cfrac{|{\bf r}_{MP}\cdot n_{p}|}{r_{MP}^{3}}dS_{P}\leq  \cfrac{2\pi C_{2}}{(1-\cos\theta^{\ast})^{3/2}}+[\cfrac{C_{0}}{(\sin\theta^{\ast})^{\alpha}}
+\cfrac{2C_{\alpha}}{(\sin\theta^{\ast})^{1+\alpha}}]C_{1}+\cfrac{S(\partial\Omega)}{(\delta_{1}^{\ast}-\delta_{2}^{\ast})^{2}}. \]In the same way, we can obtain following,\ $\forall M\in\Omega_{2}$,
\[\int_{\partial\Omega}\cfrac{|{\bf r}_{MP}\cdot n_{p}|}{r_{MP}^{3}}dS_{P}\leq  \cfrac{2\pi C_{2}}{(1-\cos\theta^{\ast})^{3/2}}+[\cfrac{C_{0}}{(\sin\theta^{\ast})^{\alpha}}
+\cfrac{2C_{\alpha}}{(\sin\theta^{\ast})^{1+\alpha}}]C_{1}+\cfrac{S(\partial\Omega)}{(\delta_{1}^{\ast}-\delta_{2}^{\ast})^{2}}, \]
which proves the statement.\qed \\
\begin{corollary} If\ $\Omega$\ is bounded,\ $\partial\Omega\in C^{1,\ \alpha},\ 0<\alpha\leq 1,$\ then double layer potential
\[u(M)=\int_{\partial\Omega}v(P)\cfrac{\partial\cfrac{1}{r_{PM}}}{\partial n_{p}}dS_{P},\] is continuous on\ $R^{3}\setminus \partial\Omega$, where\ $v(P)\in C(\partial\Omega)$, moreover\ $\forall P_{0}\in \partial\Omega$, we have\begin{eqnarray}&& \lim_{M\rightarrow P_{0}+}u(M)=u(P_{0})-2\pi v(P_{0}),\\&&  \lim_{M\rightarrow P_{0}-}u(M)=u(P_{0})+2\pi v(P_{0}),\end{eqnarray}where\ $M\rightarrow P_{0}+$\ means\ $M$\ is near to\ $P_{0}$\ from the interior of\ $\Omega$\ and\ $M\rightarrow P_{0}-$\ means\ $M$\ is near to\ $P_{0}$\ from the exterior of\ $\Omega$. \end{corollary}
{\it Proof of corollary 2.1}. From lemma 2.2 and the previous Theorem, we may get that\ $\forall P_{0}\in \partial\Omega$,
$$u_{0}(M)=\int_{\partial\Omega}(v(P)-v(P_{0}))\cfrac{\partial\cfrac{1}{r_{PM}}}{\partial n_{p}}dS_{P}$$ is continuous at\ $P_{0}$. From the potential theory, we know$$\int_{\partial\Omega}\cfrac{\partial\cfrac{1}{r_{PM}}}{\partial n_{p}}dS_{P}
=\begin{cases}-4\pi,\ M\in \Omega,\\-2\pi,\ M\in \partial\Omega,\\0,\ M\in R^{3}\setminus\overline{\Omega}.\end{cases}
$$Hence, the statement holds.\qed \\
\begin{theorem} \label{Theorem2-3} If\ $\Omega$\ is bounded,\ $\partial\Omega\in C^{1,\ \alpha},\ 0<\alpha\leq 1,$\ then simple layer potential
\[u(M)=\int_{\partial\Omega}\cfrac{v(P)}{r_{PM}}dS_{P},\] where\ $v(P)\in C(\partial\Omega)$, satisfies the following,\ $\forall P_{0}\in \partial\Omega$,
\begin{eqnarray}&& \cfrac{\partial u(P_{0})}{\partial n_{p_{0}}^{+}}=\int_{\partial\Omega}v(P)\cfrac{\partial\cfrac{1}{r_{PP_{0}}}}{\partial n_{p_{0}}}dS_{P}-2\pi v(P_{0}),\\&&  \cfrac{\partial u(P_{0})}{\partial n_{p_{0}}^{-}}=\int_{\partial\Omega}v(P)\cfrac{\partial\cfrac{1}{r_{PP_{0}}}}{\partial n_{p_{0}}}dS_{P}+2\pi v(P_{0}),\end{eqnarray} where\[ \cfrac{\partial u(P_{0})}{\partial n_{p_{0}}^{+}}=\lim_{M\rightarrow n_{p_{0}}^{+}}\cfrac{u(M)-u(P_{0})}{r_{MP_{0}}},\ \cfrac{\partial u(P_{0})}{\partial n_{p_{0}}^{-}}=\lim_{M\rightarrow n_{p_{0}}^{-}}\cfrac{u(P_{0})-u(M)}{r_{P_{0}M}},\]
 where\ $M\rightarrow n_{p_{0}}^{+}$\ means\ $M$\ is near to\ $P_{0}$\ along\ $n_{p_{0}}$\ from the exterior of\ $\Omega$\ and\ $M\rightarrow n_{p_{0}}^{-}$\ means\ $M$\ is near to\ $P_{0}$\ along\ $n_{p_{0}}$\ from the interior of\ $\Omega$. \end{theorem}
{\it Proof of theorem 2.3}. We refer to the proof on pages 190 to page 193 of [6]. If\ $M\in n_{p_{0}}\setminus \partial\Omega$, then a directional deriviative\ $\partial u(M)/\partial n_{p_{0}}$\ exists, and we can work it out through the integral as follows.
\[\cfrac{\partial u(M)}{\partial n_{p_{0}}}=\int_{\partial\Omega}v(P)\cfrac{\partial\cfrac{1}{r_{PM}}}{\partial n_{p_{0}}}dS_{P}
=-\int_{\partial\Omega}v(P)\cfrac{\cos({\bf r}_{PM},\ n_{p_{0}})}{r_{PM}^{2}}dS_{P},\] where\ $({\bf r}_{PM},\ n_{p_{0}})$\ is the angle between\ ${\bf r}_{PM}$\ and\ $n_{p_{0}}$. Together with double potential\[u_{1}(M)=\int_{\partial\Omega}v(P)\cfrac{\partial\cfrac{1}{r_{PM}}}{\partial n_{p}}dS_{P}=\int_{\partial\Omega}v(P)\cfrac{\cos({\bf r}_{PM},\ n_{p})}{r_{PM}^{2}}dS_{P},\] we have
\[\cfrac{\partial u(M)}{\partial n_{p_{0}}}+u_{1}(M)=\int_{\partial\Omega}v(P)\cfrac{\cos({\bf r}_{PM},\ n_{p})-\cos({\bf r}_{PM},\ n_{p_{0}})}{r_{PM}^{2}}dS_{P}.\]
We want to prove the right side of (2.40) is continuous when\ $M$\ is near to\ $P_{0}$\ along\ $n_{P_{0}}$. We only need to prove\ $\forall \epsilon>0,\ \exists\delta>0$, such that\ $\forall M\in U(P_{0},\ \delta)\cap n_{p_{0}}$,\[\mid\int_{(\partial\Omega)_{\delta}}v(P)\cfrac{\cos({\bf r}_{PM},\ n_{p})-\cos({\bf r}_{PM},\ n_{p_{0}})}{r_{PM}^{2}}dS_{P}\mid\leq\epsilon, \]
where\ $(\partial\Omega)_{\delta}=U(P_{0},\ \delta)\cap\partial\Omega$.\\ If we assume\[\max_{P\in \partial\Omega}\mid v(P)\mid=C_{3},\]
then we have\begin{eqnarray}\mid v(P)\cfrac{\cos({\bf r}_{PM},\ n_{p})-\cos({\bf r}_{PM},\ n_{p_{0}})}{r_{PM}^{2}}\mid &\leq& C_{3}\cfrac{\mid \cos({\bf r}_{PM},\ n_{p})-\cos({\bf r}_{PM},\ n_{p_{0}})\mid}{r_{PM}^{2}}\\&\leq& 2C_{3}\cfrac{\mid \sin\cfrac{({\bf r}_{PM},\ n_{p})-({\bf r}_{PM},\ n_{p_{0}})}{2}\mid}{r_{PM}^{2}}\\
&\leq&2C_{3}\cfrac{\mid \sin\cfrac{(n_{p},\ n_{p_{0}})}{2}\mid}{r_{PM}^{2}}=C_{3}\cfrac{\mid n_{p}-n_{p_{0}}\mid}{r_{PM}^{2}},\end{eqnarray}
where $$\mid \cfrac{({\bf r}_{PM},\ n_{p})-({\bf r}_{PM},\ n_{p_{0}})}{2}\mid\\
\leq\mid \cfrac{(n_{p},\ n_{p_{0}})}{2}\mid, $$ is obtained from the sum of two angles of the trihedral being no less than the third one, and\ $(n_{p},\ n_{p_{0}})/2$\ being in\ $[0,\ \pi/2]$.\\
From lemma 2.3 and (2.24), we have\ $\forall P=(x,\ y,\ z)^{T}\in U(P_{0},\ \delta_{1}^{\ast})\cap\partial\Omega,\ \forall M\in n_{p_{0}},\ |M-P_{0}|\leq \delta_{2}^{\ast}$, \[\mid v(P)\cfrac{\cos({\bf r}_{PM},\ n_{p})-\cos({\bf r}_{PM},\ n_{p_{0}})}{r_{PM}^{2}}\mid \leq \cfrac{C_{0}C_{3}}{(\sin\theta^{\ast})^{\alpha}}\cfrac{1}{r_{PM}^{2-\alpha}}.\]
From lemma 2.2, we know (2.41) is true and the right side of (2.40) is continuous when\ $M$\ is near to\ $P_{0}$\ along\ $n_{P_{0}}$. \\
Since the continuity, we have the following,\begin{eqnarray}&&\lim_{M^{\prime}\rightarrow n_{p_{0}}^{+}}(\cfrac{\partial u(M^{\prime})}{\partial n_{p_{0}}}+u_{1}(M^{\prime}))=\lim_{M^{\prime\prime}\rightarrow n_{p_{0}}^{-}}(\cfrac{\partial u(M^{\prime\prime})}{\partial n_{p_{0}}}+u_{1}(M^{\prime\prime}))\\&&=\int_{\partial\Omega}v(P)\cfrac{\cos({\bf r}_{PP_{0}},\ n_{p})-\cos({\bf r}_{PP_{0}},\ n_{p_{0}})}{r_{PP_{0}}^{2}}dS_{P}.\end{eqnarray}
From corollary 2.1, we have\begin{eqnarray}&&\lim_{M^{\prime}\rightarrow n_{p_{0}}^{+}}u_{1}(M^{\prime})=\lim_{M^{\prime}\rightarrow P_{0}-}u_{1}(M^{\prime})
=\int_{\partial\Omega}v(P)\cfrac{\partial\cfrac{1}{r_{PP_{0}}}}{\partial n_{p}}dS_{P}+2\pi v(P_{0})\\&&=\int_{\partial\Omega}v(P)\cfrac{\cos({\bf r}_{PP_{0}},\ n_{p})}{r_{PP_{0}}^{2}}dS_{P}+2\pi v(P_{0}),\end{eqnarray}\begin{eqnarray}&&\lim_{M^{\prime\prime}\rightarrow n_{p_{0}}^{+}}u_{1}(M^{\prime\prime})=\lim_{M^{\prime\prime}\rightarrow P_{0}-}u_{1}(M^{\prime\prime})
=\int_{\partial\Omega}v(P)\cfrac{\partial\cfrac{1}{r_{PP_{0}}}}{\partial n_{p}}dS_{P}-2\pi v(P_{0})\\&&=\int_{\partial\Omega}v(P)\cfrac{\cos({\bf r}_{PP_{0}},\ n_{p})}{r_{PP_{0}}^{2}}dS_{P}-2\pi v(P_{0}).\end{eqnarray}
And from (2.47) and (2.48) we can push out the limits\[\cfrac{\partial u(P_{0})}{\partial n_{p_{0}}^{+}}=\lim_{M^{\prime}\rightarrow n_{p_{0}}^{+}}\cfrac{\partial u(M^{\prime})}{\partial n_{p_{0}}},\ \cfrac{\partial u(P_{0})}{\partial n_{p_{0}}^{-}}=\lim_{M^{\prime\prime}\rightarrow n_{p_{0}}^{-}}\cfrac{\partial u(M^{\prime\prime})}{\partial n_{p_{0}}},\]and integrals\[\int_{\partial\Omega}v(P)\cfrac{-\cos({\bf r}_{PP_{0}},\ n_{p_{0}})}{r_{PP_{0}}^{2}}dS_{P}\]all exist, moreover\begin{eqnarray}&&\cfrac{\partial u(P_{0})}{\partial n_{p_{0}}^{+}}=\lim_{M^{\prime}\rightarrow n_{p_{0}}^{+}}\cfrac{\partial u(M^{\prime})}{\partial n_{p_{0}}}=\int_{\partial\Omega}v(P)\cfrac{-\cos({\bf r}_{PP_{0}},\ n_{p_{0}})}{r_{PP_{0}}^{2}}dS_{P}-2\pi v(P_{0})\\&&=\int_{\partial\Omega}v(P)\cfrac{\partial\cfrac{1}{r_{PP_{0}}}}{\partial n_{p_{0}}}dS_{P}-2\pi v(P_{0}),\\&&\cfrac{\partial u(P_{0})}{\partial n_{p_{0}}^{-}}=\lim_{M^{\prime\prime}\rightarrow n_{p_{0}}^{-}}\cfrac{\partial u(M^{\prime\prime})}{\partial n_{p_{0}}}=\int_{\partial\Omega}v(P)\cfrac{-\cos({\bf r}_{PP_{0}},\ n_{p_{0}})}{r_{PP_{0}}^{2}}dS_{P}+2\pi v(P_{0})\\&&=\int_{\partial\Omega}v(P)\cfrac{\partial\cfrac{1}{r_{PP_{0}}}}{\partial n_{p_{0}}}dS_{P}+2\pi v(P_{0}).\end{eqnarray}
That's the end of proof.\qed \\These last two theorems are the main results in Lyapunov's potential theory on pages 173 to 193 of [6]. We see that\ $\partial\Omega\in C^{1,\ \alpha},\ 0<\alpha\leq 1,$ can play the role of Lyapunov's surface.\\
Finally, we see whether\ $\partial\Omega\in C^{1,\ \alpha},\ 0<\alpha\leq 1$, is the Hopf's surface for the strong maximum principle. From the example on page 35 of [2], we know it is not the Hopf's surface if domain\ $\Omega$\ only satisfies an interior cone condition. However, we can obtain a stronger condition as follows.\begin{theorem} \label{Theorem2-4} If\ $\Omega$\ is bounded,\ $\partial\Omega\in C^{1,\ \alpha},\ 0<\alpha\leq 1,$\ then domain\ $\Omega$\ satisfies a uniform interior oblate spheroid condition, that is,\ $\exists \delta>0,\ \forall P_{0}\in \partial\Omega$, there exists a finite right oblate spheroid as follows,\[K_{\delta}(P_{0})=\{P: {\bf r}_{P_{0}P}\cdot(-n_{p_{0}})\geq 4C_{\alpha}r_{P_{0}P}^{1+\alpha},\ r_{P_{0}P}\leq \delta\},\]
where\ $C_{\alpha}$\ is defined in the same way as in lemma 2.3,\ $K_{\delta}(P_{0})\cap(R^{3}\setminus\Omega)=P_{0}$. \end{theorem}
{\it Proof of theorem 2.4}. If we select\ $\delta\leq\delta_{1}^{\ast}$, then we will get\ $K_{\delta}(P_{0})\cap(R^{3}\setminus\Omega)=P_{0}$.\\ If there exists\ $P\in K_{\delta}(P_{0})\cap(R^{3}\setminus\Omega)$, and\ $P\neq P_{0}$, then from lemma 2.1, we obtain\[ {\bf r}_{P_{0}P}\cdot(-n_{p_{0}})\leq 2C_{\alpha}r_{P_{0}P}^{1+\alpha}.\] This contradicts (2.59).\qed \\ Now we can see that domain\ $\Omega$\ satisfies not only a uniform exterior and interior cone condition, but also a uniform exterior and interior oblate spheroid condition, if\ $\Omega$\ is bounded,\ $\partial\Omega\in C^{1,\ \alpha},\ 0<\alpha\leq 1$. Therefore, we arrive at the following.\begin{theorem} \label{Theorem2-5} If\ $\Omega$\ is bounded,\ $\partial\Omega\in C^{1,\ \alpha},\ 0<\alpha\leq 1$, supposing\ $u\in C^{2}(\Omega)$\ and\ $\triangle u=0$\ in\ $\Omega$, letting\ $P_{0}\in \partial\Omega$\ be such that\ $u$\ is continuous at\ $P_{0}$,\ $u(P_{0})>u(P)$\ for all\ $P\in \Omega$, then the exterior normal derivative of\ $u$\ at\ $P_{0}$, if it exists, satisfies the strict inequality\[ \cfrac{\partial u}{\partial n}(P_{0})=\lim_{P\rightarrow n_{p_{0}}^{-}}\cfrac{u(P_{0})-u(P)}{r_{P_{0}P}}>0.\]\end{theorem}
{\it Proof of theorem 2.5}. From the previous theorem, we know\ $\forall \delta\in (0,\ \delta_{1}^{\ast}]$, there exists an oblate spheroid\ $K_{\delta}(P_{0})\subset \Omega$. We introduce an auxiliary function\ $v_{1}$\ by defining\[ v_{1}(P)=e^{-(4C_{\alpha}r_{P_{0}P}^{1+\alpha})^{\gamma}}-e^{-({\bf r}_{P_{0}P}\cdot(-n_{p_{0}}))^{\gamma}}+e^{-(4C_{\alpha}r_{P_{0}P}^{1+\alpha})}-e^{-({\bf r}_{P_{0}P}\cdot(-n_{p_{0}}))},\]
where\ $\gamma\in(1,\ 1+\alpha)$. \\ If we assume\ $P_{0}(x_{0},\ y_{0},\ z_{0}),\ P(x,\ y,\ z),\ n_{p_{0}}=(n_{01},\ n_{02},\ n_{03})^{T},\ n_{01}^{2}+n_{02}^{2}+n_{03}^{2}=1$, then direct calculation gives
\begin{eqnarray*}\cfrac{\partial v_{1}}{\partial n}(P_{0})&=&\lim_{P\rightarrow n_{p_{0}}^{-}}\cfrac{v_{1}(P_{0})-v_{1}(P)}{r_{P_{0}P}}\\
&=&\lim_{r_{P_{0}P}\rightarrow 0}\cfrac{-e^{-(4C_{\alpha}r_{P_{0}P}^{1+\alpha})^{\gamma}}+e^{-(r_{P_{0}P})^{\gamma}}-e^{-(4C_{\alpha}r_{P_{0}P}^{1+\alpha})}+e^{-r_{P_{0}P}}}{r_{P_{0}P}}\ (L'Hopital)\\
&=&\lim_{r_{P_{0}P}\rightarrow 0}e^{-(4C_{\alpha}r_{P_{0}P}^{1+\alpha})^{\gamma}}(4C_{\alpha})^{\gamma}(\gamma+\gamma\alpha)r_{P_{0}P}^{\gamma+\gamma\alpha-1}
-e^{-(r_{P_{0}P})^{\gamma}}\gamma(r_{P_{0}P})^{\gamma-1}\\&&+e^{-(4C_{\alpha}r_{P_{0}P}^{1+\alpha})}4C_{\alpha}(1+\alpha)r_{P_{0}P}^{\alpha}-e^{-r_{P_{0}P}}\\&=&-1,
\end{eqnarray*}
\begin{eqnarray*}{\bf r}_{P_{0}P}\cdot(-n_{p_{0}})&=&(x-x_{0})(-n_{01})+(y-y_{0})(-n_{02})+(z-z_{0})(-n_{03}),\\ r_{P_{0}P}&=&\sqrt{(x-x_{0})^{2}+(y-y_{0})^{2}+(z-z_{0})^{2}},\\
\cfrac{\partial r_{P_{0}P}}{\partial x}&=&\cfrac{x-x_{0}}{r_{P_{0}P}},
\end{eqnarray*}
\begin{eqnarray*}\cfrac{\partial v_{1}}{\partial x}
&=&e^{-(4C_{\alpha}r_{P_{0}P}^{1+\alpha})^{\gamma}}[-(4C_{\alpha})^{\gamma}](\gamma+\gamma\alpha)r_{P_{0}P}^{\gamma+\gamma\alpha-2}(x-x_{0})\\
&&-e^{-({\bf r}_{P_{0}P}\cdot(-n_{p_{0}}))^{\gamma}}(-\gamma)({\bf r}_{P_{0}P}\cdot(-n_{p_{0}}))^{\gamma-1}(-n_{01})\\
&&+e^{-(4C_{\alpha}r_{P_{0}P}^{1+\alpha})}[-(4C_{\alpha})](1+\alpha)r_{P_{0}P}^{1+\alpha-2}(x-x_{0})
-e^{-({\bf r}_{P_{0}P}\cdot(-n_{p_{0}}))}n_{01},
\end{eqnarray*}
\begin{eqnarray*}\cfrac{\partial^{2} v_{1}}{\partial x^{2}}
&=&e^{-(4C_{\alpha}r_{P_{0}P}^{1+\alpha})^{\gamma}}[(4C_{\alpha})^{2\gamma}(\gamma+\gamma\alpha)^{2}r_{P_{0}P}^{2\gamma+2\gamma\alpha-4}(x-x_{0})^{2}\\
&&-(4C_{\alpha})^{\gamma}(\gamma+\gamma\alpha)(r_{P_{0}P}^{\gamma+\gamma\alpha-2}+(\gamma+\gamma\alpha-2)r_{P_{0}P}^{\gamma+\gamma\alpha-4}(x-x_{0})^{2})]\\
&&-e^{-({\bf r}_{P_{0}P}\cdot(-n_{p_{0}}))^{\gamma}}[(-\gamma)^{2}({\bf r}_{P_{0}P}\cdot(-n_{p_{0}}))^{2\gamma-2}(-n_{01})^{2}-\gamma(\gamma-1)({\bf r}_{P_{0}P}\cdot(-n_{p_{0}}))^{\gamma-2}(-n_{01})^{2}]\\
&&+e^{-(4C_{\alpha}r_{P_{0}P}^{1+\alpha})}[(4C_{\alpha})^{2}(1+\alpha)^{2}r_{P_{0}P}^{2+2\alpha-4}(x-x_{0})^{2}\\&&
-(4C_{\alpha})(1+\alpha)(r_{P_{0}P}^{1+\alpha-2}+(1+\alpha-2)r_{P_{0}P}^{1+\alpha-4}(x-x_{0})^{2})]-e^{-({\bf r}_{P_{0}P}\cdot(-n_{p_{0}}))}(-1)^{2}(-n_{01})^{2},
\end{eqnarray*}
\begin{eqnarray*}\triangle v_{1}&=&e^{-(4C_{\alpha}r_{P_{0}P}^{1+\alpha})^{\gamma}}[(4C_{\alpha})^{2\gamma}(\gamma+\gamma\alpha)^{2}r_{P_{0}P}^{2\gamma+2\gamma\alpha-2}
-(4C_{\alpha})^{\gamma}(\gamma+\gamma\alpha)(\gamma+\gamma\alpha+1)r_{P_{0}P}^{\gamma+\gamma\alpha-2}]\\ &&-e^{-({\bf r}_{P_{0}P}\cdot(-n_{p_{0}}))^{\gamma}}[\gamma^{2}({\bf r}_{P_{0}P}\cdot(-n_{p_{0}}))^{2\gamma-2}-\gamma(\gamma-1)({\bf r}_{P_{0}P}\cdot(-n_{p_{0}}))^{\gamma-2}]\\&&+e^{-(4C_{\alpha}r_{P_{0}P}^{1+\alpha})}[(4C_{\alpha})^{2}(1+\alpha)^{2}r_{P_{0}P}^{2\alpha}
-(4C_{\alpha})(1+\alpha)(2+\alpha) r_{P_{0}P}^{\alpha-1}]-e^{-({\bf r}_{P_{0}P}\cdot(-n_{p_{0}}))}.\end{eqnarray*}
The main items are\ $r_{P_{0}P}^{\gamma+\gamma\alpha-2}$,\ $({\bf r}_{P_{0}P}\cdot(-n_{p_{0}}))^{\gamma-2}$, and\ $r_{P_{0}P}^{\alpha-1}$. From$$ \lim_{r_{P_{0}P}\rightarrow 0}({\bf r}_{P_{0}P}\cdot(-n_{p_{0}}))^{2-\gamma}r_{P_{0}P}^{\alpha-1}=\lim_{r_{P_{0}P}\rightarrow 0}({\bf r}_{P_{0}P}\cdot(-n_{p_{0}}))^{2-\gamma}r_{P_{0}P}^{\gamma+\gamma\alpha-2}=0,$$we obtain$$\lim_{r_{P_{0}P}\rightarrow 0}\triangle v_{1}=+\infty.$$
Hence, there exists\ $ \delta_{1}\in (0,\ \delta_{1}^{\ast}]$, such that\ $\triangle v_{1}>0$, throughout oblate spheroid\ $K_{\delta_{1}}(P_{0})$.\\
At the point\ $P$\ which\ ${\bf r}_{P_{0}P}\cdot(-n_{p_{0}})=4C_{\alpha}r_{P_{0}P}^{1+\alpha}$, we have\ $v_{1}(P)=0$. At the point\ $P$\ which\ $r_{P_{0}P}=\delta_{1}$,
$v_{1}(P)$\ is not\ $0$\ but bounded. Therefore, there exists a constant\ $\epsilon>0$, for which\ $u-u(P_{0})+\epsilon v_{1}\leq 0$, on\ $\partial K_{\delta_{1}}(P_{0})$. The maximum principle now implies that\ $u-u(P_{0})+\epsilon v_{1}\leq 0$, throughout oblate spheroid\ $K_{\delta_{1}}(P_{0})$. Taking the exterior normal derivative of at\ $P_{0}$, we obtain, as required,
\[\cfrac{\partial u}{\partial n}(P_{0})\geq -\epsilon\cfrac{\partial v_{1}}{\partial n}(P_{0})=\epsilon>0.\] That's the end of the proof.
\qed \\
\section{Equivalence} \setcounter{equation}{0}
We transform Eqs.(1.1) and (1.2) into the following,
 \begin{eqnarray}
 \cfrac{\partial u_{1}}{\partial x}&=&-\cfrac{\partial u_{2}}{\partial y}-\cfrac{\partial u_{3}}{\partial z},\\
 \cfrac{\partial u}{\partial
 t}&=& (u_{xx}+u_{yy}+u_{zz})-\tau grad p+v,\end{eqnarray}
 where $$\tau=\cfrac{1}{\rho},\ v=(\mu-1) (u_{xx}+u_{yy}+u_{zz})-u_{1}\cfrac{\partial u}{\partial
 x}-u_{2}\cfrac{\partial u}{\partial
 y}-u_{3}\cfrac{\partial u}{\partial
 z}+F,\ \cfrac{\partial u}{\partial
 x}=(-\cfrac{\partial u_{2}}{\partial y}-\cfrac{\partial u_{3}}{\partial z},\ \cfrac{\partial u_{2}}{\partial
 x},\ \cfrac{\partial u_{3}}{\partial
 x})^{T}.$$ Let's introduce\ $Z=(u,\ p,\ \partial u\setminus\ u_{1x},\ \partial^{2} u, grad p,\ v)^{T}$,\ $\partial u=(u_{x},\ u_{y},\ u_{z})^{T},\ \partial u\setminus\ u_{1x}=(u_{2x},\ u_{3x},\ u_{y},\ u_{z})^{T}$,\ $u_{jx}=\partial u_{j}/\partial x,\ j=1,\ 2,\ 3$,\ $\partial^{2}u=(u_{xx},\ u_{xy},\ u_{xz},\ u_{yy},\ u_{yz},\ u_{zz})^{T}$.\\Then Eq.(3.1) and Eqs.(3.2) is equivalent to \begin{eqnarray}
  \cfrac{\partial u}{\partial
 t}&=&\alpha_{1}^{T}Z,\\ \cfrac{\partial u_{1}}{\partial x}&=&\alpha_{2}^{T}Z,\end{eqnarray} where
                                                                                                           \begin{eqnarray*}\alpha_{1}^{T}&=&(
                                                                                                             0_{3\times3},\ 0_{3\times1},\ 0_{3\times8},\ E,\ 0_{3\times3},\ 0_{3\times3},\ E,\ 0_{3\times3},\ E,\ -\tau E,,\ E),\\ \alpha_{2}^{T}&=&(
                                                                                                             0_{1\times3},\ 0_{1\times1},\ (0,\ 0,\ 0,\ -1,\ 0,\ 0,\ 0,\ -1),\ 0_{1\times3},\ 0_{1\times3},\ 0_{1\times3},\ 0_{1\times3},\ 0_{1\times3},\ 0_{1\times3},\ 0_{1\times3},\ 0_{1\times3}),\end{eqnarray*}\ $E$\ is three order identity matrix.\\
 We should discuss\ $Z$\ as follows,
 \begin{eqnarray*}&&\cfrac{\partial u}{\partial
 t}=\alpha_{1}^{T}Z,\ \cfrac{\partial u_{1}}{\partial x}=\alpha_{2}^{T}Z,\ u=E_{1}^{T}Z,\ p=e_{4}^{T}Z,\ \cfrac{\partial u_{2}}{\partial x}=e_{5}^{T}Z,\ \cfrac{\partial u_{3}}{\partial x}=e_{6}^{T}Z,\ u_{y}=E_{3}^{T}Z,\\&& u_{z}=E_{4}^{T}Z,\ u_{xx}=E_{5}^{T}Z,\ u_{xy}=E_{6}^{T}Z,\ u_{xz}=E_{7}^{T}Z,\ u_{yy}=E_{8}^{T}Z,\ u_{yz}=E_{9}^{T}Z,\\&& u_{zz}=E_{10}^{T}Z,\ grad p=E_{11}^{T}Z,\ v=E_{12}^{T}Z,\ \mbox{where}\ E_{j}=(e_{3j-2},\ e_{3j-1},\ e_{3j}),\ 1\leq j\leq12,\\&& e_{k}\ \mbox{is the\ $k$th 36 dimensional unit coordinate vector},\ 1\leq k\leq 36.\end{eqnarray*}
We can obtain Eq.(3.1) and Eqs.(3.2) are equivalent to the following system respect to\ $Z$,\begin{eqnarray}&&\cfrac{\partial E_{1}^{T}Z}{\partial t}=\alpha_{1}^{T}Z,\\&& \cfrac{\partial E_{1}^{T}Z}{\partial x}=\left(
                                                                                                                                                      \begin{array}{c}
                                                                                                                                                        \alpha_{2}^{T}Z \\
                                                                                                                                                       e_{5}^{T}Z \\
                                                                                                                                                        e_{6}^{T}Z \\
                                                                                                                                                      \end{array}
                                                                                                                                                    \right)
 ,\\&& \cfrac{\partial E_{1}^{T}Z}{\partial y}=E_{3}^{T}Z,\\&& \cfrac{\partial E_{1}^{T}Z}{\partial z}=E_{4}^{T}Z,\\&& \cfrac{\partial^{2} E_{1}^{T}Z}{\partial x^{2}}=E_{5}^{T}Z,\\&& \cfrac{\partial^{2} E_{1}^{T}Z}{\partial x \partial y}=E_{6}^{T}Z,\\&& \cfrac{\partial^{2} E_{1}^{T}Z}{\partial x \partial z}=E_{7}^{T}Z,\\&& \cfrac{\partial^{2} E_{1}^{T}Z}{\partial y^{2}}=E_{8}^{T}Z,\end{eqnarray}\begin{eqnarray}&& \cfrac{\partial^{2} E_{1}^{T}Z}{\partial y \partial z}=E_{9}^{T}Z,\\&& \cfrac{\partial^{2} E_{1}^{T}Z}{\partial z^{2}}=E_{10}^{T}Z,\\&& grad(e_{4}^{T}Z)=E_{11}^{T}Z, \\&& (\mu-1)(E_{5}^{T}Z+E_{8}^{T}Z+E_{10}^{T}Z) -e_{1}^{T}Z\left(
                                                                                                                                                      \begin{array}{c}
                                                                                                                                                        \alpha_{2}^{T}Z \\
                                                                                                                                                       e_{5}^{T}Z \\
                                                                                                                                                        e_{6}^{T}Z \\
                                                                                                                                                      \end{array}
                                                                                                                                                    \right)
 \nonumber\\&&-e_{2}^{T}ZE_{3}^{T}Z-e_{3}^{T}ZE_{4}^{T}Z+F=E_{12}^{T}Z=v,\end{eqnarray}where\ $Z\in C(\overline{\Omega}\times [0,\ T]),\  E_{1}^{T}Z\in C^{2}(\overline{\Omega})\cap C^{1}[0,\ T],\  e_{4}^{T}Z\in C^{1}(\overline{\Omega})\cap C[0,\ T]$.\\We have the first equivalent result as follows.
\begin{theorem} \label{Theorem3-1} Eq.(3.1) and Eqs(3.2) is equivalent to the system from Eqs(3.5) to (3.16) with respect to\ $Z$. \end{theorem}{\it Proof of theorem 3.1}. If\ $u,\ p$\ satisfies Eq.(3.1) and Eqs(3.2), then letting\ $Z=T_{1}((u,\ p)^{T})=(u,\ p,\ \partial u\setminus\ u_{1x},\ \partial^{2} u, grad p,\ v)^{T}$, we obtain\begin{eqnarray*}&&\cfrac{\partial u}{\partial
 t}=\alpha_{1}^{T}Z,\ \cfrac{\partial u_{1}}{\partial x}=\alpha_{2}^{T}Z,\ u=E_{1}^{T}Z,\ p=e_{4}^{T}Z,\ \cfrac{\partial u_{2}}{\partial x}=e_{5}^{T}Z,\ \cfrac{\partial u_{3}}{\partial x}=e_{6}^{T}Z,\ u_{y}=E_{3}^{T}Z,\\&& u_{z}=E_{4}^{T}Z,\ u_{xx}=E_{5}^{T}Z,\ u_{xy}=E_{6}^{T}Z,\ u_{xz}=E_{7}^{T}Z,\ u_{yy}=E_{8}^{T}Z,\ u_{yz}=E_{9}^{T}Z,\\&& u_{zz}=E_{10}^{T}Z,\ grad p=E_{11}^{T}Z,\ v=E_{12}^{T}Z.\end{eqnarray*} Hence,\ $Z$\ satisfies Eqs(3.5) to (3.16). \\If\ $Z$\ satisfies Eqs(3.5) to (3.16), then letting\ $(u,\ p)^{T}=T_{2}(Z)=(E_{1},\ e_{4})^{T}Z$, we obtain
\begin{eqnarray*}&&\cfrac{\partial u}{\partial
 t}=\alpha_{1}^{T}Z,\ \cfrac{\partial u_{1}}{\partial x}=\alpha_{2}^{T}Z,\ u=E_{1}^{T}Z,\ p=e_{4}^{T}Z,\ \cfrac{\partial u_{2}}{\partial x}=e_{5}^{T}Z,\ \cfrac{\partial u_{3}}{\partial x}=e_{6}^{T}Z,\ u_{y}=E_{3}^{T}Z,\\&& u_{z}=E_{4}^{T}Z,\ u_{xx}=E_{5}^{T}Z,\ u_{xy}=E_{6}^{T}Z,\ u_{xz}=E_{7}^{T}Z,\ u_{yy}=E_{8}^{T}Z,\ u_{yz}=E_{9}^{T}Z,\\&& u_{zz}=E_{10}^{T}Z,\ grad p=E_{11}^{T}Z,\ v=E_{12}^{T}Z.\end{eqnarray*}It follows that\ $u,\ p$\ satisfies Eq.(3.1) and Eqs(3.2).\\ Obviously\ $T_{1},\ T_{2}$\ are continuous. Moreover\ $T_{1}(T_{2}(Z))=Z,\ T_{2}(T_{1}((u,\ p)^{T}))=(u,\ p)^{T}.$ From definition 1.1, we know the statement stands.\qed\\We notice that Eqs(3.5) to (3.15) are good, because they are the second order linear partial differential equations with constant coefficients. Eq.(3.16) could be considered as complicated, but if we assume
$$ Z=\left(
       \begin{array}{c}
         Z_{1} \\
         Z_{2} \\
       \end{array}
     \right),\ \mbox{where}\ Z_{1}\ \mbox{is the first\ $33$\ componenets of}\ Z,$$ then we will obtain\ $Z_{2}=\psi(Z_{1})$\ from Eq.(3.16).  Next we want to get\ $Z_{1}=T_{0}(Z_{2})$\ from Eq.(3.5) to (3.15). This is not difficult. From our experience, we guess that\ $T_{0}$\ should be related to the integral equations. We will obtain\ $T_{0}$\ by the Fourier transform on\ $\overline{\Omega}\times [0,\ T]$. At last we will transform Eq.(3.5) to (3.16) into the equivalent generalized integral equations\ $Z_{1}=T_{0}(\psi(Z_{1}))$.\\ We apply the Fourier transform on\ $\overline{\Omega}\times [0,\ T]$\ on both sides from Eq.(3.5) to (3.15) as follows, \begin{eqnarray*}&&\int_{0}^{T}dt\int_{\overline{\Omega}}e^{-it\xi_{0}-ix\xi_{1}-iy\xi_{2}-iz\xi_{3}}\cfrac{\partial E_{1}^{T}Z}{\partial t}dxdydz=\int_{0}^{T}dt\int_{\overline{\Omega}}e^{-it\xi_{0}-ix\xi_{1}-iy\xi_{2}-iz\xi_{3}}\alpha_{1}^{T}Zdxdydz,\\&& \int_{0}^{T}dt\int_{\overline{\Omega}}e^{-it\xi_{0}-ix\xi_{1}-iy\xi_{2}-iz\xi_{3}} \cfrac{\partial E_{1}^{T}Z}{\partial x}dxdydz=\int_{0}^{T}dt\int_{\overline{\Omega}}e^{-it\xi_{0}-ix\xi_{1}-iy\xi_{2}-iz\xi_{3}}\left(
                                                                                                                                                      \begin{array}{c}
                                                                                                                                                        \alpha_{2}^{T}Z \\
                                                                                                                                                       e_{5}^{T}Z \\
                                                                                                                                                        e_{6}^{T}Z \\
                                                                                                                                                      \end{array}
                                                                                                                                                    \right)dxdydz,
                                                                                                                                                    \\&& \int_{0}^{T}dt\int_{\overline{\Omega}}e^{-it\xi_{0}-ix\xi_{1}-iy\xi_{2}-iz\xi_{3}}\cfrac{\partial E_{1}^{T}Z}{\partial y}dxdydz=\int_{0}^{T}dt\int_{\overline{\Omega}}e^{-it\xi_{0}-ix\xi_{1}-iy\xi_{2}-iz\xi_{3}}E_{3}^{T}Zdxdydz,\\&&\int_{0}^{T}dt\int_{\overline{\Omega}}e^{-it\xi_{0}-ix\xi_{1}-iy\xi_{2}-iz\xi_{3}} \cfrac{\partial E_{1}^{T}Z}{\partial z}dxdydz=\int_{0}^{T}dt\int_{\overline{\Omega}}e^{-it\xi_{0}-ix\xi_{1}-iy\xi_{2}-iz\xi_{3}}E_{4}^{T}Zdxdydz,\\&&\int_{0}^{T}dt\int_{\overline{\Omega}}e^{-it\xi_{0}-ix\xi_{1}-iy\xi_{2}-iz\xi_{3}}\cfrac{\partial^{2} E_{1}^{T}Z}{\partial x^{2}}dxdydz=\int_{0}^{T}dt\int_{\overline{\Omega}}e^{-it\xi_{0}-ix\xi_{1}-iy\xi_{2}-iz\xi_{3}}E_{5}^{T}Zdxdydz,\\&& \int_{0}^{T}dt\int_{\overline{\Omega}}e^{-it\xi_{0}-ix\xi_{1}-iy\xi_{2}-iz\xi_{3}} \cfrac{\partial^{2} E_{1}^{T}Z}{\partial x\partial y}dxdydz=\int_{0}^{T}dt\int_{\overline{\Omega}}e^{-it\xi_{0}-ix\xi_{1}-iy\xi_{2}-iz\xi_{3}}E_{6}^{T}Zdxdydz,\\&& \int_{0}^{T}dt\int_{\overline{\Omega}}e^{-it\xi_{0}-ix\xi_{1}-iy\xi_{2}-iz\xi_{3}}\cfrac{\partial^{2} E_{1}^{T}Z}{\partial x\partial z}dxdydz=\int_{0}^{T}dt\int_{\overline{\Omega}}e^{-it\xi_{0}-ix\xi_{1}-iy\xi_{2}-iz\xi_{3}}E_{7}^{T}Zdxdydz,\\&&\int_{0}^{T}dt\int_{\overline{\Omega}}e^{-it\xi_{0}-ix\xi_{1}-iy\xi_{2}-iz\xi_{3}} \cfrac{\partial^{2} E_{1}^{T}Z}{\partial y^{2}}dxdydz=\int_{0}^{T}dt\int_{\overline{\Omega}}e^{-it\xi_{0}-ix\xi_{1}-iy\xi_{2}-iz\xi_{3}}E_{8}^{T}Zdxdydz,\\&& \int_{0}^{T}dt\int_{\overline{\Omega}}e^{-it\xi_{0}-ix\xi_{1}-iy\xi_{2}-iz\xi_{3}}\cfrac{\partial^{2} E_{1}^{T}Z}{\partial y\partial z}dxdydz=\int_{0}^{T}dt\int_{\overline{\Omega}}e^{-it\xi_{0}-ix\xi_{1}-iy\xi_{2}-iz\xi_{3}}E_{9}^{T}Zdxdydz,\\&&\int_{0}^{T}dt\int_{\overline{\Omega}}e^{-it\xi_{0}-ix\xi_{1}-iy\xi_{2}-iz\xi_{3}} \cfrac{\partial^{2} E_{1}^{T}Z}{\partial z^{2}}dxdydz=\int_{0}^{T}dt\int_{\overline{\Omega}}e^{-it\xi_{0}-ix\xi_{1}-iy\xi_{2}-iz\xi_{3}}E_{10}^{T}Zdxdydz,\\&&\int_{0}^{T}dt\int_{\overline{\Omega}}e^{-it\xi_{0}-ix\xi_{1}-iy\xi_{2}-iz\xi_{3}} grad (e_{4}^{T}Z)dxdydz=\int_{0}^{T}dt\int_{\overline{\Omega}}e^{-it\xi_{0}-ix\xi_{1}-iy\xi_{2}-iz\xi_{3}}E_{11}^{T}Zdxdydz.\end{eqnarray*}
The reason why we apply the Fourier transform on\ $\overline{\Omega}\times[0,\ T]$\ instead of\ $R^{4}$\ is that Eqs(1.1) and (1.2) only holds on\ $\overline{\Omega}\times[0,\ T]$. It is possible that they will not stand outside of\ $\overline{\Omega}\times[0,\ T]$. We can't apply the Fourier transform on\ $R^{4}$\ to both sides of Eq(3.5) to (3.15). In order to denote this easily, we define the Fourier transform on\ $\overline{\Omega}\times[0,\ T]$\ as
follows.\begin{definition}\label{definition} $\forall\ f(x,\ y,\ z,\ t)\in L^{2}(\overline{\Omega}\times[0,\ T]),\ \forall (\xi_{1},\ \xi_{2},\ \xi_{3},\ \xi_{0})^{T}\in R^{4}$,
\begin{eqnarray*}FI(f(x,\ y,\ z,\ t))&=&\int_{0}^{T}dt\int_{\overline{\Omega}}f(x,\ y,\ z,\ t)e^{-it\xi_{0}-ix\xi_{1}-iy\xi_{2}-iz\xi_{3}}dxdydz\\&=&F(f(x,\ y,\ z,\ t)I_{\overline{\Omega}\times[0,\ T]}(x,\ y,\ z,\ t)),\end{eqnarray*}
where\ $F$\ means the Fourier transform and\ $I_{\overline{\Omega}\times[0,\ T]}(x,\ y,\ z,\ t)$\
is the characteristic function. In the following, we write\ $I_{\overline{\Omega}\times[0,\ T]}(x,\ y,\ z,\ t)$\ into\ $I_{\overline{\Omega}\times[0,\ T]}$.\end{definition}
By using divergence theorem, we obtain\begin{eqnarray*}
FI(\cfrac{\partial E_{1}^{T}Z}{\partial t})&=&\int_{0}^{T}dt\int_{\overline{\Omega}}e^{-it\xi_{0}-ix\xi_{1}-iy\xi_{2}-iz\xi_{3}}(\cfrac{\partial E_{1}^{T}Z}{\partial t})dxdydz
\\&=&\int_{\overline{\Omega}}(E_{1}^{T}Z)e^{-it\xi_{0}}|_{t=0}^{t=T}e^{-ix\xi_{1}-iy\xi_{2}-iz\xi_{3}}dxdydz
+\\&&i\xi_{0}\int_{0}^{T}dt\int_{\overline{\Omega}}e^{-it\xi_{0}-ix\xi_{1}-iy\xi_{2}-iz\xi_{3}}(E_{1}^{T}Z)dxdydz
\\&=&f_{0}+ i\xi_{0} FI(E_{1}^{T}Z),\end{eqnarray*}\begin{eqnarray*}
FI(\cfrac{\partial E_{1}^{T}Z}{\partial x})&=&\int_{0}^{T}dt\int_{\overline{\Omega}}e^{-it\xi_{0}-ix\xi_{1}-iy\xi_{2}-iz\xi_{3}}(\cfrac{\partial E_{1}^{T}Z}{\partial x})dxdydz
\\&=&\int_{0}^{T}\int_{\partial\Omega}(E_{1}^{T}Z)n_{1}e^{-it\xi_{0}-ix\xi_{1}-iy\xi_{2}-iz\xi_{3}}dS
+\\&&i\xi_{1}\int_{0}^{T}dt\int_{\overline{\Omega}}e^{-it\xi_{0}-ix\xi_{1}-iy\xi_{2}-iz\xi_{3}}(E_{1}^{T}Z)dxdydz
\\&=&f_{1}+ i\xi_{1} FI(E_{1}^{T}Z),\end{eqnarray*} \begin{eqnarray*}
FI(\cfrac{\partial E_{1}^{T}Z}{\partial y})&=&\int_{0}^{T}dt\int_{\overline{\Omega}}e^{-it\xi_{0}-ix\xi_{1}-iy\xi_{2}-iz\xi_{3}}(\cfrac{\partial E_{1}^{T}Z}{\partial y})dxdydz
\\&=&\int_{0}^{T}\int_{\partial\Omega}(E_{1}^{T}Z)n_{2}e^{-it\xi_{0}-ix\xi_{1}-iy\xi_{2}-iz\xi_{3}}dS
+\\&&i\xi_{2}\int_{0}^{T}dt\int_{\overline{\Omega}}e^{-it\xi_{0}-ix\xi_{1}-iy\xi_{2}-iz\xi_{3}}(E_{1}^{T}Z)dxdydz
\\&=&f_{2}+ i\xi_{2} FI(E_{1}^{T}Z),\end{eqnarray*} \begin{eqnarray*}
FI(\cfrac{\partial E_{1}^{T}Z}{\partial z})&=&\int_{0}^{T}dt\int_{\overline{\Omega}}e^{-it\xi_{0}-ix\xi_{1}-iy\xi_{2}-iz\xi_{3}}(\cfrac{\partial E_{1}^{T}Z}{\partial z})dxdydz
\\&=&\int_{0}^{T}\int_{\partial\Omega}(E_{1}^{T}Z)n_{3}e^{-it\xi_{0}-ix\xi_{1}-iy\xi_{2}-iz\xi_{3}}dS
+\\&&i\xi_{3}\int_{0}^{T}dt\int_{\overline{\Omega}}e^{-it\xi_{0}-ix\xi_{1}-iy\xi_{2}-iz\xi_{3}}(E_{1}^{T}Z)dxdydz
\\&=&f_{3}+ i\xi_{3} FI(E_{1}^{T}Z),\end{eqnarray*}
\begin{eqnarray*}
FI(\cfrac{\partial^{2} E_{1}^{T}Z}{\partial x^{2}})&=&\int_{0}^{T}dt\int_{\overline{\Omega}}e^{-it\xi_{0}-ix\xi_{1}-iy\xi_{2}-iz\xi_{3}}(\cfrac{\partial^{2} E_{1}^{T}Z}{\partial x^{2}})dxdydz
\\&=&\int_{0}^{T}\int_{\partial\Omega}(\cfrac{\partial E_{1}^{T}Z}{\partial x})n_{1}e^{-it\xi_{0}-ix\xi_{1}-iy\xi_{2}-iz\xi_{3}}dS
+\\&&i\xi_{1}\int_{0}^{T}dt\int_{\overline{\Omega}}e^{-it\xi_{0}-ix\xi_{1}-iy\xi_{2}-iz\xi_{3}}(\cfrac{\partial E_{1}^{T}Z}{\partial x})dxdydz
\\&=&f_{11}+ i\xi_{1} FI(\cfrac{\partial E_{1}^{T}Z}{\partial x})=f_{11}+ i\xi_{1}(f_{1}+ i\xi_{1} FI(E_{1}^{T}Z)),\end{eqnarray*}
\begin{eqnarray*}
FI(\cfrac{\partial^{2} E_{1}^{T}Z}{\partial x\partial y})&=&\int_{0}^{T}dt\int_{\overline{\Omega}}e^{-it\xi_{0}-ix\xi_{1}-iy\xi_{2}-iz\xi_{3}}(\cfrac{\partial^{2} E_{1}^{T}Z}{\partial x\partial y})dxdydz
\\&=&\int_{0}^{T}\int_{\partial\Omega}(\cfrac{\partial E_{1}^{T}Z}{\partial y})n_{1}e^{-it\xi_{0}-ix\xi_{1}-iy\xi_{2}-iz\xi_{3}}dS
+\\&&i\xi_{1}\int_{0}^{T}dt\int_{\overline{\Omega}}e^{-it\xi_{0}-ix\xi_{1}-iy\xi_{2}-iz\xi_{3}}(\cfrac{\partial E_{1}^{T}Z}{\partial y})dxdydz
\\&=&f_{21}+ i\xi_{1} FI(\cfrac{\partial E_{1}^{T}Z}{\partial y})=f_{21}+ i\xi_{1}(f_{2}+ i\xi_{2} FI(E_{1}^{T}Z)),\end{eqnarray*}
\begin{eqnarray*}
FI(\cfrac{\partial^{2} E_{1}^{T}Z}{\partial x\partial z})&=&\int_{0}^{T}dt\int_{\overline{\Omega}}e^{-it\xi_{0}-ix\xi_{1}-iy\xi_{2}-iz\xi_{3}}(\cfrac{\partial^{2} E_{1}^{T}Z}{\partial x\partial z})dxdydz
\\&=&\int_{0}^{T}\int_{\partial\Omega}(\cfrac{\partial E_{1}^{T}Z}{\partial z})n_{1}e^{-it\xi_{0}-ix\xi_{1}-iy\xi_{2}-iz\xi_{3}}dS
+\\&&i\xi_{1}\int_{0}^{T}dt\int_{\overline{\Omega}}e^{-it\xi_{0}-ix\xi_{1}-iy\xi_{2}-iz\xi_{3}}(\cfrac{\partial E_{1}^{T}Z}{\partial z})dxdydz
\\&=&f_{31}+ i\xi_{1} FI(\cfrac{\partial E_{1}^{T}Z}{\partial z})=f_{31}+ i\xi_{1}(f_{3}+ i\xi_{3} FI(E_{1}^{T}Z)),\end{eqnarray*}
\begin{eqnarray*}
FI(\cfrac{\partial^{2} E_{1}^{T}Z}{\partial y^{2}})&=&\int_{0}^{T}dt\int_{\overline{\Omega}}e^{-it\xi_{0}-ix\xi_{1}-iy\xi_{2}-iz\xi_{3}}(\cfrac{\partial^{2} E_{1}^{T}Z}{\partial y^{2}})dxdydz
\\&=&\int_{0}^{T}\int_{\partial\Omega}(\cfrac{\partial E_{1}^{T}Z}{\partial y})n_{2}e^{-it\xi_{0}-ix\xi_{1}-iy\xi_{2}-iz\xi_{3}}dS
+\\&&i\xi_{2}\int_{0}^{T}dt\int_{\overline{\Omega}}e^{-it\xi_{0}-ix\xi_{1}-iy\xi_{2}-iz\xi_{3}}(\cfrac{\partial E_{1}^{T}Z}{\partial y})dxdydz
\\&=&f_{22}+ i\xi_{2} FI(\cfrac{\partial E_{1}^{T}Z}{\partial y})=f_{22}+ i\xi_{2}(f_{2}+ i\xi_{2} FI(E_{1}^{T}Z)),\end{eqnarray*}
\begin{eqnarray*}
FI(\cfrac{\partial^{2} E_{1}^{T}Z}{\partial y\partial z})&=&\int_{0}^{T}dt\int_{\overline{\Omega}}e^{-it\xi_{0}-ix\xi_{1}-iy\xi_{2}-iz\xi_{3}}(\cfrac{\partial^{2} E_{1}^{T}Z}{\partial y\partial z})dxdydz
\\&=&\int_{0}^{T}\int_{\partial\Omega}(\cfrac{\partial E_{1}^{T}Z}{\partial z})n_{2}e^{-it\xi_{0}-ix\xi_{1}-iy\xi_{2}-iz\xi_{3}}dS
+\\&&i\xi_{2}\int_{0}^{T}dt\int_{\overline{\Omega}}e^{-it\xi_{0}-ix\xi_{1}-iy\xi_{2}-iz\xi_{3}}(\cfrac{\partial E_{1}^{T}Z}{\partial z})dxdydz
\\&=&f_{32}+ i\xi_{2} FI(\cfrac{\partial E_{1}^{T}Z}{\partial z})=f_{32}+ i\xi_{2}(f_{3}+ i\xi_{3} FI(E_{1}^{T}Z)),\end{eqnarray*}
\begin{eqnarray*}
FI(\cfrac{\partial^{2} E_{1}^{T}Z}{\partial z^{2}})&=&\int_{0}^{T}dt\int_{\overline{\Omega}}e^{-it\xi_{0}-ix\xi_{1}-iy\xi_{2}-iz\xi_{3}}(\cfrac{\partial^{2} E_{1}^{T}Z}{\partial z^{2}})dxdydz
\\&=&\int_{0}^{T}\int_{\partial\Omega}(\cfrac{\partial E_{1}^{T}Z}{\partial z})n_{3}e^{-it\xi_{0}-ix\xi_{1}-iy\xi_{2}-iz\xi_{3}}dS
+\\&&i\xi_{3}\int_{0}^{T}dt\int_{\overline{\Omega}}e^{-it\xi_{0}-ix\xi_{1}-iy\xi_{2}-iz\xi_{3}}(\cfrac{\partial E_{1}^{T}Z}{\partial z})dxdydz
\\&=&f_{33}+ i\xi_{3} FI(\cfrac{\partial E_{1}^{T}Z}{\partial z})=f_{33}+ i\xi_{3}(f_{3}+ i\xi_{3} FI(E_{1}^{T}Z)),\end{eqnarray*}
\begin{eqnarray*}
FI(\cfrac{\partial e_{4}^{T}Z}{\partial x})&=&\int_{0}^{T}dt\int_{\overline{\Omega}}e^{-it\xi_{0}-ix\xi_{1}-iy\xi_{2}-iz\xi_{3}}(\cfrac{\partial e_{4}^{T}Z}{\partial x})dxdydz
\\&=&\int_{0}^{T}\int_{\partial\Omega}(e_{4}^{T}Z)n_{1}e^{-it\xi_{0}-ix\xi_{1}-iy\xi_{2}-iz\xi_{3}}dS
+\\&&i\xi_{1}\int_{0}^{T}dt\int_{\overline{\Omega}}e^{-it\xi_{0}-ix\xi_{1}-iy\xi_{2}-iz\xi_{3}}(e_{4}^{T}Z)dxdydz
\\&=&g_{1}+ i\xi_{1} FI(e_{4}^{T}Z),\end{eqnarray*} \begin{eqnarray*}
FI(\cfrac{\partial e_{4}^{T}Z}{\partial y})&=&\int_{0}^{T}dt\int_{\overline{\Omega}}e^{-it\xi_{0}-ix\xi_{1}-iy\xi_{2}-iz\xi_{3}}(\cfrac{\partial e_{4}^{T}Z}{\partial y})dxdydz
\\&=&\int_{0}^{T}\int_{\partial\Omega}(e_{4}^{T}Z)n_{2}e^{-it\xi_{0}-ix\xi_{1}-iy\xi_{2}-iz\xi_{3}}dS
+\\&&i\xi_{2}\int_{0}^{T}dt\int_{\overline{\Omega}}e^{-it\xi_{0}-ix\xi_{1}-iy\xi_{2}-iz\xi_{3}}(e_{4}^{T}Z)dxdydz
\\&=&g_{2}+ i\xi_{2} FI(e_{4}^{T}Z),\end{eqnarray*} \begin{eqnarray*}
FI(\cfrac{\partial e_{4}^{T}Z}{\partial z})&=&\int_{0}^{T}dt\int_{\overline{\Omega}}e^{-it\xi_{0}-ix\xi_{1}-iy\xi_{2}-iz\xi_{3}}(\cfrac{\partial e_{4}^{T}Z}{\partial z})dxdydz
\\&=&\int_{0}^{T}\int_{\partial\Omega}(e_{4}^{T}Z)n_{3}e^{-it\xi_{0}-ix\xi_{1}-iy\xi_{2}-iz\xi_{3}}dS
+\\&&i\xi_{3}\int_{0}^{T}dt\int_{\overline{\Omega}}e^{-it\xi_{0}-ix\xi_{1}-iy\xi_{2}-iz\xi_{3}}(e_{4}^{T}Z)dxdydz
\\&=&g_{3}+ i\xi_{3} FI(e_{4}^{T}Z),\end{eqnarray*}
where\begin{eqnarray*}f_{0}&=&\int_{\overline{\Omega}}(A_{2}e^{-iT\xi_{0}}-A_{1})e^{-ix\xi_{1}-iy\xi_{2}-iz\xi_{3}}dxdydz,\\
f_{1}&=&\int_{0}^{T}dt\int_{\partial \Omega}A_{3}n_{1}e^{-it\xi_{0}-ix\xi_{1}-iy\xi_{2}-iz\xi_{3}}dS,\\
f_{2}&=&\int_{0}^{T}dt\int_{\partial \Omega}A_{3}n_{2}e^{-it\xi_{0}-ix\xi_{1}-iy\xi_{2}-iz\xi_{3}}dS,\\
f_{3}&=&\int_{0}^{T}dt\int_{\partial \Omega}A_{3}n_{3}e^{-it\xi_{0}-ix\xi_{1}-iy\xi_{2}-iz\xi_{3}}dS,\end{eqnarray*}\begin{eqnarray*}
f_{11}&=&\int_{0}^{T}\int_{\partial\Omega}A_{4}n_{1}e^{-it\xi_{0}-ix\xi_{1}-iy\xi_{2}-iz\xi_{3}}dS,\\
f_{21}&=&\int_{0}^{T}\int_{\partial\Omega}A_{5}n_{1}e^{-it\xi_{0}-ix\xi_{1}-iy\xi_{2}-iz\xi_{3}}dS,\\
f_{31}&=&\int_{0}^{T}\int_{\partial\Omega}A_{6}n_{1}e^{-it\xi_{0}-ix\xi_{1}-iy\xi_{2}-iz\xi_{3}}dS,\\
f_{22}&=&\int_{0}^{T}\int_{\partial\Omega}A_{5}n_{2}e^{-it\xi_{0}-ix\xi_{1}-iy\xi_{2}-iz\xi_{3}}dS,\\
f_{32}&=&\int_{0}^{T}\int_{\partial\Omega}A_{6}n_{2}e^{-it\xi_{0}-ix\xi_{1}-iy\xi_{2}-iz\xi_{3}}dS,\\
f_{33}&=&\int_{0}^{T}\int_{\partial\Omega}A_{6}n_{3}e^{-it\xi_{0}-ix\xi_{1}-iy\xi_{2}-iz\xi_{3}}dS,\\
g_{1}&=&\int_{0}^{T}dt\int_{\partial \Omega}A_{9}n_{1}e^{-it\xi_{0}-ix\xi_{1}-iy\xi_{2}-iz\xi_{3}}dS,\\
g_{2}&=&\int_{0}^{T}dt\int_{\partial \Omega}A_{9}n_{2}e^{-it\xi_{0}-ix\xi_{1}-iy\xi_{2}-iz\xi_{3}}dS,\\
g_{3}&=&\int_{0}^{T}dt\int_{\partial \Omega}A_{9}n_{3}e^{-it\xi_{0}-ix\xi_{1}-iy\xi_{2}-iz\xi_{3}}dS,\end{eqnarray*}\begin{eqnarray*}
A_{1}&=&u|_{t=0},\ A_{2}=u|_{t=T},\ A_{3}=u|_{\partial\Omega\times(0,\ T)},\
A_{4}=u_{x}|_{\partial\Omega\times(0,\ T)},
\ A_{5}=u_{y}|_{\partial\Omega\times(0,\ T)},\\ A_{6}&=&u_{z}|_{\partial\Omega\times(0,\ T)},\ A_{7}=\cfrac{\partial u}{\partial n}|_{\partial\Omega\times[0,\ T]},\ A_{8}=(\cfrac{\partial u}{\partial n}+\sigma u)|_{\partial\Omega\times[0,\ T]},\ A_{9}=p|_{\partial\Omega\times(0,\ T)},
\end{eqnarray*}\ $n_{k}$\ is the\ $k$th component of
the normal vector to\ $\partial \Omega,\ k=1,\ 2,\ 3$.\\Now we have transformed the equations (3.5) to (3.15) into the following.\[\label
{A-Problem}BFI(Z)=\beta_{1},\]where$$B=\left(
                                         \begin{array}{c}
                                        i\xi_{0}E_{1}^{T}-\alpha_{1}^{T} \\
                                         i\xi_{1}E_{1}^{T}-F_{0} \\
                                         i\xi_{2}E_{1}^{T}-E_{3}^{T} \\
                                         i\xi_{3}E_{1}^{T}-E_{4}^{T} \\
                                         (i\xi_{1})^{2}E_{1}^{T}-E_{5}^{T} \\
                                         i\xi_{1}i\xi_{2}E_{1}^{T}-E_{6}^{T} \\
                                         i\xi_{1}i\xi_{3}E_{1}^{T}-E_{7}^{T} \\
                                         (i\xi_{2})^{2}E_{1}^{T}-E_{8}^{T} \\
                                         i\xi_{2}i\xi_{3}E_{1}^{T}-E_{9}^{T} \\
                                         (i\xi_{3})^{2}E_{1}^{T}-E_{10}^{T} \\
                                          i\xi_{1}e_{4}^{T}-e_{31}^{T} \\
                                         i\xi_{2}e_{4}^{T}-e_{32}^{T} \\
                                         i\xi_{3}e_{4}^{T}-e_{33}^{T} \\
                                         \end{array}
                                       \right)_{33\times36}
=(B_{1},\ -B_{2}),\ B_{2}=\left(
                                        \begin{array}{c}
                                         E \\
                                          0_{30\times3} \\
                                        \end{array}
                                      \right)_{33\times 3},\ F_{0}=\left(
                                                                                                                                                      \begin{array}{c}
                                                                                                                                                        \alpha_{2}^{T} \\
                                                                                                                                                       e_{5}^{T} \\
                                                                                                                                                        e_{6}^{T} \\
                                                                                                                                                      \end{array}
                                                                                                                                                    \right),$$
$$B_{1}=\left(
          \begin{array}{ccccccccccc}
            i\xi_{0}E, & 0, &0, &0, & -E, & 0, & 0, & -E, & 0, &  -E,&\tau E\\
            i\xi_{1}E, & F_{1}, &F_{2}, & F_{3}, & 0, & 0, & 0, & 0, & 0, & 0,&0 \\
            i\xi_{2}E, & 0, &-E, & 0, & 0, & 0, & 0, & 0, & 0, & 0,&0 \\
            i\xi_{3}E, & 0, &0, & -E, & 0, & 0, & 0, & 0, & 0, & 0,&0 \\
            (i\xi_{1})^{2}E, &  0, & 0, & 0, & -E, & 0, & 0, & 0, & 0, & 0,&0 \\
            i\xi_{1}i\xi_{2}E, &  0, & 0, & 0, & 0, &-E, & 0, & 0, & 0, & 0,&0 \\
            i\xi_{1}i\xi_{3}E, &  0, & 0, & 0, & 0, &0, & -E, & 0, & 0, & 0,&0 \\
            (i\xi_{2})^{2}E, & 0, &  0, & 0, & 0, & 0, & 0, &-E, & 0, & 0,&0 \\
           i\xi_{2}i\xi_{3}E, & 0, &  0, & 0, & 0, & 0, & 0, &0, & -E, & 0,&0 \\
            (i\xi_{3})^{2}E, & 0, & 0, & 0, & 0, & 0, & 0, &0, & 0, & -E,&0 \\
            0, & F_{4}, & 0, & 0, & 0, & 0, & 0, &0, & 0, & 0,&-E \\
          \end{array}
        \right)_{33\times33},$$
        $$F_{1}=\left(
                  \begin{array}{ccc}
                    0 & 0 & 0 \\
                    0 & -1 & 0 \\
                    0 & 0 & -1 \\
                  \end{array}
                \right),\ F_{2}=\left(
                  \begin{array}{ccc}
                    0 & 1 & 0 \\
                    0 & 0 & 0 \\
                    0 & 0 & 0 \\
                  \end{array}
                \right),\ F_{3}=\left(
                  \begin{array}{ccc}
                    0 & 0 & 1 \\
                    0 & 0 & 0 \\
                    0 & 0 & 0 \\
                  \end{array}
                \right),\ F_{4}=\left(
                  \begin{array}{ccc}
                    i\xi_{1} & 0 & 0 \\
                    i\xi_{2} & 0 & 0 \\
                    i\xi_{3} & 0 & 0 \\
                  \end{array}
                \right),$$\begin{eqnarray*}
\beta_{1}&=&(-f_{0}^{T},\ -f_{1}^{T},\ -f_{2}^{T},\ -f_{3}^{T},\ -f_{11}^{T}-i\xi_{1}f_{1}^{T},\ -f_{21}^{T}-i\xi_{1}f_{2}^{T},\ -f_{31}^{T}-i\xi_{1}f_{3}^{T},\\&& -f_{22}^{T}-i\xi_{2}f_{2}^{T},\ -f_{32}^{T}-i\xi_{2}f_{3}^{T},\ -f_{33}^{T}-i\xi_{3}f_{3}^{T},\ -g_{1},\ -g_{2},\ -g_{3})^{T}.\end{eqnarray*}
We assume
$$ Z=\left(
       \begin{array}{c}
         Z_{1} \\
         Z_{2} \\
       \end{array}
     \right),\ \mbox{where}\ Z_{1}\ \mbox{is the first\ $33$\ componenets of}\ Z,$$ then we can get$$B_{1}FI(Z_{1})=\beta_{1}+B_{2}FI(Z_{2}).$$By the primary row block transformations\ $R_{1}-R_{5}-R_{8}-R_{10}+\tau R_{11}$, we obtain $$ det(B_{1})=\tau(i\xi_{0}-((i\xi_{1})^{2}+(i\xi_{2})^{2}+(i\xi_{3})^{2}))^{2}((i\xi_{1})^{2}+(i\xi_{2})^{2}+(i\xi_{3})^{2})=a_{01}.$$ Also by the primary row block transformations, we can work out
     $$B_{1}^{-1}=-\left(
          \begin{array}{c}
           R_{1}\\
            R_{2}-i\xi_{1}F_{1}R_{1}\\
             E_{3}^{T}+i\xi_{2}R_{1} \\
             E_{4}^{T}+i\xi_{3}R_{1}\\
             E_{5}^{T}+(i\xi_{1})^{2}R_{1} \\
             E_{6}^{T}+i\xi_{1}i\xi_{2}R_{1} \\
             E_{7}^{T}+i\xi_{1}i\xi_{3}R_{1} \\
             E_{8}^{T}+(i\xi_{2})^{2}R_{1} \\
            E_{9}^{T}+i\xi_{2}i\xi_{3}R_{1} \\
             E_{10}^{T}+(i\xi_{3})^{2}R_{1}\\
              E_{11}^{T}+F_{4}R_{2}\\
          \end{array}
        \right)_{33\times33},\ B_{1}^{-1}B_{2}=\left(
          \begin{array}{c}
           B_{01}\\
            \left(
              \begin{array}{c}
                B_{03} \\
                0_{1\times3} \\
                 0_{1\times3} \\
              \end{array}
            \right)-
            i\xi_{1}F_{1}B_{01}\\
             i\xi_{2}B_{01} \\
             i\xi_{3}B_{01}\\
            (i\xi_{1})^{2}B_{01} \\
             i\xi_{1}i\xi_{2}B_{01}\\
            i\xi_{1}i\xi_{3}B_{01} \\
             (i\xi_{2})^{2}B_{01} \\
            i\xi_{2}i\xi_{3}B_{01} \\
            (i\xi_{3})^{2}B_{01}\\
            \left(
              \begin{array}{c}
                i\xi_{1}B_{03} \\
                i\xi_{2}B_{03} \\
                i\xi_{3}B_{03} \\
              \end{array}
            \right)
          \end{array}
        \right)_{33\times3},$$where\ $R_{1}=(-B_{01},\ -B_{02}e_{1}^{T},\ -B_{02}e_{2}^{T},\ -B_{02}e_{3}^{T},\ B_{01},\ 0_{3\times3},\ 0_{3\times3},\ B_{01},\ 0_{3\times3},\ B_{01},\ -\tau B_{01})_{3\times 33},$ \begin{eqnarray*}R_{2}=\left(
                                                                                                       \begin{array}{ccccccccccc}
                                                                                                         -B_{03} & -B_{04}e_{1}^{T} & -B_{04}e_{2}^{T} & -B_{04}e_{3}^{T} & B_{03} & 0_{1\times3} & 0_{1\times3} &  B_{03} & 0_{1\times3} &  B_{03} & -\tau B_{03} \\
                                                                                                         0_{1\times3} & e_{2}^{T} & 0_{1\times3} & 0_{1\times3} & 0_{1\times3} & 0_{1\times3} & 0_{1\times3} & 0_{1\times3} & 0_{1\times3} & 0_{1\times3} & 0_{1\times3} \\
                                                                                                         0_{1\times3} & e_{3}^{T} & 0_{1\times3} & 0_{1\times3} & 0_{1\times3} & 0_{1\times3} & 0_{1\times3} & 0_{1\times3} & 0_{1\times3} & 0_{1\times3} & 0_{1\times3} \\
                                                                                                       \end{array}
                                                                                                     \right)_{3\times 33},&&
        \end{eqnarray*}\ $e_{1},\ e_{2},\ e_{3}$\ in\ $R_{1},\ R_{2}$\ are all three dimensional unit coordinate vectors,
        \begin{eqnarray*}&&B_{01}=\left(
                                    \begin{array}{ccc}
                                      \cfrac{1}{a}-\cfrac{(i\xi_{1})^{2}}{a^{2}a_{1}} & -\cfrac{i\xi_{1}i\xi_{2}}{a^{2}a_{1}} & -\cfrac{i\xi_{1}i\xi_{3}}{a^{2}a_{1}} \\
                                      -\cfrac{i\xi_{1}i\xi_{2}}{a^{2}a_{1}} & \cfrac{1}{a}-\cfrac{(i\xi_{2})^{2}}{a^{2}a_{1}} & -\cfrac{i\xi_{2}i\xi_{3}}{a^{2}a_{1}} \\
                                      -\cfrac{i\xi_{1}i\xi_{3}}{a^{2}a_{1}} & -\cfrac{i\xi_{2}i\xi_{3}}{a^{2}a_{1}} & \cfrac{1}{a}-\cfrac{(i\xi_{3})^{2}}{a^{2}a_{1}} \\
                                    \end{array}
                                  \right),\ B_{02}=\left(
                                                     \begin{array}{c}
                                                       \cfrac{i\xi_{1}}{aa_{1}} \\
                                                       \cfrac{i\xi_{2}}{aa_{1}} \\
                                                       \cfrac{i\xi_{3}}{aa_{1}} \\
                                                     \end{array}
                                                   \right),\end{eqnarray*}\begin{eqnarray*}&& B_{03}=(\cfrac{i\xi_{1}}{\tau aa_{1}},\ \cfrac{i\xi_{2}}{\tau aa_{1}},\ \cfrac{i\xi_{3}}{\tau aa_{1}}),\ B_{04}=-\cfrac{1}{\tau a_{1}},\ \tau a^{3}a_{1}=a_{01}, \\&&a=i\xi_{0}-((i\xi_{1})^{2}+(i\xi_{2})^{2}+(i\xi_{3})^{2}),\ a_{1}=\cfrac{(i\xi_{1})^{2}+(i\xi_{2})^{2}+(i\xi_{3})^{2}}{a},\\&&\left(
                                                                                                                       \begin{array}{cc}
                                                                                                                         B_{01} & B_{02} \\
                                                                                                                         B_{03} & B_{04} \\
                                                                                                                       \end{array}
                                                                                                                     \right)=B_{0}^{-1},\ B_{0}=\left(
                                                                                                                                                  \begin{array}{cccc}
                                                                                                                                                    a & 0 & 0 & \tau i\xi_{1} \\
                                                                                                                                                    0 & a & 0 & \tau i\xi_{2} \\
                                                                                                                                                   0 & 0 & a & \tau i\xi_{3} \\
                                                                                                                                                    i\xi_{1} & i\xi_{2} & i\xi_{3} & 0 \\
                                                                                                                                                  \end{array}
                                                                                                                                                \right).
        \end{eqnarray*}
        The correct of\ $B_{1}^{-1}$\ is very important for the discussion as follows. We have tested it by the products of the block matrices. The following is useful, assuming that\ $F_{4,\ 1}$\ is the first column of\ $F_{4}$,$$ aB_{01}=E-\tau F_{4,\ 1}B_{03},\ aB_{02}=-\tau F_{4,\ 1}B_{04},\ (i\xi_{1},\ i\xi_{2},\ i\xi_{3})B_{01}=0,\ (i\xi_{1},\ i\xi_{2},\ i\xi_{3})B_{02}=1,\ F_{4}F_{1}=0.$$
        If we assume\ $C=\{\xi|a_{01}=0\}$, where\ $\xi=(\xi_{1},\ \xi_{2},\ \xi_{3},\ \xi_{0})^{T}$, then the measure of\ $C$\ is\ $0$. And we obtain
 \[FI(Z_{1})(1-I_{C}(\xi))=B_{1}^{-1}\beta_{1}(1-I_{C}(\xi))+B_{1}^{-1}B_{2}FI(Z_{2})(1-I_{C}(\xi)).\]
We need some lemmas in [4].
\begin{lemma} \label{lemma1}(Plancherel Theorem) If\ $f(x,\ y,\ z,\ t)\in L^{2}(R^{4})$, then\ $F(f(x,\ y,\ z,\ t))$\ exists, moreover\\(1)$\parallel F(f(x,\ y,\ z,\ t))\parallel_{L^{2}}=\parallel f(x,\ y,\ z,\ t)\parallel_{L^{2}}$,\\(2)$F^{-1}[F(f(x,\ y,\ z,\ t))=f(x,\ y,\ z,\ t)$.\end{lemma}
\begin{lemma} \label{lemma1} If\ $f(x,\ y,\ z,\ t)\in L^{2}(R^{4})$,\ $C\subset R^{4}$, the measure of\ $C$\ is\ $0$, then $$F^{-1}([F(f(x,\ y,\ z,\ t))(1-I_{C}(\xi)))=f(x,\ y,\ z,\ t).$$\end{lemma}
{\it Proof of lemma 3.2}. From the lemma 3.1, we know\ $F(f(x,\ y,\ z,\ t))\in L^{2}(R^{4})$. Therefore,
$$\int_{C}F(f(x,\ y,\ z,\ t))e^{it\xi_{0}+ix\xi_{1}+iy\xi_{2}+iz\xi_{3}}d\xi_{0}d\xi_{1}d\xi_{2}d\xi_{3}=0.$$And we obtain
$$ F^{-1}([F(f(x,\ y,\ z,\ t))](1-I_{C}(\xi)))=F^{-1}[F(f(x,\ y,\ z,\ t))]=f(x,\ y,\ z,\ t).$$\\
From these two lemmas, we obtain$$F^{-1}[FI(Z_{1})(1-I_{C}(\xi))]=Z_{1}I_{\Omega\times(0,\ T)}.$$
Now we know\ $B_{1}^{-1}\beta_{1}(1-I_{C}(\xi))+B_{1}^{-1}B_{2}FI(Z_{2})(1-I_{C}(\xi))$\ is the Fourier transform of continuous functions on\ $\overline{\Omega}\times[0,\ T]$. We have the Paley-Wiener-Schwartz theorem from [1] as follows,
\begin{lemma} \label{lemma1}(Paley-Wiener-Schwartz) Let\ $K$\ be a convex compact set of\ $R^{n}$\ with support function\ $H(\xi)=\sup_{x\in K}<x,\ \xi>,\ \forall \xi\in R^{n}$. If\ $u$\ is a distribution with support contained in\ $K$,\ then there exists\ $C>0,\ N$\ is a positive whole number, such that\[|F(u)(\zeta)|\leq C(1+|\zeta|)^{N}e^{H(Im \zeta)}, \forall\ \zeta\in C^{n}.\]
Conversely, every entire analytic function in\ $C^{n}$\ satisfying an estimate of the form (3.19) is the Fourier-Laplace transform of a distribution with support contained in\ $K$.\end{lemma}We may read the proof on pages 181 to 182 in [1]. H$\ddot{o}$rmander theorem from [1] is as follows,
\begin{lemma} \label{lemma1}(H$\ddot{o}$rmander) If\ $F(u)(\zeta)$\ is an entire analytic function in\ $C^{n}$\ satisfying an estimate of the form (3.19),\ $p(\zeta)$\ is a polynomial,\ $F(u)(\zeta)/p(\zeta)$\ is an entire function, then\ $F(u)(\zeta)/p(\zeta)$\ satisfies an estimate of the form (3.19), too.\end{lemma}We can see the proof on page 183 of [1]. \\
Because\ $a^{2}a_{1}B_{1}^{-1}$\ is entire, from lemma 3.2, Paley-Wiener-Schwartz theorem and H$\ddot{o}$rmander theorem in [1], we know that \[B_{1}^{-1}\beta_{1}+B_{1}^{-1}B_{2}FI(Z_{2})\] should be entire, too. \\
Under this condition, we can obtain$$Z_{1}I_{\overline{\Omega}\times[0,\ T]}=F^{-1}(B_{1}^{-1}\beta_{1}+B_{1}^{-1}B_{2}FI(Z_{2})).$$If we assume\ $b_{1}=a^{2}a_{1}$, then we can work out the following, \begin{eqnarray*}B_{1}^{-1}\beta_{1}&=&(b_{1}^{-1})(b_{1}B_{1}^{-1}\beta_{1}),\\
B_{1}^{-1}B_{2}&=&(b_{1}^{-1})(b_{1}B_{1}^{-1}B_{2}),
\\F^{-1}(b_{1}^{-1})&=&F^{-1}[(i\xi_{0}-((i\xi_{1})^{2}+(i\xi_{2})^{2}+(i\xi_{3})^{2}))^{-1}((i\xi_{1})^{2}+(i\xi_{2})^{2}+(i\xi_{3})^{2})^{-1}]
\\&=&\cfrac{-I_{\{t>0\}}}{4\pi(2\sqrt{\pi t})^{3}}\int_{R^{3}}e^{-\cfrac{(x-x_{1})^{2}+(y-y_{1})^{2}+(z-z_{1})^{2}}{4t}}\cfrac{1}{\sqrt{x_{1}^{2}+y_{1}^{2}+z_{1}^{2}}}dx_{1}dy_{1}dz_{1}.\end{eqnarray*}
We find that\ $F^{-1}(b_{1}^{-1})\in S^{\prime},\ b_{1}B_{1}^{-1}\beta_{1}\ \mbox{and}\ b_{1}B_{1}^{-1}B_{2}$\ satisfy the Paley-Wiener-Schwartz form (3.19). Hence\ $F^{-1}(b_{1}B_{1}^{-1}\beta_{1})\in \varepsilon^{\prime},\ F^{-1}(b_{1}B_{1}^{-1}B_{2})\in \varepsilon^{\prime}$, and their compact support is contained in\ $\overline{\Omega}\times[0,\ T]$, where\ $S^{\prime}$\ is the dual space of the Schwartz space, and\ $\varepsilon^{\prime}$\ is the dual space of\ $C^{\infty}(R^{4})$. We need a lemma as follows.
\begin{lemma} \label{lemma1} If\ $v_{1}\in S^{\prime},\ v_{2}\in \varepsilon^{\prime}$, it follows that\ $v_{1}\ast v_{2}\in S^{\prime}$\ and that
$$F(v_{1}\ast v_{2})=F(v_{1})F(v_{2}).$$\end{lemma} We can see the proof on page 166 in [1]. \\From this lemma we know\begin{eqnarray*}F^{-1}(B_{1}^{-1}\beta_{1})&=&F^{-1}(b_{1}^{-1}).\ast F^{-1}(b_{1}B_{1}^{-1}\beta_{1}),\\ F^{-1}(B_{1}^{-1}B_{2})&=&F^{-1}(b_{1}^{-1}).\ast F^{-1}(b_{1}B_{1}^{-1}B_{2}),\end{eqnarray*} all exist and are in\ $S^{\prime}$, where\ $.\ast$\ is the matrix convolution.\\
If we assume\begin{eqnarray*}w_{1}(x,\ y,\ z,\ t)&=&F^{-1}(B_{1}^{-1}\beta_{1}),\\ w_{2}(x,\ y,\ z,\ t)&=&F^{-1}(B_{1}^{-1}B_{2}),\end{eqnarray*}
then we obtain\[Z_{1}I_{\overline{\Omega}\times[0,\ T]}=w_{1}(x,\ y,\ z,\ t)+w_{2}(x,\ y,\ z,\ t).\ast Z_{2}I_{\overline{\Omega}\times[0,\ T]},\]
where\begin{eqnarray*}Z_{1}I_{\overline{\Omega}\times[0,\ T]}&=&(E_{j}^{T}ZI_{\overline{\Omega}\times[0,\ T]},\ 1\leq j\leq 11,)^{T},\\Z_{2}I_{\overline{\Omega}\times[0,\ T]}&=&E_{12}^{T}ZI_{\overline{\Omega}\times[0,\ T]}=(\mu-1)(E_{5}^{T}ZI_{\overline{\Omega}\times[0,\ T]}+E_{8}^{T}ZI_{\overline{\Omega}\times[0,\ T]}+E_{10}^{T}ZI_{\overline{\Omega}\times[0,\ T]})\\&& -e_{1}^{T}ZI_{\overline{\Omega}\times[0,\ T]}\left(
                                                                                                                                                      \begin{array}{c}
                                                                                                                                                        \alpha_{2}^{T}ZI_{\overline{\Omega}\times[0,\ T]} \\
                                                                                                                                                       e_{5}^{T}ZI_{\overline{\Omega}\times[0,\ T]} \\
                                                                                                                                                        e_{6}^{T}ZI_{\overline{\Omega}\times[0,\ T]} \\
                                                                                                                                                      \end{array}
                                                                                                                                                    \right)
-e_{2}^{T}ZI_{\overline{\Omega}\times[0,\ T]}E_{3}^{T}ZI_{\overline{\Omega}\times[0,\ T]}\\&&-e_{3}^{T}ZI_{\overline{\Omega}\times[0,\ T]}E_{4}^{T}ZI_{\overline{\Omega}\times[0,\ T]}+FI_{\overline{\Omega}\times[0,\ T]}.\end{eqnarray*}It is obvious\ $\exists\ \psi$, such that\ $Z_{2}I_{\overline{\Omega}\times[0,\ T]}=\psi(Z_{1}I_{\overline{\Omega}\times[0,\ T]})$. Therefore, we attain
\[ Z_{1}I_{\overline{\Omega}\times[0,\ T]}=w_{1}(x,\ y,\ z,\ t)+w_{2}(x,\ y,\ z,\ t).\ast(\psi(Z_{1}I_{\overline{\Omega}\times[0,\ T]})).\]
Now we arrive at the second equivalent result as follows,\begin{theorem} \label{Theorem3-2}\ $ w_1,\ w_2,\ \psi,$\ as we described, then Eq.(3.1) and Eqs(3.2) is equivalent to Eq.(3.22). \end{theorem}
 {\it Proof of theorem 3.2}. If\ $u$\ satisfies Eq.(3.1) and Eqs(3.2), then from the theorem 3.1,\\ $Z=(u,\ p,\ \partial u\setminus\ u_{1x},\ \partial^{2} u, grad p,\ v)^{T}$\ satisfies Eq.(3.5) to Eq.(3.16). Hence we can obtain the following by Fourier transform on\ $\overline{\Omega}\times[0,\ T]$,
 $$B FI(Z)=\beta_{1},\ B_{1}FI(Z_{1})=\beta_{1}+B_{2}FI(Z_{2}),\ FI(Z_{1})=B_{1}^{-1}\beta_{1}+B_{1}^{-1}B_{2}FI(Z_{2}).$$ After we do the inverse Fourier transform, we obtain\ $Z_{1}=S_{1}((u,\ p)^{T})=(u,\ p,\ \partial u\setminus\ u_{1x},\ \partial^{2} u, grad p)^{T}$\ satisfies Eq.(3.22). \\If\ $Z_{1}$\ satisfies Eq.(3.22), then letting\ $Z_{2}=\psi(Z_{1}),\ Z=(Z_{1},\ Z_{2})^{T}$,\ we obtain the following by the Fourier transform,
 $$FI(Z_{1})=B_{1}^{-1}\beta_{1}+B_{1}^{-1}B_{2}FI(Z_{2}),\ B_{1}FI(Z_{1})=\beta_{1}+B_{2}FI(Z_{2}),\ B FI(Z)=\beta_{1}.$$ After we do the inverse Fourier transform, we obtain\ $Z$\ satisfies Eq.(3.5) to Eq.(3.16) on\ $\overline{\Omega}\times[0,\ T]$. From the theorem 3.1, we know\ $(u,\ p)^{T}=S_{2}(Z_{1})=E_{1}^{T}(Z_{1},\ \psi(Z_{1}))^{T}$\ is the solution of Eq.(3.1) and Eqs(3.2).\\ Obviously\ $S_{1},\ S_{2}$\ are continuous. Moreover\ $S_{1}(S_{2}( Z_{1}))=Z_{1},\ S_{2}(S_{1}((u,\ p)^{T}))=(u,\ p)^{T}.$\ From definition 1.1, we know the the statement stands.\qed\\
We denote Eqs(3.22) as\ $Z_{1}=T_{0}(Z_{1})$, where\ $T_{0}(Z_{1})=w_{1}+w_{2}.\ast\psi(Z_{1})$. We can get a necessary and sufficient condition for there exist\ $u\in C^{2}(\overline{\Omega})\cap C^{1}[0,\ T]$,\ $p\in C^{1}(\overline{\Omega})\cap C[0,\ T]$\ satisfy Eqs(3.1) and (3.2) is that there exists\ $Z_{1}\in C(\overline{\Omega}\times[0,\ T])$\ satisfies\ $Z_{1}=T_{0}(Z_{1})$.\\Next we determine all the initial conditions and boundary conditions.\\ If we put\ $-\tau grad p$\ into\ $v$, then we will get\ $det(B_{1})=0$. So we can't do that. However, we may do that to determinate all the initial conditions and boundary conditions.\\ From Eqs(3.2), we can obtain that there exists continuous function\ $\psi_{1}(Z_{1})=-\tau E_{11}^{T}(Z_{1},\ \psi(Z_{1}))^{T}+\psi(Z_{1})$, such that
\[ u_{t}-\triangle u=\psi_{1}(Z_{1}).\]After we do Fourier transform on\ $\overline{\Omega}\times[0,\ T]$, we can get the following,
\[(i\xi_{0}-(i\xi_{1})^{2}-(i\xi_{2})^{2}-(i\xi_{3})^{2})FI(u)=-f_{0}+f_{11}+f_{22}+f_{33}+i\xi_{1}f_{1}+i\xi_{2}f_{2}+i\xi_{3}f_{3}+FI(\psi_{1}(Z_{1})).\]
We get that\[(i\xi_{0}-(i\xi_{1})^{2}-(i\xi_{2})^{2}-(i\xi_{3})^{2})^{-1}(-f_{0}+f_{11}+f_{22}+f_{33}+i\xi_{1}f_{1}+i\xi_{2}f_{2}+i\xi_{3}f_{3}+FI(\psi_{1}(Z_{1})))\]should be entire. Under this condition, we obtain that
\begin{eqnarray*}uI_{\overline{\Omega}\times[0,\ T]}&=&F^{-1}[(-a)^{-1}(f_{0}-(f_{11}+f_{22}+f_{33})\\&&-(i\xi_{1}f_{1}+i\xi_{2}f_{2}+i\xi_{3}f_{3})-FI(\psi_{1}(Z_{1})(x,\ y,\ z,\ t)))]\ \mbox{a.e.}\\F^{-1}[a^{-1}]&=&h(x,\ y,\ z,\ t)=\cfrac{1}{(2\sqrt{\pi t})^{3}}\ e^{-\cfrac{x^{2}+y^{2}+z^{2}}{4t}}I_{\{t>0\}},\end{eqnarray*}\begin{eqnarray*}
        F^{-1}((-a)^{-1}f_{0})&=&(\cfrac{1}{2\pi})^{4}\int_{R^{4}}e^{i\xi_{0}t+i\xi_{1}x+i\xi_{2}y+i\xi_{3}z}(-a)^{-1}
        (\int_{\overline{\Omega}}(A_{2}(x_{1},\ y_{1},\ z_{1})e^{-iT\xi_{0}}-A_{2}(x_{1},\ y_{1},\ z_{1}))\\&&e^{-ix_{1}\xi_{1}-iy_{1}\xi_{2}-iz_{1}\xi_{3}}dx_{1}dy_{1}dz_{1})d\xi_{0}d\xi_{1}d\xi_{2}d\xi_{3}\\&=&\int_{\overline{\Omega}}(-h(x-x_{1},\ y-y_{1},\ z-z_{1},\ t-T)A_{2}(x_{1},\ y_{1},\ z_{1})+\\&&h(x-x_{1},\ y-y_{1},\ z-z_{1},\ t)A_{2}(x_{1},\ y_{1},\ z_{1}))dx_{1}dy_{1}dz_{1},
        \end{eqnarray*}where\ $a.e.$\ means almost everywhere, $$(\cfrac{1}{2\pi})^{4}\int_{R^{4}}e^{i\xi_{0}(t-T)+i\xi_{1}(x-x_{1})+i\xi_{2}(y-y_{1})+i\xi_{3}(z-z_{1})}
        (-a)^{-1}d\xi_{0}d\xi_{1}d\xi_{2}d\xi_{3}=-h(x-x_{1},\ y-y_{1},\ z-z_{1},\ t-T),$$\begin{eqnarray*}(\cfrac{1}{2\pi})^{4}\int_{R^{4}}e^{i\xi_{0}t+i\xi_{1}(x-x_{1})+i\xi_{2}(y-y_{1})+i\xi_{3}(z-z_{1})}
        (-a)^{-1}d\xi_{0}d\xi_{1}d\xi_{2}d\xi_{3}&=&-h(x-x_{1},\ y-y_{1},\ z-z_{1},\ t),\end{eqnarray*}
         \begin{eqnarray*}F^{-1}[a^{-1}(f_{11}+f_{22}+f_{33})]&=&\int_{0}^{T}d\tau\int_{\partial\Omega}h(x-x_{1},\ y-y_{1},\ z-z_{1},\ t-\tau)A_{7}(x_{1},\ y_{1},\ z_{1},\ \tau)dS,\end{eqnarray*}
         $$F^{-1}[a^{-1}(i\xi_{1}f_{1}+i\xi_{2}f_{2}+i\xi_{3}f_{3})]=\int_{0}^{T}d\tau\int_{\partial\Omega}\cfrac{\partial h(x-x_{1},\ y-y_{1},\ z-z_{1},\ t-\tau)}{\partial n(x_{1},\ y_{1},\ z_{1})}A_{3}(x_{1},\ y_{1},\ z_{1},\ \tau)dS,$$
         $$F^{-1}[a^{-1}FI(\psi_{1}(Z_{1})(x,\ y,\ z,\ t))]=\int_{0}^{T}d\tau\int_{\overline{\Omega}} h(x-x_{1},\ y-y_{1},\ z-z_{1},\ t-\tau)\psi_{1}(Z_{1})(x_{1},\ y_{1},\ z_{1},\ \tau)dx_{1}dy_{1}dz_{1},$$
         and\begin{eqnarray*} \cfrac{\partial h(x-x_{1},\ y-y_{1},\ z-z_{1},\ t-\tau)}{\partial n(x_{1},\ y_{1},\ z_{1})}&=&\cfrac{\partial h(x-x_{1},\ y-y_{1},\ z-z_{1},\ t-\tau)}{\partial x}n_{1}(x_{1},\ y_{1},\ z_{1})+\\&&\cfrac{\partial h(x-x_{1},\ y-y_{1},\ z-z_{1},\ t-\tau)}{\partial y} n_{2}(x_{1},\ y_{1},\ z_{1})+\\&&\cfrac{\partial h(x-x_{1},\ y-y_{1},\ z-z_{1},\ t-\tau)}{\partial z}n_{3}(x_{1},\ y_{1},\ z_{1}).\end{eqnarray*}
         We denote it in an easy way as follows,\begin{eqnarray*} F^{-1}((-a)^{-1}f_{0})&=&-h(t-T).\ast_{\overline{\Omega}}A_{2}+h.\ast_{\overline{\Omega}}A_{1},
         \\F^{-1}[a^{-1}(f_{11}+f_{22}+f_{33})]&=&h.\ast_{\partial\Omega}A_{7},
         \\F^{-1}[a^{-1}(i\xi_{1}f_{1}+i\xi_{2}f_{2}+i\xi_{3}f_{3})]&=&\cfrac{\partial h}{\partial n_{p_{1}}}.\ast_{\partial\Omega}A_{3},\\F^{-1}[a^{-1}FI(\psi_{1}(Z_{1})(x,\ y,\ z,\ t))]&=&h.\ast \psi_{1}(Z_{1}),\end{eqnarray*}where\ $P_{1}=(x_{1},\ y_{1},\ z_{1})^{T}$, then we get the following,\[uI_{\overline{\Omega}\times[0,\ T]}=-h(t-T).\ast_{\overline{\Omega}}A_{2}+h.\ast_{\overline{\Omega}}A_{1}+h.\ast_{\partial\Omega}A_{7}+\cfrac{\partial h}{\partial n_{p_{1}}}.\ast_{\partial\Omega}A_{3}+h.\ast \psi_{1}(Z_{1}).\ \mbox{a.e.}\]From$$h(x,\ y,\ z,\ t)=\cfrac{1}{(2\sqrt{\pi t})^{3}}\ e^{-\cfrac{x^{2}+y^{2}+z^{2}}{4t}}I_{\{t>0\}},$$ we find\ $h(t-T).\ast_{\overline{\Omega}}A_{2}=0$, if\ $t\in[0,\ T]$. So we can obtain
         \[ uI_{\overline{\Omega}\times[0,\ T]}=h.\ast_{\overline{\Omega}}A_{1}+h.\ast_{\partial\Omega}A_{7}+\cfrac{\partial h}{\partial n_{p_{1}}}.\ast_{\partial\Omega}A_{3}+h.\ast \psi_{1}(Z_{1}),\ \mbox{a.e.}\]It is easy to think that all items in the right hand side of (3.27) are continuous on\ $\overline{\Omega}\times[0,\ T]$, and that we don't need\ $a.e.$. This is true except\ $\partial h/\partial  n_{p_{1}}.\ast_{\partial\Omega}A_{3}$. We can elaborate on this.\begin{eqnarray*}\cfrac{\partial h(M-P_{1},\ t-\tau)}{\partial n_{p_{1}}} &=&\nabla_{M}h(M-P_{1},\ t-\tau)\cdot n_{p_{1}}\\
         &=&-\nabla_{P_{1}}h(M-P_{1},\ t-\tau)\cdot n_{p_{1}},\end{eqnarray*}
         where\ $M=(x,\ y,\ z)^{T}$, and\begin{eqnarray*}\nabla_{M}h(M-P_{1},\ t-\tau)&=&(\cfrac{\partial h(M-P_{1},\ t-\tau)}{\partial x},\ \cfrac{\partial h(M-P_{1},\ t-\tau)}{\partial y},\ \cfrac{\partial h(M-P_{1},\ t-\tau)}{\partial z})^{T},\\ \nabla_{P_{1}}h(M-P_{1},\ t-\tau)&=&(\cfrac{\partial h(M-P_{1},\ t-\tau)}{\partial x_{1}},\ \cfrac{\partial h(M-P_{1},\ t-\tau)}{\partial y_{1}},\ \cfrac{\partial h(M-P_{1},\ t-\tau)}{\partial z_{1}})^{T}.\end{eqnarray*}We see\ $\nabla_{P_{1}}h(M-P_{1},\ t-\tau)\cdot n_{p_{1}}$\ in some books. That will increase a sign of minus.
         \begin{theorem} \label{Theorem3-3} If\ $\Omega$\ is bounded,\ $\partial\Omega\in C^{1,\ \beta},\ 0<\beta\leq 1,\ A_{3}(P_{1},\ t)\in C^{2}(\partial\Omega)\cap C^{1}[0,\ T]$, then\ $\partial h/\partial n_{p_{1}}.\ast_{\partial\Omega}A_{3}$\ is continuous on\ $(R^{3}\setminus \partial\Omega)\times[0,\ T]$\ and\ $\partial\Omega\times[0,\ T]$, but\ $\forall P_{0}\in \partial\Omega,\ \forall t\in(0,\ T]$, we have\begin{eqnarray}&& \lim_{M\rightarrow P_{0}+}(\cfrac{\partial h}{\partial n_{p_{1}}}.\ast_{\partial\Omega}A_{3})(M,\ t)=(\cfrac{\partial h}{\partial n_{p_{1}}}.\ast_{\partial\Omega}A_{3})(P_{0},\ t)+\cfrac{1}{2} A_{3}(P_{0},\ t),\\&&  \lim_{M\rightarrow P_{0}-}(\cfrac{\partial h}{\partial n_{p_{1}}}.\ast_{\partial\Omega}A_{3})(M,\ t)=(\cfrac{\partial h}{\partial n_{p_{1}}}.\ast_{\partial\Omega}A_{3})(P_{0},\ t)-\cfrac{1}{2} A_{3}(P_{0},\ t),\end{eqnarray}where$$(\cfrac{\partial h}{\partial n_{p_{1}}}.\ast_{\partial\Omega}A_{3})(M,\ t)=\int_{0}^{t}d\tau\int_{\partial\Omega}[\nabla_{M}h(M-P_{1},\ t-\tau)\cdot n_{p_{1}}]A_{3}(P_{1},\ \tau) dS_{P_{1}}, $$
         $M\rightarrow P_{0}+$\ means\ $M$\ is near to\ $P_{0}$\ from the interior of\ $\Omega$\ and\ $M\rightarrow P_{0}-$\ means\ $M$\ is near to\ $P_{0}$\ from the exterior of\ $\Omega$. \end{theorem}
{\it Proof of theorem 3.3}. It is obvious that\ $\partial h/\partial n_{p_{1}}.\ast_{\partial\Omega}A_{3}$\ is continuous on\ $(R^{3}\setminus \partial\Omega)\times[0,\ T]$. If\ $M=(x,\ y,\ z)^{T}\in \partial\Omega,\ P_{1}=(x_{1},\ y_{1},\ z_{1})^{T}\in \partial\Omega$, then we can work out the following, \begin{eqnarray*}\int_{0}^{t}d\tau\int_{\partial\Omega}\cfrac{\partial h(M-P_{1},\ t-\tau)}{\partial n_{p_{1}}} A_{3}(P_{1},\ \tau) dS_{P_{1}}&=&\int_{0}^{t}d\tau\int_{\partial\Omega} C_{4}\ e^{-\cfrac{a_{0}}{(t-\tau)}}\cfrac{b_{0}A_{3}(P_{1},\ \tau)}{(t-\tau)^{5/2}} dS_{P_{1}}\\ (v=\cfrac{a_{0}}{(t-\tau)},\ \tau=t-\cfrac{a_{0}}{v},\ d\tau=\cfrac{a_{0}dv}{v^{2}})&=&\int_{\partial\Omega}dS_{P_{1}}\int_{\cfrac{a_{0}}{t}}^{+\infty}C_{4}e^{-v}b_{0}A_{3}(P_{1},\ t-\cfrac{a_{0}}{v}) a_{0}^{-5/2}v^{5/2}\cfrac{a_{0}dv}{v^{2}}\\&=&\int_{\partial\Omega}b_{0} a_{0}^{-3/2}dS_{P_{1}}\int_{\cfrac{a_{0}}{t}}^{+\infty}C_{4}e^{-v}v^{1/2}A_{3}(P_{1},\ t-\cfrac{a_{0}}{v})dv,\end{eqnarray*}
where\ $b_{0}={\bf r}_{MP_{1}}\cdot n_{p_{1}}$,\ $a_{0}=\mid M-P_{1}\mid^{2}/4=[(x-x_{1})^{2}+(y-y_{1})^{2}+(z-z_{1})^{2}]/4,\ C_{4}=1/(16(\sqrt{\pi })^{3})$.\\
From Theorem 2.1, we know there exists a neighborhood\ $U(M,\ \delta(M)),\ \delta(M)>0$,\ and\ $C(M)>0$, such that\ $|b_{0}|\leq C(M) a_{0}^{(1+\beta)/2}$,\ $\forall P_{1}\in U(M,\ \delta(M))\cap \partial\Omega$. And from lemma 2.2, we can get that\ $\partial h/\partial n_{p_{1}}.\ast_{\partial\Omega}A_{3}$\ is continuous on\ $\partial\Omega\times[0,\ T]$.\\
But if\ $M=(x,\ y,\ z)^{T}\in R^{3}\setminus \partial\Omega,\ P_{1}=(x_{1},\ y_{1},\ z_{1})^{T}\in \partial\Omega,\ t\in(0,\ T]$, then we can work out the following, \begin{eqnarray*}\int_{0}^{t}d\tau\int_{\partial\Omega}\cfrac{\partial h(M-P_{1},\ t-\tau)}{\partial n_{p_{1}}} A_{3}(P_{1},\ \tau) dS_{P_{1}}&=&\int_{0}^{t}d\tau\int_{\partial\Omega} C_{4}\ e^{-\cfrac{a_{0}}{(t-\tau)}}\cfrac{b_{0}A_{3}(P_{1},\ \tau)}{(t-\tau)^{5/2}} dS_{P_{1}}\\ (v=\cfrac{a_{0}}{(t-\tau)},\ \tau=t-\cfrac{a_{0}}{v},\ d\tau=\cfrac{a_{0}dv}{v^{2}})&=&\int_{\partial\Omega}dS_{P_{1}}\int_{\cfrac{a_{0}}{t}}^{+\infty}C_{4}e^{-v}b_{0}A_{3}(P_{1},\ t-\cfrac{a_{0}}{v}) a_{0}^{-5/2}v^{5/2}\cfrac{a_{0}dv}{v^{2}}\\&=&\int_{\partial\Omega}b_{0} a_{0}^{-3/2}dS_{P_{1}}\int_{\cfrac{a_{0}}{t}}^{+\infty}C_{4}e^{-v}v^{1/2}A_{3}(P_{1},\ t-\cfrac{a_{0}}{v})dv\\&=&I_{1}+I_{2},\end{eqnarray*}where\ $b_{0}={\bf r}_{MP_{1}}\cdot n_{p_{1}}$,\ $a_{0}=\mid M-P_{1}\mid^{2}/4=[(x-x_{1})^{2}+(y-y_{1})^{2}+(z-z_{1})^{2}]/4,\ C_{4}=1/(16(\sqrt{\pi })^{3})$.
\begin{eqnarray*}I_{1}&=&\int_{\partial\Omega}b_{0} a_{0}^{-3/2}dS_{P_{1}}\int_{\cfrac{a_{0}}{t}}^{+\infty}C_{4}e^{-v}v^{1/2}A_{3}(P_{1},\ t)dv\\&=&\int_{\partial\Omega}b_{0} a_{0}^{-3/2}A_{3}(P_{1},\ t)dS_{P_{1}}(\int_{0}^{+\infty}C_{4}e^{-v}v^{1/2}dv-\int_{0}^{\cfrac{a_{0}}{t}}C_{4}e^{-v}v^{1/2}dv)\\&=&C_{4}\Gamma(3/2)\int_{\partial\Omega}b_{0} a_{0}^{-3/2}A_{3}(P_{1},\ t)dS_{P_{1}}-\int_{\partial\Omega}b_{0} a_{0}^{-3/2}A_{3}(P_{1},\ t)dS_{P_{1}}\int_{0}^{\cfrac{a_{0}}{t}}C_{4}e^{-v}v^{1/2}dv\\&=& I_{3}-I_{4},\end{eqnarray*} $\Gamma(z)$\ is the Euler gamma function,
\begin{eqnarray*}I_{3}&=& C_{4}\Gamma(3/2)\int_{\partial\Omega}b_{0} a_{0}^{-3/2}A_{3}(P_{1},\ t)dS_{P_{1}},\\ I_{4}&=&\int_{\partial\Omega}b_{0} a_{0}^{-3/2}A_{3}(P_{1},\ t)dS_{P_{1}}\int_{0}^{\cfrac{a_{0}}{t}}C_{4}e^{-v}v^{1/2}dv\\
 I_{2}&=&\int_{\partial\Omega}b_{0} a_{0}^{-3/2}dS_{P_{1}}\int_{\cfrac{a_{0}}{t}}^{+\infty}C_{4}e^{-v}v^{1/2}[A_{3}(P_{1},\ t-\cfrac{a_{0}}{v})-A_{3}(P_{1},\ t)]dv.
\end{eqnarray*}
If we let$$C_{5}=\max_{v\in [0,\ +\infty)}e^{-v}v^{1/2},\ C_{6}=\max_{P_{1}\in\partial\Omega,\ t\in[0,\ T]}|\cfrac{\partial A_{3}(P_{1},\ t)}{\partial t}|,$$
then we can obtain the following\begin{eqnarray*} |I_{4}|&\leq &C_{4}C_{5}\int_{\partial\Omega}\cfrac{b_{0} a_{0}^{-1/2}}{t}|A_{3}(P_{1},\ t)|dS_{P_{1}},\\
 |I_{2}|&\leq &\int_{\partial\Omega}b_{0} a_{0}^{-3/2}dS_{P_{1}}\int_{\cfrac{a_{0}}{t}}^{+\infty}C_{4}C_{6}e^{-v}v^{1/2}\cfrac{a_{0}}{v}dv\\
&\leq & C_{4}C_{6}\Gamma(1/2)\int_{\partial\Omega}b_{0} a_{0}^{-1/2}dS_{P_{1}}.\end{eqnarray*} From lemma 2.2, we know\ $\forall P_{0}\in \partial\Omega,\ \forall t\in(0,\ T],\ I_{4}$\ and\ $I_{2}$\ are continuous at\ $P_{0}$. But from corollary 2.1 we have the following,
\[ \lim_{M\rightarrow P_{0}+} I_{3}(M,\ t)=I_{3}(P_{0},\ t)+\cfrac{1}{2}A_{3}(P_{0},\ t),\ \lim_{M\rightarrow P_{0}-} I_{3}(M,\ t)=I_{3}(P_{0},\ t)-\cfrac{1}{2}A_{3}(P_{0},\ t),\]
which completes the statement.\qed\\
\begin{corollary} If\ $\Omega$\ is bounded,\ $\partial\Omega\in C^{1,\ \beta},\ 0<\beta\leq 1,\ A_{3}(P_{1},\ t)\in C(\partial\Omega\times[0,\ T])$, then\ $\partial h/\partial n_{p_{1}}.\ast_{\partial\Omega}A_{3}$\ is continuous on\ $(R^{3}\setminus \partial\Omega)\times[0,\ T]$\ and\ $\partial\Omega\times[0,\ T]$, moreover (3.28) and (3.29) stand.\end{corollary}
{\it Proof of corollary 3.1}. If\ $A_{3}(P_{1},\ t)\in C(\partial\Omega\times[0,\ T])$, then from Weierstrass theorem, we can make\ $A_{3k}(P_{1},\ t)\in C^{2}(\partial\Omega)\cap C^{1}[0,\ T])$, such that\[\lim_{k\rightarrow +\infty}A_{3k}(P_{1},\ t)=A_{3}(P_{1},\ t),\ \mbox{uniformlly on}\ \partial\Omega\times[0,\ T].\] Then\ $A_{3k}(P_{1},\ t)$\ satisfies theorem 3.3. And from theorem 2.2, we know\[\int_{\partial\Omega}|b_{0} a_{0}^{-3/2}|dS_{P_{1}}<+\infty.\]So we obtain the following,
\begin{eqnarray*} \lim_{k\rightarrow +\infty}\int_{0}^{+\infty}e^{-v}v^{1/2}| A_{3k}(P_{1},\ t-\cfrac{a_{0}}{v})-A_{3}(P_{1},\ t-\cfrac{a_{0}}{v})|dv&=&0,\ \mbox{uniformally},\\ \lim_{k\rightarrow +\infty}\int_{\partial\Omega}|b_{0} a_{0}^{-3/2}[ A_{3k}(P_{1},\ t)-A_{3}(P_{1},\ t)]|dS_{P_{1}}&=&0,\ \mbox{uniformally}.\end{eqnarray*}  Hence the statement holds.\qed\\
\begin{corollary} If\ $\Omega$\ is bounded,\ $\partial\Omega\in C^{1,\ \beta},\ 0<\beta\leq 1,\ A_{3}(P_{1},\ t)\in C(\partial\Omega\times[0,\ T])$, then\ $\partial h/\partial n_{p_{0}}.\ast_{\partial\Omega}A_{3}$\ is continuous on\ $(R^{3}\setminus \partial\Omega)\times[0,\ T]$\ and\ $\partial\Omega\times[0,\ T]$, but\ $\forall P_{0}\in \partial\Omega,\ \forall t\in(0,\ T]$, we have\begin{eqnarray}&& \lim_{M\rightarrow n_{p_{0}}^{+}}(\cfrac{\partial h}{\partial n_{p_{0}}}.\ast_{\partial\Omega}A_{3})(M,\ t)=(\cfrac{\partial h}{\partial n_{p_{0}}}.\ast_{\partial\Omega}A_{3})(P_{0},\ t)-\cfrac{1}{2} A_{3}(P_{0},\ t),\\&&  \lim_{M\rightarrow n_{p_{0}}^{-}}(\cfrac{\partial h}{\partial n_{p_{0}}}.\ast_{\partial\Omega}A_{3})(M,\ t)=(\cfrac{\partial h}{\partial n_{p_{0}}}.\ast_{\partial\Omega}A_{3})(P_{0},\ t)+\cfrac{1}{2} A_{3}(P_{0},\ t),\end{eqnarray}
 where$$(\cfrac{\partial h}{\partial n_{p_{0}}}.\ast_{\partial\Omega}A_{3})(M,\ t)=\int_{0}^{t}d\tau\int_{\partial\Omega}[\nabla_{M}h(M-P_{1},\ t-\tau)\cdot n_{p_{0}}]A_{3}(P_{1},\ \tau) dS_{P_{1}}, $$
 $M\rightarrow n_{p_{0}}^{+}$\ means\ $M$\ is near to\ $P_{0}$\ along\ $n_{P_{0}}$\ from the exterior of\ $\Omega$\ and\ $M\rightarrow n_{p_{0}}^{-}$\ means\ $M$\ is near to\ $P_{0}$\ along\ $n_{P_{0}}$\ from the interior of\ $\Omega$.\end{corollary}
{\it Proof of corollary 3.2}. We can obtain the proof from theorem 2.3 and the previous corollary.\qed\\ From theorem 2.1, we know\ $\partial h/\partial n_{p_{1}}.\ast_{\partial\Omega}A_{3}$\ is only continuous on\ $\Omega\times[0,\ T]$. So from (3.27), we only get\[ uI_{\Omega\times[0,\ T]}=[h.\ast_{\overline{\Omega}}A_{2}+h.\ast_{\partial\Omega}A_{7}+\cfrac{\partial h}{\partial n_{p_{1}}}.\ast_{\partial\Omega}A_{3}+h.\ast \psi_{1}(Z_{1})]I_{\Omega\times[0,\ T]}.\] This is enough. We can get\ $A_{2}$\ as follows,
\[ A_{2}I_{\Omega}=[h.\ast_{\overline{\Omega}}A_{2}+h.\ast_{\partial\Omega}A_{7}+\cfrac{\partial h}{\partial n_{p_{1}}}.\ast_{\partial\Omega}A_{3}+h.\ast \psi_{1}(Z_{1})]I_{\Omega}|_{t=T}.\] We need it to test whether (3.25) is entire or not. Next we see\ $h.\ast \psi_{1}(Z_{1})$.\\
\begin{lemma} \label{lemma1} If\ $\Omega$\ is bounded,\ $\partial\Omega\in C^{1,\ \beta},\ 0<\beta\leq 1,\ f$\ is continuous, then\ $w=h.\ast f\in C^{1}(\overline{\Omega})$.
\end{lemma}We can see the proof on page 54 of [2].\\Now we can see that \[uI_{\Omega\times[0,\ T]}=[h.\ast_{\overline{\Omega}}A_{1}+h.\ast_{\partial\Omega}A_{7}+\cfrac{\partial h}{\partial n_{p_{1}}}.\ast_{\partial\Omega}A_{3}+h.\ast \psi_{1}(Z_{1})]I_{\Omega\times[0,\ T]}\in C^{1}(\Omega)\cap C[0,\ T].\]
But we find that there is only one of\ $A_{3},\ A_{7}$\ is known in (3.37). We should determinate both of these. Following the Lyapunov's potential theory on pages 173 to 201 of [6], we discuss three boundary conditions as follows.\\
(1)Dirichlet problem. If\ $A_{1},\ A_{3}$\ are known, then from (3.37) we can get the following,\ $\forall M\in \Omega$, \[ \cfrac{\partial u(M,\ t)I_{\Omega\times[0,\ T]}}{\partial n_{p_{0}}}=(\cfrac{\partial g_{01}}{\partial n_{p_{0}}}+\cfrac{\partial h(M,\ t)}{\partial n_{p_{0}}}.\ast_{\partial\Omega}A_{7}+\cfrac{\partial h(M,\ t)}{\partial n_{p_{0}}}.\ast \psi_{1}(Z_{1}))I_{\Omega\times[0,\ T]},\] where$$ g_{01}=h.\ast_{\overline{\Omega}}A_{1}+\cfrac{\partial h}{\partial n_{p_{1}}}.\ast_{\partial\Omega}A_{3}.$$
As\ $u(M,\ t)\in C^{2}(\overline{\Omega})\cap C^{1}[0,\ T]$, and from Corollary 2.2, we can get the following,
$$ \lim_{M\rightarrow n_{p_{0}}^{-}} \cfrac{\partial u(M,\ t)I_{\Omega\times[0,\ T]}}{\partial n_{p_{0}}}=A_{7},\ \lim_{M\rightarrow n_{p_{0}}^{-}}\cfrac{\partial h(M,\ t)}{\partial n_{p_{0}}}.\ast_{\partial\Omega}A_{7}=\cfrac{\partial h(P_{0},\ t)}{\partial n_{p_{0}}}.\ast_{\partial\Omega}A_{7}+\cfrac{1}{2}A_{7},\ t\in(0,\ T]. $$ Hence there exists\ $g_{02}\in C^{1}(\partial\Omega\times [0,\ T])$, such that
\[  \lim_{M\rightarrow n_{p_{0}}^{-}} \cfrac{\partial g_{01}(M,\ t)}{\partial n_{p_{0}}}=g_{02}(P_{0},\ t),\ \forall P_{0}\in \partial\Omega.\] From the continuity, we get the second type of linear Fredholm integral equations that\ $A_{7}$\ should satisfy,
\[ \cfrac{1}{2}A_{7}=\cfrac{\partial h(P_{0},\ t)}{\partial n_{p_{0}}}.\ast_{\partial\Omega}A_{7}+g_{02}+\cfrac{\partial h(P_{0},\ t)}{\partial n_{p_{0}}}.\ast \psi_{1}(Z_{1}).\]
We will prove that there is only\ $0$\ for the homogeneous equations as follows, \[ \cfrac{1}{2}A_{7}=\cfrac{\partial h(P_{0},\ t)}{\partial n_{p_{0}}}.\ast_{\partial\Omega}A_{7}.\] If\ $A_{7}$\ satisfies (3.41), then we may let\[ W(M,\ t)=h(M,\ t).\ast_{\partial\Omega}A_{7}=\int_{0}^{t}d\tau\int_{\partial\Omega}h(M-P_{1},\ t-\tau)A_{7}(P_{1},\ \tau)dS_{P_{1}},\]
and we can obtain that\ $W$\ is continuous on\ $R^{3}$, moreover
\[W_{t}-\triangle W=0,\ \mbox{on}\ (R^{3}\setminus\partial\Omega)\times[0,\ T],\ W|_{t=0}=0.\] From the fact that\ $A_{7}$\ satisfies (3.41), we can get
\[ \cfrac{\partial W}{\partial  n_{p_{0}}^{+}}|_{\partial \Omega}=\lim_{M\rightarrow n_{p_{0}}^{+}}\cfrac{\partial  W(M,\ t)}{\partial n_{p_{0}}} =\cfrac{\partial h(P_{0},\ t)}{\partial n_{p_{0}}}.\ast_{\partial\Omega}A_{7}-\cfrac{1}{2}A_{7}=0,\ t\in(0,\ T].\]
From (3.42), we know\ $W$\ is rapid descent in\ $R^{3}$. By the uniqueness of the solution of Eqs(3.43), we get\ $W\equiv0$, on\ $(R^{3}\setminus\overline{\Omega})\times[0,\ T]$. From the continuity we can see that\ $W|_{\partial\Omega}=0$. Hence, we can deduce that\[W\equiv0\ \mbox{on}\ R^{3}\times[0,\ T]. \] This means that\[ \cfrac{\partial W}{\partial  n_{p_{0}}^{+}}|_{\partial \Omega}-\cfrac{\partial W}{\partial  n_{p_{0}}^{-}}|_{\partial \Omega}=-A_{7}=0.\]
Hence, the solution of (3.41) is only\ $0$. From Fredholm integral equation theory, we know there is only one solution for (3.40). Moreover, there exists an entire analytic function\ $\Gamma_{1}$\ that is only related to\ $\partial h(P_{0},\ t)/\partial n_{p_{0}}$, such that
\begin{eqnarray}\cfrac{A_{7}(P_{0},\ t)}{2}&=&g_{02}(P_{0},\ t)+\cfrac{\partial h(P_{0},\ t)}{\partial n_{p_{0}}}.\ast \psi_{1}(Z_{1})+\nonumber\\&&\int_{0}^{T}d\tau\int_{\partial\Omega}\Gamma_{1}(P_{0},\ t,\ P_{1},\ \tau)[g_{02}(P_{1},\ \tau)+\cfrac{\partial h(P_{1},\ \tau)}{\partial n_{p_{1}}}.\ast \psi_{1}(Z_{1})]dS_{P_{1}}.\end{eqnarray}
However, it is true that\ $\partial h(P_{0},\ t)/\partial n_{p_{0}}(P_{0}-P_{1},\ t-\tau)$\ is not continuous, if\ $P_{0}=P_{1},\ t=\tau$. Therefore, it's lucky that \[(\sqrt{|P_{0}-P_{1}|^{2}+(t-\tau)^{2}})^{(\epsilon_{0}+5/2)}\cfrac{\partial h(P_{0},\ t)}{\partial n_{p_{0}}}(P_{0}-P_{1},\ t-\tau)\]
is continuous, if\ $\epsilon_{0}\in(0,\ 0.5]$. Hence,\ $\partial h(P_{0},\ t)/\partial n_{p_{0}}(P_{0}-P_{1},\ t-\tau)$\ is a weak singular kernel. From Fredholm theorem, we know (3.47) still stands.\\
From\ $ A_{7}\in C^{1}(\partial\Omega\times [0,\ T])$, we know\ $g_{02}$\ should be in\ $C^{1}(\partial\Omega\times [0,\ T])$.\\
(2)Neumann problem. If\ $A_{1},\ A_{7}$\ are known, then from (3.37) we can get the following,\ $\forall M\in \Omega$, \[  u(M,\ t)I_{\Omega\times[0,\ T]}=(g_{03}+\cfrac{\partial h(M,\ t)}{\partial n_{p_{1}}}.\ast_{\partial\Omega}A_{3}+ h(M,\ t).\ast \psi_{1}(Z_{1}))I_{\Omega\times[0,\ T]},\] where$$ g_{03}=h.\ast_{\overline{\Omega}}A_{1}+h.\ast_{\partial\Omega}A_{7}.$$
Because\ $u(M,\ t)\in C^{2}(\overline{\Omega})\cap C^{1}[0,\ T]$, and from Theorem 3.3, we can get as follows,
$$ \lim_{M\rightarrow P_{0}+}  u(M,\ t)I_{\Omega\times[0,\ T]}=A_{3},\ \lim_{M\rightarrow P_{0}+}\cfrac{\partial h(M,\ t)}{\partial n_{p_{1}}}.\ast_{\partial\Omega}A_{3}=\cfrac{\partial h(P_{0},\ t)}{\partial n_{p_{1}}}.\ast_{\partial\Omega}A_{3}+\cfrac{1}{2}A_{3},\ t\in(0,\ T]. $$ Hence, there exists\ $g_{04}\in C^{2}(\partial\Omega)\cap C^{1}[0,\ T])$, such that
\[  \lim_{M\rightarrow P_{0}+} g_{03}(M,\ t)=g_{04}(P_{0},\ t),\ \forall P_{0}\in \partial\Omega.\] From the continuity, we get the second type of linear Fredholm integral equations that\ $A_{3}$\ should satisfy,
\[ \cfrac{1}{2}A_{3}=\cfrac{\partial h(P_{0},\ t)}{\partial n_{p_{1}}}.\ast_{\partial\Omega}A_{3}+g_{04}+ h(P_{0},\ t).\ast \psi_{1}(Z_{1}).\]
We will prove that there is only \ $0$\ for the homogeneous equations as follows, \[ \cfrac{1}{2}A_{3}=\cfrac{\partial h(P_{0},\ t)}{\partial n_{p_{1}}}.\ast_{\partial\Omega}A_{3}.\] We write the detail of (3.52) as follows,
\begin{eqnarray*} \cfrac{1}{2}A_{3}(P_{0},\ t)&=&\cfrac{\partial h(P_{0},\ t)}{\partial n_{p_{1}}}.\ast_{\partial\Omega}A_{3}\\
&=&\int_{0}^{T}d\tau\int_{\partial\Omega}\cfrac{1}{16(\sqrt{\pi})^{3}}(t-\tau)^{-5/2}e^{-|P_{0}-P_{1}|^{2}/[4(t-\tau)]}I_{\{t>\tau\}}[-(P_{0}-P_{1})\cdot n_{p_{1}}]A_{3}(P_{1},\ \tau)dS_{P_{1}}.\end{eqnarray*}
From (3.48), we know\ $\partial h(P_{0},\ t)/\partial n_{p_{1}}(P_{0}-P_{1},\ t-\tau)$\ is also a weak singular kernel.
Hence, from Fredholm theorem, we may transpose the equation (3.52) as follows,
\begin{eqnarray}\cfrac{1}{2}A_{3}(P_{0},\ t)&=&\int_{0}^{T}d\tau\int_{\partial\Omega}\cfrac{1}{16(\sqrt{\pi})^{3}}(\tau-t)^{-5/2}e^{-|P_{0}-P_{1}|^{2}/[4(\tau-t)]}I_{\{\tau>t\}}[(P_{0}-P_{1})\cdot n_{p_{0}}]A_{3}(P_{1},\ \tau)dS_{P_{1}}\nonumber\\&=& -\cfrac{\partial h_{1}(P_{0},\ t)}{\partial n_{p_{0}}}.\ast_{\partial\Omega}A_{3},\end{eqnarray}where\ $h_{1}(M,\ t)=(2\sqrt{\pi})^{-3}(-t)^{-3/2}e^{|M|^{2}/4t}I_{\{t<0\}},\ |M|^{2}=x^{2}+y^{2}+z^{2}$. \\If\ $h_{1}$\ satisfies\ $(h_{1})_{t}-\triangle h_{1}=0$, then from case (1), we can see that there is only\ $0$\ for (3.53). But something amazing happened, and\ $h_{1}$\ satisfies\ $(h_{1})_{t}+\triangle h_{1}=0$. We can't use the method in case (1) unless there is only\ $0$\ for\[(W_{2})_{t}+\triangle W_{2}=0,\ \mbox{on}\ \overline{\Omega}\times [0,\ T],\ W_{2}|_{t=0}=0,\ \mbox{one of}\ W_{2}|_{\partial\Omega},\ \cfrac{\partial W_{2}}{\partial n}|_{\partial\Omega},\ \cfrac{\partial W_{2}}{\partial n}+\sigma W_{2}|_{\partial\Omega},\ \sigma>0,\ \mbox{is}\ 0.\] We are not sure about (3.54). It looks that we are stumped by the finish line. But if we prove there is only\ $0$\ for (3.53), then (3.54) should be true.\\ The answer is not very complex. Letting$$\tau=T-\tau_{1},\ t=T-t_{1},\ A_{3}^{\prime}(P_{0},\ t_{1})=A_{3}(P_{0},\ T-t_{1}),$$
then we can obtain,\[\cfrac{1}{2}A_{3}^{\prime}(P_{0},\ t_{1})=-\cfrac{\partial h(P_{0},\ t_{1})}{\partial n_{p_{0}}}.\ast_{\partial\Omega}A_{3}^{\prime}(P_{0},\ t_{1}).\]  Following the idea in case (1), letting\[ W_{1}(M,\ t_{1})=h(M,\ t_{1}).\ast_{\partial\Omega}A_{3}^{\prime}=\int_{0}^{T}d\tau\int_{\partial\Omega}h(M-P_{1},\ t_{1}-\tau)A_{3}^{\prime}(P_{1},\ \tau)dS_{P_{1}},\]
then we can obtain that\ $W_{1}$\ is continuous on\ $R^{3}$, moreover
\[(W_{1})_{t_{1}}-\triangle W_{1}=0,\ \mbox{on}\ (R^{3}\setminus\partial\Omega)\times[0,\ T],\ W_{1}|_{t_{1}=0}=0.\] From\ $A_{3}^{\prime}$\ satisfying (3.55), we can get
\[ \cfrac{\partial W_{1}}{\partial  n_{p_{0}}^{-}}|_{\partial \Omega}=\lim_{M\rightarrow n_{p_{0}}^{-}}\cfrac{\partial  W_{1}(M,\ t_{1})}{\partial n_{p_{0}}} =\cfrac{\partial h(P_{0},\ t_{1})}{\partial n_{p_{0}}}.\ast_{\partial\Omega}A_{3}^{\prime}+\cfrac{1}{2}A_{3}^{\prime}=0,\ t_{1}\in(0,\ T].\]
By the uniqueness of the solution of Eqs(3.57), we get\ $W_{1}\equiv0$, on\ $\Omega\times[0,\ T]$. From the continuity we can get that\ $W_{1}|_{\partial\Omega}=0$. From (3.56), we know\ $W_{1}$\ is rapid descent in\ $R^{3}$. Hence we can get that\[W_{1}\equiv0\ \mbox{on}\ R^{3}\times[0,\ T]. \]  This means that\[ \cfrac{\partial W_{1}}{\partial  n_{p_{0}}^{+}}|_{\partial \Omega}-\cfrac{\partial W_{1}}{\partial  n_{p_{0}}^{-}}|_{\partial \Omega}=-A_{3}^{\prime}=0.\]
So there is only\ $0$\ for (3.55). Moreover, there is only\ $0$\ for (3.53). From Fredholm theorem, we know there is only\ $0$\ for (3.52). Hence there exists an entire analytic function\ $\Gamma_{2}$\ that is only related to\ $\partial h(P_{0},\ t)/\partial n_{p_{1}}$, such that
\begin{eqnarray} \cfrac{A_{3}(P_{0},\ t)}{2}&=&g_{04}(P_{0},\ t)+h(P_{0},\ t).\ast \psi_{1}(Z_{1})+\nonumber\\&&\int_{0}^{T}d\tau\int_{\partial\Omega}\Gamma_{2}(P_{0},\ t,\ P_{1},\ \tau)[g_{04}(P_{1},\ \tau)+ h(P_{1},\ \tau).\ast \psi_{1}(Z_{1})]dS_{P_{1}}.\end{eqnarray}
(3)Robin problem. If\ $A_{1},\ A_{8}=A_{7}+\sigma A_{3},$\ are known, then from (3.37) we can get the following,\ $\forall M\in \Omega$, \[  u(M,\ t)I_{\Omega\times[0,\ T]}=(g_{05}-h.\ast_{\partial\Omega}(\sigma A_{3})+\cfrac{\partial h(M,\ t)}{\partial n_{p_{1}}}.\ast_{\partial\Omega}A_{3}+ h(M,\ t).\ast \psi_{1}(Z_{1}))I_{\Omega\times[0,\ T]},\] where$$ g_{05}=h.\ast_{\overline{\Omega}}A_{1}+h.\ast_{\partial\Omega}A_{8}.$$
Since\ $u(M,\ t)\in C^{2}(\overline{\Omega})\cap C^{1}[0,\ T]$, and from Theorem 2.1, we can get,
$$ \lim_{M\rightarrow P_{0}+}  u(M,\ t)I_{\Omega\times[0,\ T]}=A_{3},\ \lim_{M\rightarrow P_{0}+}\cfrac{\partial h(M,\ t)}{\partial n_{p_{1}}}.\ast_{\partial\Omega}A_{3}=\cfrac{\partial h(P_{0},\ t)}{\partial n_{p_{1}}}.\ast_{\partial\Omega}A_{3}+\cfrac{1}{2}A_{3},\ t\in(0,\ T]. $$ Hence, there exists\ $g_{06}\in C^{2}(\partial\Omega)\cap C^{1}[0,\ T])$, such that
\[  \lim_{M\rightarrow P_{0}+} g_{05}(M,\ t)=g_{06}(P_{0},\ t),\ \forall P_{0}\in \partial\Omega.\] From the continuity, we get the second type of linear Fredholm integral equations that\ $A_{3}$\ should satisfy,
\[ \cfrac{1}{2}A_{3}=-h(P_{0},\ t).\ast_{\partial\Omega}(\sigma A_{3})+\cfrac{\partial h(P_{0},\ t)}{\partial n_{p_{1}}}.\ast_{\partial\Omega}A_{3}+g_{06}+ h(P_{0},\ t).\ast \psi_{1}(Z_{1}).\]
We will prove that there is only \ $0$\ for the homogeneous equations as follows, \[ \cfrac{1}{2}A_{3}=-h(P_{0},\ t).\ast_{\partial\Omega}(\sigma A_{3})+\cfrac{\partial h(P_{0},\ t)}{\partial n_{p_{1}}}.\ast_{\partial\Omega}A_{3}.\] We write the detail of (3.65) as follows,
\begin{eqnarray*} \cfrac{1}{2}A_{3}(P_{0},\ t)&=&-h(P_{0},\ t).\ast_{\partial\Omega}(\sigma A_{3})+\cfrac{\partial h(P_{0},\ t)}{\partial n_{p_{1}}}.\ast_{\partial\Omega}A_{3}\\
&=&\int_{0}^{T}d\tau\int_{\partial\Omega}\cfrac{-1}{(2\sqrt{\pi})^{3}}(t-\tau)^{-3/2}e^{-|P_{0}-P_{1}|^{2}/[4(t-\tau)]}I_{\{t>\tau\}}\sigma(P_{1},\ \tau)A_{3}(P_{1},\ \tau)dS_{P_{1}}+
\\&&\int_{0}^{T}d\tau\int_{\partial\Omega}\cfrac{1}{16(\sqrt{\pi})^{3}}(t-\tau)^{-5/2}e^{-|P_{0}-P_{1}|^{2}/[4(t-\tau)]}I_{\{t>\tau\}}[-(P_{0}-P_{1})\cdot n_{p_{1}}]A_{3}(P_{1},\ \tau)dS_{P_{1}}.\end{eqnarray*}
From Fredholm theorem, we may transpose equation (3.65) as follows,
\begin{eqnarray}\cfrac{1}{2}A_{3}(P_{0},\ t)&=&\int_{0}^{T}d\tau\int_{\partial\Omega}\cfrac{-1}{(2\sqrt{\pi})^{3}}(\tau-t)^{-3/2}e^{-|P_{0}-P_{1}|^{2}/[4(\tau-t)]}I_{\{\tau>t\}}\sigma(P_{0},\ t)A_{3}(P_{1},\ \tau)dS_{P_{1}}+\nonumber
\\&&\int_{0}^{T}d\tau\int_{\partial\Omega}\cfrac{1}{16(\sqrt{\pi})^{3}}(\tau-t)^{-5/2}e^{-|P_{0}-P_{1}|^{2}/[4(\tau-t)]}I_{\{\tau>t\}}[(P_{0}-P_{1})\cdot n_{p_{0}}]A_{3}(P_{1},\ \tau)dS_{P_{1}}\nonumber\\&=& \sigma(P_{0},\ t)( -h_{1}(P_{0},\ t).\ast_{\partial\Omega}A_{3})-\cfrac{\partial h_{1}(P_{0},\ t)}{\partial n_{p_{0}}}.\ast_{\partial\Omega}A_{3},\end{eqnarray}where\ $h_{1}(M,\ t)=(2\sqrt{\pi})^{-3}(-t)^{-3/2}e^{|M|^{2}/4t}I_{\{t<0\}},\ |M|^{2}=x^{2}+y^{2}+z^{2}$. \\Letting$$\tau=T-\tau_{1},\ t=T-t_{1},\ A_{3}^{\prime}(P_{0},\ t_{1})=A_{3}(P_{0},\ T-t_{1}),\ \sigma^{\prime}(P_{0},\ t_{1})=\sigma(P_{0},\ T-t_{1}),$$
then we can get as follows,\[\cfrac{1}{2}A_{3}^{\prime}(P_{0},\ t_{1})=\sigma^{\prime}(P_{0},\ t_{1})( -h(P_{0},\ t_{1}).\ast_{\partial\Omega}A_{3}^{\prime})-\cfrac{\partial h(P_{0},\ t_{1})}{\partial n_{p_{0}}}.\ast_{\partial\Omega}A_{3}^{\prime}.\]
Letting\[ W^{\prime}(M,\ t_{1})=h(M,\ t_{1}).\ast_{\partial\Omega}A_{3}^{\prime}=\int_{0}^{T}d\tau\int_{\partial\Omega}h(M-P_{1},\ t_{1}-\tau)A_{3}^{\prime}(P_{1},\ \tau)dS_{P_{1}},\]
then we can obtain that\ $W^{\prime}$\ is continuous on\ $R^{3}$, moreover,
\[W^{\prime}_{t_{1}}-\triangle W^{\prime}=0,\ \mbox{on}\ (R^{3}\setminus\partial\Omega)\times[0,\ T],\ W^{\prime}|_{t_{1}=0}=0.\] We obtain
\[ \cfrac{\partial W^{\prime}}{\partial  n_{p_{0}}^{-}}|_{\partial \Omega}=\lim_{M\rightarrow n_{p_{0}}^{-}}\cfrac{\partial  W^{\prime}(M,\ t_{1})}{\partial n_{p_{0}}} =\cfrac{\partial h(P_{0},\ t_{1})}{\partial n_{p_{0}}}.\ast_{\partial\Omega}A_{3}^{\prime}+\cfrac{1}{2}A_{3}^{\prime},\ t_{1}\in(0,\ T].\]
Since\ $A_{3}^{\prime}$\ satisfies (3.67), we can get that
\[\cfrac{\partial W^{\prime}}{\partial  n_{p_{0}}^{-}}+\sigma^{\prime}(P_{0},\ t_{1})W^{\prime}|_{\partial \Omega}=0.\]From\ $W^{\prime}|_{t_{1}=0}=0$, we get\ $W^{\prime}\equiv0$, on\ $\Omega\times[0,\ T]$. So\ $W^{\prime}|_{\partial\Omega}=0$. We know\ $W^{\prime}$\ is also rapid descent in\ $R^{3}$. Hence, we can get that \[W^{\prime}\equiv0\ \mbox{on}\ R^{3}\times[0,\ T].\] This means that\[ \cfrac{\partial W^{\prime}}{\partial  n_{p_{0}}^{+}}-\cfrac{\partial W^{\prime}}{\partial  n_{p_{0}}^{-}}=-A_{3}^{\prime}=0.\]
Hence, the solution of (3.67) is only\ $0$. From Fredholm integral equation theory, we know there is only one solution for (3.64). Moreover, there exists an entire analytic function\ $\Gamma_{3}$\ that is only related to\ $\partial h(P_{0},\ t)/\partial n_{p_{0}}$\ and\ $h(P_{0},\ t)$, such that
\begin{eqnarray} \cfrac{A_{3}(P_{0},\ t)}{2}&=&g_{06}(P_{0},\ t)+ h(P_{0},\ t).\ast \psi_{1}(Z_{1})+\nonumber\\&&\int_{0}^{T}d\tau\int_{\partial\Omega}\Gamma_{3}(P_{0},\ t,\ P_{1},\ \tau)[g_{06}(P_{1},\ \tau)+h(P_{1},\ \tau).\ast \psi_{1}(Z_{1})]dS_{P_{1}}.\end{eqnarray}
Now we can introduce an easy way to get\ $A_{4},\ A_{5},\ A_{6}$, which we will use in the next section.\\
For any\ $P_{0}\in \partial\Omega$, we discuss a smooth curve on\ $\partial\Omega$\ that passes through\ $P_{0}$. The parameter coordinates of the point\ $P$\ on the curve are\ $(x(\theta),\ y(\theta),\ z(\theta))$. Moreover, the parameter coordinates of \ $P_{0}$\ are\ $(x(\theta_{0}),\ y(\theta_{0}),\ z(\theta_{0}))$. We assume the tangent vector at\ $P_{0}$\ is as follows,
$$ s_{0}=(x^{\prime}(\theta_{0}),\ y^{\prime}(\theta_{0}),\ z^{\prime}(\theta_{0}))^{T}.$$
If we select a\ $P(x(\theta),\ y(\theta),\ z(\theta))$\ on the curve that is near to\ $P_{0}$, then we can get\[ u(P,\ t)-u(P_{0},\ t)=u_{x}(P_{0},\ t)(x(\theta)-x(\theta_{0}))+u_{y}(P_{0},\ t)(y(\theta)-y(\theta_{0}))+u_{z}(P_{0},\ t)(z(\theta)-z(\theta_{0}))+o(\rho),\] where\ $\rho=|P-P_{0}|=\sqrt{(x(\theta)-x(\theta_{0}))^{2}+(y(\theta)-y(\theta_{0}))^{2}+(z(\theta)-z(\theta_{0}))^{2}}$.\\
From$$u(P,\ t)|_{P\in\partial\Omega}=A_{3}((P,\ t),$$ we get,\[A_{3}(P,\ t)-A_{3}(P_{0},\ t)=A_{3x}(P_{0},\ t)(x(\theta)-x(\theta_{0}))+A_{3y}(P_{0},\ t)(y(\theta)-y(\theta_{0}))+A_{3z}(P_{0},\ t)(z(\theta)-z(\theta_{0}))+o(\rho),\]where
$$A_{3x}=\cfrac{\partial A_{3}}{\partial x},\ A_{3y}=\cfrac{\partial A_{3}}{\partial y},\ A_{3z}=\cfrac{\partial A_{3}}{\partial z}.$$
By the subtraction of (3.75) and (3.76), we can obtain \begin{eqnarray}&&(u_{x}(P_{0},\ t)-A_{3x}(P_{0},\ t))(x(\theta)-x(\theta_{0}))+(u_{y}(P_{0},\ t)-A_{3y}(P_{0},\ t))(y(\theta)-y(\theta_{0}))\nonumber\\&&+(u_{z}(P_{0},\ t)-A_{3z}(P_{0},\ t))(z(\theta)-z(\theta_{0}))+o(\rho)=0.\end{eqnarray}
If we divide\ $\rho$\ on the left hand side of (3.77), and let\ $\theta\rightarrow \theta_{0}$, then from the following,
\begin{eqnarray*}\lim_{\theta\rightarrow \theta_{0}}\cfrac{x(\theta)-x(\theta_{0})}{\rho}&=&\lim_{\theta\rightarrow \theta_{0}}\cfrac{(x(\theta)-x(\theta_{0}))/(\theta-\theta_{0})}{\rho/(\theta-\theta_{0})}=\cfrac{x^{\prime}(\theta_{0})}{\parallel s_{0}\parallel},\\
\lim_{\theta\rightarrow \theta_{0}}\cfrac{y(\theta)-y(\theta_{0})}{\rho}&=&\lim_{\theta\rightarrow \theta_{0}}\cfrac{(y(\theta)-y(\theta_{0}))/(\theta-\theta_{0})}{\rho/(\theta-\theta_{0})}=\cfrac{y^{\prime}(\theta_{0})}{\parallel s_{0}\parallel},\\
\lim_{\theta\rightarrow \theta_{0}}\cfrac{z(\theta)-z(\theta_{0})}{\rho}&=&\lim_{\theta\rightarrow \theta_{0}}\cfrac{(z(\theta)-z(\theta_{0}))/(\theta-\theta_{0})}{\rho/(\theta-\theta_{0})}=\cfrac{z^{\prime}(\theta_{0})}{\parallel s_{0}\parallel},\end{eqnarray*}
where\ $\parallel s_{0}\parallel=\sqrt{(x^{\prime}(\theta_{0}))^{2}+(y^{\prime}(\theta_{0}))^{2}+(z^{\prime}(\theta_{0}))^{2}}$, we get,\[(u_{x}(P_{0},\ t)-A_{3x}(P_{0},\ t))x^{\prime}(\theta_{0})+(u_{y}(P_{0},\ t)-A_{3y}(P_{0},\ t))y^{\prime}(\theta_{0})+(u_{z}(P_{0},\ t)-A_{3z}(P_{0},\ t))z^{\prime}(\theta_{0})=0.\]
(3.78) will still stand even if\ $\parallel s_{0}\parallel=0$. From the arbitrary nature of\ $s_{0}$, we know that$$((u_{x}(P_{0},\ t)-A_{3x}(P_{0},\ t))_{k},\ (u_{y}(P_{0},\ t)-A_{3y}(P_{0},\ t))_{k},\ (u_{z}(P_{0},\ t)-A_{3z}(P_{0},\ t))_{k})^{T},\ 1\leq k\leq m,$$should be parallel to the exterior normal vector\ $n_{p_{0}}=(n_{1}(P_{0}),\ n_{2}(P_{0}),\ n_{3}(P_{0}))^{T}$, where\ $(\ast)_{k}$\ is the\ $k$th component of the vector,\ $1\leq k\leq m$. \\ Hence there exists\ $\lambda(P_{0},\ t)$, such that\begin{eqnarray*}u_{x}(P_{0},\ t)-A_{3x}(P_{0},\ t)&=&\lambda(P_{0},\ t)n_{1}(P_{0}),\\u_{y}(P_{0},\ t)-A_{3y}(P_{0},\ t)&=&\lambda(P_{0},\ t)n_{2}(P_{0}),\\u_{z}(P_{0},\ t)-A_{3z}(P_{0},\ t)&=&\lambda(P_{0},\ t)n_{3}(P_{0}).\end{eqnarray*}
From\ $n_{1}^{2}(P_{0})+n_{2}^{2}(P_{0})+n_{3}^{2}(P_{0})=1$, we can get,
\begin{eqnarray*}\lambda(P_{0},\ t)&=&(u_{x}(P_{0},\ t)-A_{3x}(P_{0},\ t))n_{1}(P_{0})+(u_{y}(P_{0},\ t)-A_{3y}(P_{0},\ t))n_{2}(P_{0})+\\&&(u_{z}(P_{0},\ t)-A_{3z}(P_{0},\ t))n_{3}(P_{0})\\&=&A_{7}(P_{0},\ t)-\cfrac{\partial A_{3}}{\partial n_{p_{0}}}(P_{0},\ t).\end{eqnarray*}
Hence, we can work out\ $A_{4},\ A_{5},\ A_{6}$\ as follows, \begin{eqnarray}A_{4}&=&u_{x}|_{\partial\Omega}=\cfrac{\partial A_{3}}{\partial x}+(A_{7}-\cfrac{\partial A_{3}}{\partial n})n_{1},\\A_{5}&=&u_{y}|_{\partial\Omega}=\cfrac{\partial A_{3}}{\partial y}+(A_{7}-\cfrac{\partial A_{3}}{\partial n})n_{2},\\A_{6}&=&u_{z}|_{\partial\Omega}=\cfrac{\partial A_{3}}{\partial z}+(A_{7}-\cfrac{\partial A_{3}}{\partial n})n_{3},\end{eqnarray}
where\ $n=(n_{1},\ n_{2},\ n_{3})^{T}$\ is the exterior normal vector to\ $\partial\Omega$.\\
Finally, we determine\ $A_{9}=p|_{\partial\Omega\times[0,\ T]}$. From the following
\begin{eqnarray*}(B^{-1}_{1})_{4}&=&(B_{03},\ B_{04}e_{1}^{T},\ B_{04}e_{2}^{T},\ B_{04}e_{3}^{T},\ -B_{03},\ 0_{1\times3},\ 0_{1\times3},\ -B_{03},\ 0_{1\times3},\ -B_{03},\ \tau B_{03}),\\ (B^{-1}_{1}B_{2})_{4}&=&(B_{03}),
        \end{eqnarray*} where\ $e_{1},\ e_{2},\ e_{3}$\ are all three dimensional unit coordinate vectors,\ $(B^{-1}_{1})_{4},\ (B^{-1}_{1}B_{2})_{4}$\ are the\ $4$th rows of\ $(B^{-1}_{1}),\ (B^{-1}_{1}B_{2})$, and
        \begin{eqnarray*}&& B_{03}=(\cfrac{i\xi_{1}}{\tau aa_{1}},\ \cfrac{i\xi_{2}}{\tau aa_{1}},\ \cfrac{i\xi_{3}}{\tau aa_{1}}),\ B_{04}=-\cfrac{1}{\tau a_{1}},\ \tau a^{3}a_{1}=a_{01}, \\&&a=i\xi_{0}-((i\xi_{1})^{2}+(i\xi_{2})^{2}+(i\xi_{3})^{2}),\ a_{1}=\cfrac{(i\xi_{1})^{2}+(i\xi_{2})^{2}+(i\xi_{3})^{2}}{a},
        \end{eqnarray*} we can obtain that there exist\ $\psi_{2}$,\ such that
        \begin{eqnarray*}p(M,\ t)I_{\Omega\times[0,\ T]}&=&(\psi_{2}(A_{1},\ A_{2},\ A_{3},\ A_{7},\ Z_{1})+\int_{\partial\Omega}\cfrac{\partial h_{2}(M-P_{1})}{\partial n_{p_{1}}}A_{9}(P_{1},\ t)dS_{P_{1}})I_{\Omega\times[0,\ T]},\ a.e.
        \\ \psi_{2}(A_{1},\ A_{2},\ A_{3},\ A_{7},\ Z_{1})&=&F^{-1}((B^{-1}_{1})_{4}\beta_{1}+(B^{-1}_{1}B_{2})_{4}FI(Z_{2}))-\int_{\partial\Omega}\cfrac{\partial h_{2}(M-P_{1})}{\partial n_{p_{1}}}A_{9}(P_{1},\ t)dS_{P_{1}}\\
        \\h_{2}(x,\ y,\ z)&=&\cfrac{1}{4\pi\sqrt{x^{2}+y^{2}+z^{2}}}\ .\end{eqnarray*}
If we let\ $M\rightarrow P_{0}^{+}$, from Corollary 2.1, we obtain \[ \cfrac{A_{9}(P_{0},\ t)}{2}=\psi_{3}(A_{1},\ A_{2},\ A_{3},\ A_{7},\ Z_{1})(P_{0},\ t)+\int_{\partial\Omega}\cfrac{\partial h_{2}(P_{0}-P_{1})}{\partial n_{p_{1}}}A_{9}(P_{1},\ t)dS_{P_{1}},\] where$$\psi_{3}(A_{1},\ A_{2},\ A_{3},\ A_{7},\ Z_{1})(P_{0},\ t)=\lim_{M\rightarrow P_{0}^{+}}\psi_{2}(A_{1},\ A_{2},\ A_{3},\ A_{7},\ Z_{1})(M,\ t). $$We discuss (3.82) as the following.
\begin{lemma}\label{lamma1} There is only one linear independent solution for each of two transposed equations,\begin{eqnarray} \cfrac{A_{9}}{2}&=&\cfrac{\partial h_{2}}{\partial n_{p_{1}}}.\ast_{\partial\Omega}A_{9},\\ \cfrac{A^{\prime}_{9}}{2}&=&-\cfrac{\partial h_{2}}{\partial n_{p_{0}}}.\ast_{\partial\Omega}A^{\prime}_{9}.\end{eqnarray}\end{lemma}
{\it Proof of lemma 3.7}. From the potential theory, we know\ $A_{9}=1$\ is the solution of (3.83). Hence, from Fredholm theorem, we obtain there is at least one linear independent solution\ $A^{\prime}_{9}$\ for (3.84). If there exist two linear independent solutions\ $A^{\prime}_{9,\ 1}$\ and\ $A^{\prime}_{9,\ 2}$\ for (3.84), then we consider\[ v_{1}(M)=h_{2}.\ast_{\partial\Omega}A^{\prime}_{9,\ 1},\ v_{2}(M)=h_{2}.\ast_{\partial\Omega}A^{\prime}_{9,\ 2}.\]
They are all harmonic on\ $R^{3}\setminus \partial \Omega$. From (3.84), we obtain\[ \cfrac{\partial v_{i}}{\partial n_{p_{0}}^{-}}=\cfrac{\partial h_{2}}{\partial n_{p_{0}}}.\ast_{\partial\Omega}A^{\prime}_{9,\ i}+\cfrac{A^{\prime}_{9,\ i}}{2}=0,\ i=1,\ 2.\]
From the uniqueness, we obtain\ $v_{1}$\ and\ $v_{2}$\ are all constants on\ $\overline{\Omega}$.
Therefore, we may select constants\ $C_{3,\ 1}$\ and\ $C_{3,\ 2}$, moreover,\ $C_{3,\ 1}^{2}+C_{3,\ 2}^{2}\neq 0$, such that on\ $\overline{\Omega}$, we have\ $C_{3,\ 1}v_{1}+C_{3,\ 2}v_{2}=0$.\\ Now we consider
\[v_{3}=C_{3,\ 1}v_{1}+C_{3,\ 2}v_{2}=h_{2}.\ast_{\partial\Omega}(C_{3,\ 1}A^{\prime}_{9,\ 1}+C_{3,\ 2}A^{\prime}_{9,\ 2}).\]
It is\ $0$\ on\ $\overline{\Omega}$. Moreover, it is harmonic on\ $R^{3}\setminus \overline{\Omega}$, and near to\ $0$\ uniformly at infinity. Hence, it is\ $0$\ on\ $R^{3}$. Therefore we obtain\[\cfrac{\partial v_{3}}{\partial n_{p_{0}}^{+}}-\cfrac{\partial v_{3}}{\partial n_{p_{0}}^{-}}=-(C_{3,\ 1}A^{\prime}_{9,\ 1}+C_{3,\ 2}A^{\prime}_{9,\ 2})=0.\]
This contradicts that\ $A^{\prime}_{9,\ 1}$\ and\ $A^{\prime}_{9,\ 2}$\ are linear independent. Hence, there is only one linear independent solution for (3.84). From Fredholm theorem, there is also only one linear independent solution for (3.83).\qed\\
\begin{corollary} If\ $A^{\prime}_{9}$\ is a non-zero solution of (3.84), then\ $h.\ast_{\partial\Omega}A^{\prime}_{9}$\ is a non-zero constant on\ $\overline{\Omega}$. \end{corollary}
{\it Proof of corollary 3.3}. If\ $h_{2}.\ast_{\partial\Omega}A^{\prime}_{9}$\ is\ $0$\ on\ $\overline{\Omega}$, then from the previous lemma we obtain\ $A^{\prime}_{9}=0$, which completes the proof.\qed\\ Hence, we can select\ $A^{\prime *}_{9}$\ is a non-zero solution of (3.84), such that\ $h_{2}.\ast_{\partial\Omega}A^{\prime *}_{9}=1$,\  on\ $\overline{\Omega}$.\\
From Fredholm theorem, we obtain a sufficient and necessary condition for the existence of the solution of (3.82) as follows,
\[ \int_{\partial\Omega}\psi_{3}(A_{1},\ A_{2},\ A_{3},\ A_{7},\ Z_{1})(P,\ t)A^{\prime *}_{9}(P)dS_{P}=0,\ \forall\ t\in [0,\ T].\]
Under this condition, again from Fredholm theorem, there exists an entire analytic function\ $\Gamma_{4}$\ which is only related to\ $\partial h_{2}(P_{0})/\partial n_{p_{1}}$, such that
\[ \cfrac{A_{9}(P_{0},\ t)}{2}=\psi_{3}(A_{1},\ A_{2},\ A_{3},\ A_{7},\ Z_{1})+\int_{\partial\Omega}\Gamma_{4}(P_{0},\ P_{1})[\psi_{3}(A_{1},\ A_{2},\ A_{3},\ A_{7},\ Z_{1})(P_{1},\ t)]dS_{P_{1}}+C_{0}^{*},\]where\ $C_{0}^{*}$\ is a constant.\\ Maybe you will point out\ $\partial h_{2}(P_{0})/\partial n_{p_{1}}(P_{0}-P_{1})$\ is not continuous, if\ $P_{0}=P_{1}$. Yes, that's true. It's lucky that $$|P_{0}-P_{1}|^{(\epsilon_{0}+5/2)}\cfrac{\partial h_{2}(P_{0})}{\partial n_{p_{1}}}(P_{0}-P_{1})$$
is continuous, if\ $\epsilon_{0}\in(0,\ 0.5]$. Hence\ $\partial h_{2}(P_{0})/\partial n_{p_{1}}(P_{0}-P_{1})$\ is weak singular kernel. From Fredholm theorem, we know (3.90) still stands.\\
It looks that the classical solution of Eqs(3.22) would be locally exist and unique. So were the Eq(3.1) and Eqs(3.2). But\ $F^{-1}[f_{1}]$\ will cause some trouble.\\
          We can work out$$F^{-1}[f_{1}]=\int_{0}^{T}d\tau\int_{\partial\Omega}\delta(x-x_{1},\ y-y_{1},\ z-z_{1},\ t-\tau)A_{3}(x_{1},\ y_{1},\ z_{1},\ \tau)n_{1}(x_{1},\ y_{1},\ z_{1})dS,$$where\ $\delta$\ is the Dirac function. We don't know what the next is.\\ Because\ $\partial\Omega\in C^{1,\ \beta}$, we can get the following from Theorem 1.1, $$\partial\Omega=\bigcup_{k=1}^{N}\partial\Omega_{k},$$ where\ $\partial\Omega_{k}$\ is the graph of a\ $C^{1,\ \beta}$\ function of two of the coordinates\ $x,\ y,\ z$.\\ Without loss of the generality, we assume\ $\partial\Omega_{k}$\ is the graph of a\ $C^{1,\ \beta}$\ function\ $z=f_{k}(x,\ y),\ (x,\ y)\in D_{k}\subset R^{2}$, then we may solve the inverse Fourier transform of\ $f_{1}$\ on\ $\partial\Omega_{k}$\ as follows,\begin{eqnarray*}&&F^{-1}[\int_{0}^{T}d\tau\int_{\partial\Omega_{k}}A_{3}(x_{1},\ y_{1},\ z_{1},\ \tau)n_{1}(x_{1},\ y_{1},\ z_{1})e^{-i\xi_{0}\tau-i\xi_{1}x_{1}-i\xi_{2}y_{1}-i\xi_{3}z_{1}}dS]\\
          &=&F^{-1}[\int_{0}^{T}d\tau\int_{D_{k}}A_{3}(x_{1},\ y_{1},\ f_{k}(x_{1},\ y_{1}),\ \tau)n_{1}(x_{1},\ y_{1},\ f_{k}(x_{1},\ y_{1}))e^{-i\xi_{0}\tau-i\xi_{1}x_{1}-i\xi_{2}y_{1}-i\xi_{3}f_{k}(x_{1},\ y_{1})}\\&&\sqrt{1+f^{2}_{k1}(x_{1},\ y_{1})+f^{2}_{k2}(x_{1},\ y_{1})}dx_{1}dy_{1}]\end{eqnarray*}\begin{eqnarray*}&=&A_{3}(x,\ y,\ f_{k}(x,\ y),\ t)n_{1}(x,\ y,\ f_{k}(x,\ y))\sqrt{1+f^{2}_{k1}(x,\ y)+f^{2}_{k2}(x,\ y)}I_{D_{k}\times[0,\ T]}\\&&\cfrac{1}{2\pi}\int_{-\infty}^{+\infty}e^{-i\xi_{3}f_{k}(x,\ y)}e^{i\xi_{3}z}d\xi_{3}\\&=&A_{3}(x,\ y,\ f_{k}(x,\ y),\ t)n_{1}(x,\ y,\ f_{k}(x,\ y))\sqrt{1+f^{2}_{k1}(x,\ y)+f^{2}_{k2}(x,\ y)}I_{D_{k}\times[0,\ T]}\delta[z-f_{k}(x,\ y)],\end{eqnarray*}
          where$$f_{k1}(x,\ y)=\cfrac{\partial f_{k}(x,\ y)}{\partial x},\ f_{k2}(x,\ y)=\cfrac{\partial f_{k}(x,\ y)}{\partial y}.$$ We are surprised to see that the inverse Fourier transform of\ $f_{1}$\ on\ $\partial\Omega_{k}$\ is related to\ $\delta[z-f_{k}(x,\ y)]$. This means that\ $F^{-1}[f_{1}]$\ will be related to\ $\delta(\partial\Omega)$.\\
       We regret that we lost local existence and uniqueness for the solution of\ $Z_{1}=T_{0}(Z_{1})$. However, it would be too easy if there were no\ $\delta(\partial\Omega)$. In the next section, we will discuss\ $Z_{1}=T_{0}(Z_{1})$\ by Lerry-Schauder degree and the Sobolev space\ $H^{-m_{1}}(\Omega_{1})$, where\ $\Omega_{1}=\Omega\times(0,\ T)$.\\
\section{Existence} \setcounter{equation}{0}
In this section, we will discuss the existence for classical solutions of\ $Z_{1}=T_{0}(Z_{1})$. We include remarks if it is necessary. First of all, we introduce our ideas as follows.\\
First, we will construct a norm\ $\|\cdot\|_{-m_{1}}$\ for\ $T_{0}(Z_{1})$.\\ Second, since\ $Z_{1}=T_{0}(Z_{1})$\ are generalized integral equations, we make approximate ordinary integral equations\ $Z_{1}=T_{0\epsilon}(Z_{1}),\ \forall \epsilon>0$,
such that\ $\forall Z_{1}\in C(\overline{\Omega}\times[0,\ T]),\ \|Z_{1}\|_{\infty}\leq M,\ M$\ is given, we have the following,
\[ \lim_{\epsilon\rightarrow 0}\|T_{0\epsilon}(Z_{1})-T_{0}(Z_{1})\|_{-m_{1}}=0,\ \mbox{uniformly,}\]where\[\|Z_{1}\|_{\infty}=\max_{1\leq i\leq 33}\|Z_{1,\ i}\|_{\infty},\ \|Z_{1,\ i}\|_{\infty}=\max_{X=(x,\ y,\ z,\ t)^{T}\in\overline{\Omega}\times[0,\ T]}|Z_{1,\ i}(X)|,\] $Z_{1,\ i},\ 1\leq i\leq 33$, are components of\ $Z_{1}$.\\ This will help us to prove\ $T_{0}(Z_{1k})$\ is sequentially compact under norm\ $\|\cdot\|_{-m_{1}}$, if\ $\forall Z_{1k}\in C(\overline{\Omega}\times[0,\ T]),\ \|Z_{1k}\|_{\infty}\leq M,\ k\geq 1,\ M$\ is given.\\Finally, we use some primary theorems on the Leray-Schauder degree to discuss ordinary integral equations\ $Z_{1}=T_{0\epsilon}(Z_{1}),\ \forall \epsilon>0$. If\ $Z_{1\epsilon}$\ satisfies\ $Z_{1}=T_{0\epsilon}(Z_{1}),\ \forall \epsilon>0$,\ and bounded uniformly, then there exists a sequence\ $Z_{1\epsilon_{k}},\ k\geq1,\ \epsilon_{k}\rightarrow 0$, if\ $k\rightarrow+\infty$, such that\ $Z_{1\epsilon_{k}}$\ is convergent to\ $Z_{1}^{\ast}$. We will obtain\ $Z_{1}^{\ast}=T_{0}(Z_{1}^{\ast})$, which is what we want.\\Can our imagination come true? Let's introduce our answer. The answer is not unique.
\begin{definition}\label{definition}$\forall Z_{1}=(Z_{1,\ i})_{33\times 1},\ Z_{1,\ i}\in H^{-m_{1}}(\Omega_{1}),\ 1\leq i\leq 33,\ \Omega_{1}=\Omega\times(0,\ T)$, we have
\[\|Z_{1}\|_{-m_{1}}=\max_{1\leq i\leq 33}\|Z_{1,\ i}\|_{-m_{1}},\ \|Z_{1,\ i}\|_{-m_{1}}=\sup_{\varphi\in C_{0}^{\infty}(\Omega_{1})}\cfrac{|<Z_{1,\ i},\ \varphi>|}{\|\varphi\|_{m_{1}}},\] where\ $m_{1}=6+2c,\ c=\max\{\partial(b_{1}B^{-1}_{1}),\ \partial(b_{1}B^{-1}_{1}B_{2})\}$, here\ $\partial(\cdot)$\ means the highest degree,\ $X=(x,\ y,\ z,\ t)^{T}$,\ $<Z_{1,\ i},\ \varphi>$\ is the value of the generalized function\ $Z_{1,\ i}$\ on\ $\varphi$, if\ $Z_{1,\ i}$\ is locally integrable, then \begin{eqnarray}<Z_{1,\ i},\ \varphi>&=&\int_{\Omega_{1}}Z_{1,\ i}(X)\varphi(X) dX,\ 1\leq i\leq 33,\\ \|\varphi\|_{m_{1}}&=&(\int_{\Omega_{1}}\sum_{|\alpha|\leq m_{1}}|\partial^{\alpha}\varphi(X)|^{2}dX)^{1/2}.\end{eqnarray}\end{definition}
Next we can get that\[\|T_{0}(Z_{1})\|_{-m_{1}}=\max_{1\leq i\leq 33}\sup_{\varphi\in C_{0}^{\infty}(\Omega_{1})}\cfrac{|<T_{0,\ i}(Z_{1}),\ \varphi>|}{\|\varphi\|_{m_{1}}}<+\infty,\]in the following lemma, where\ $ T_{0,\ i}(Z_{1}),\ 1\leq i\leq 33$, are components of\ $T_{0}(Z_{1})$.\\Now we see approximate ordinary integral equations,\ $\forall \epsilon>0$,\ $\forall \epsilon>0$,
\[ T_{0\epsilon}(Z_{1})=w_{1\epsilon}+w_{2\epsilon}.\ast(\psi(Z_{1}I_{\overline{\Omega}\times[0,\ T]})),\] where$$w_{1\epsilon}=F^{-1}(\widetilde{\delta_{\epsilon}}B_{1}^{-1}\beta_{1}),\ w_{2\epsilon}=F^{-1}(\widetilde{\delta_{\epsilon}}B_{1}^{-1}B_{2}),\ \widetilde{\delta_{\epsilon}}=F(\delta_{\epsilon}),$$
       $$\delta_{\epsilon}=\cfrac{1}{(\sqrt{\pi\epsilon})^{4}}e^{-|X|^{2}/\epsilon},\ |X|=\sqrt{x^{2}+y^{2}+z^{2}+t^{2}},\ \lim_{\epsilon\rightarrow 0}\delta_{\epsilon}=\delta(X),\ \lim_{\epsilon\rightarrow 0}\widetilde{\delta_{\epsilon}}=1,  $$ $\delta(X)$\ is the Dirac function.\\
        From the previous section, we can obtain\ $A_{2}$\ as follows,\[ A_{2}I_{\Omega}=[h.\ast_{\overline{\Omega}}A_{1}+h.\ast_{\partial\Omega}A_{7}+\cfrac{\partial h}{\partial n_{p_{1}}}.\ast_{\partial\Omega}A_{3}+h.\ast \psi_{1}(Z_{1})]I_{\Omega}|_{t=T}.\]In the Dirichlet problem,\ $A_{1},\ A_{3}$\ are known, and from (3.47) we can get\ $A_{7}$\ as follows,\begin{eqnarray}\cfrac{A_{7}(P_{0},\ t)}{2}&=&g_{02}(P_{0},\ t)+\cfrac{\partial h(P_{0},\ t)}{\partial n_{p_{0}}}.\ast \psi_{1}(Z_{1})+\nonumber\\&&\int_{0}^{T}d\tau\int_{\partial\Omega}\Gamma_{1}(P_{0},\ t,\ P_{1},\ \tau)[g_{02}(P_{1},\ \tau)+\cfrac{\partial h(P_{1},\ \tau)}{\partial n_{p_{1}}}.\ast \psi_{1}(Z_{1})]dS_{P_{1}}.\end{eqnarray}In the Neumann problem,\ $A_{1},\ A_{7}$\ are known, and from (3.61) we can get\ $A_{3}$\ as follows,
        \begin{eqnarray} \cfrac{A_{3}(P_{0},\ t)}{2}&=&g_{04}(P_{0},\ t)+h(P_{0},\ t).\ast \psi_{1}(Z_{1})+\nonumber\\&&\int_{0}^{T}d\tau\int_{\partial\Omega}\Gamma_{2}(P_{0},\ t,\ P_{1},\ \tau)[g_{04}(P_{1},\ \tau)+ h(P_{1},\ \tau).\ast \psi_{1}(Z_{1})]dS_{P_{1}}.\end{eqnarray}In the Robin problem,\ $A_{1},\ A_{8}$\ are known, and from (3.74) we can get\ $A_{3}$\ as follows,\begin{eqnarray} \cfrac{A_{3}(P_{0},\ t)}{2}&=&g_{06}(P_{0},\ t)+ h(P_{0},\ t).\ast \psi_{1}(Z_{1})+\nonumber\\&&\int_{0}^{T}d\tau\int_{\partial\Omega}\Gamma_{3}(P_{0},\ t,\ P_{1},\ \tau)[g_{06}(P_{1},\ \tau)+h(P_{1},\ \tau).\ast \psi_{1}(Z_{1})]dS_{P_{1}}.\end{eqnarray}And from (3.79), (3.80), (3.81), we can also get\ $A_{4},\ A_{5},\ A_{6}$\ as follows,
         \begin{eqnarray}A_{4}&=&u_{x}|_{\partial\Omega}=\cfrac{\partial A_{3}}{\partial x}+(A_{7}-\cfrac{\partial A_{3}}{\partial n})n_{1},\\A_{5}&=&u_{y}|_{\partial\Omega}=\cfrac{\partial A_{3}}{\partial y}+(A_{7}-\cfrac{\partial A_{3}}{\partial n})n_{2},\\A_{6}&=&u_{z}|_{\partial\Omega}=\cfrac{\partial A_{3}}{\partial z}+(A_{7}-\cfrac{\partial A_{3}}{\partial n})n_{3},\end{eqnarray}
where\ $n=(n_{1},\ n_{2},\ n_{3})^{T}$\ is the exterior normal vector to\ $\partial\Omega$.\\Finally, from (3.90), we can get\ $A_{9}$\ as follows,
\[ \cfrac{A_{9}(P_{0},\ t)}{2}=\psi_{3}(A_{1},\ A_{2},\ A_{3},\ A_{7},\ Z_{1})+\int_{\partial\Omega}\Gamma_{4}(P_{0},\ P_{1})[\psi_{3}(A_{1},\ A_{2},\ A_{3},\ A_{7},\ Z_{1})(P_{1},\ t)]dS_{P_{1}}+C_{0}^{*},\]where\ $C_{0}^{*}$\ is a constant.
\\ We can see that\ $A_{4},\ A_{5},\ A_{6},\ A_{9}$\ are all compact operators on\ $Z_{1}$.\\
We know\ $v_{1}\ast v_{2}\in S^{\prime}$\ and\ $F(v_{1}\ast v_{2})=F(v_{1})F(v_{2})$\ will still hold if\ $v_{1}\in S,\ v_{2}\in S^{\prime}$. There is an example of this on pages 118 to 119 of [10]. \\ We can see that\ $\forall \epsilon>0$,\ $\delta_{\epsilon}\in S,\ F^{-1}(b_{1}^{-1})\in S^{\prime},\ F^{-1}(f_{j})\in S^{\prime},\ F^{-1}(f_{jk})\in S^{\prime},\ 1\leq j,\ k\leq 3$. Hence we can get the following,\begin{eqnarray*} F^{-1}[\widetilde{\delta_{\epsilon}}(i\xi)^{\alpha}b_{1}^{-1}]&=&(\partial^{\alpha}\delta_{\epsilon}).\ast F^{-1}(b_{1}^{-1}),\\ F^{-1}[\widetilde{\delta_{\epsilon}}(i\xi)^{\alpha}f_{j}]&=&(\partial^{\alpha}\delta_{\epsilon}).\ast_{\partial\Omega}(A_{3}n_{j}),\\ F^{-1}(\widetilde{\delta_{\epsilon}}f_{jk})&=&\delta_{\epsilon}.\ast_{\partial\Omega}(A_{j+3}n_{k}),\end{eqnarray*}where\ $1\leq j,\ k\leq 3.$\\
We may denote that\ $T_{0\epsilon}(Z_{1})=\delta_{\epsilon}.\ast T_{0}(Z_{1})$. From\ $\partial^{\alpha}\delta_{\epsilon}\in S$, we can get$$ (\partial^{\alpha}\delta_{\epsilon}).\ast F^{-1}(b_{1}^{-1})\in C(R^{4}),\ (\partial^{\alpha}\delta_{\epsilon}.\ast F^{-1}(b_{1}^{-1})).\ast(\psi(Z_{1}I_{\overline{\Omega}\times[0,\ T]}))\in C(\overline{\Omega}\times[0,\ T]).$$ We see that\ $\forall \epsilon>0$, there is no\ $\delta(\partial\Omega)$\ or\ $\partial F^{-1}(b_{1}^{-1})$\ in\ $T_{0\epsilon}$\ again. So\ $Z_{1}I_{\overline{\Omega}\times[0,\ T]}=T_{0\epsilon}(Z_{1}I_{\overline{\Omega}\times[0,\ T]})$\ are ordinary integral equations. Moreover\ $T_{0\epsilon}(Z_{1})$\ is bounded uniformly and equicontinuous if\ $Z_{1}\in C(\overline{\Omega}\times[0,\ T])$\ and bounded uniformly. From the Arzela-Ascoli theorem, we know\ $T_{0\epsilon}$\ is a compact operator on\ $Z_{1}$. And the Leray-Schauder degree can work now. \\For the preliminaries, we have a lemma as follows. \begin{lemma} \label{lemma1}(1)There exists\ $C>0$, such that\[ \|Z_{1}\|_{-m_{1}}\leq C \|Z_{1}\|_{\infty},\ \forall Z_{1}\in C(\overline{\Omega}\times[0,\ T]).\]
        (2) $\forall Z_{1}\in C(\overline{\Omega}\times[0,\ T]),\ \|Z_{1}\|_{\infty}\leq M,\ M>0,\ M$\ is given, we have
\[ \|T_{0}(Z_{1})\|_{-m_{1}}<+\infty,\ \mbox{and}\ \lim_{\epsilon\rightarrow 0}\|T_{0\epsilon}(Z_{1})-T_{0}(Z_{1})\|_{-m_{1}}=0,\ \mbox{uniformly.}\]
(3) $T_{0}(Z_{1k}),\ k\geq1,$\ is sequentially compact under norm\ $\|\cdot\|_{-m_{1}}$, if\ $Z_{1k}\in C(\overline{\Omega}\times[0,\ T]),\ \|Z_{1k}\|_{\infty}\leq M,\ k\geq 1,\ M>0,\ M$\ is given.\\(4)$\forall \epsilon_{0}>0,\ \forall Z_{1}\in C(\overline{\Omega}\times[0,\ T]),\ \|Z_{1}\|_{\infty}\leq M,\ M>0,\ M$\ is given, we have\[ \lim_{\epsilon\rightarrow\epsilon_{0}}\|[T_{0\epsilon}(Z_{1})-T_{0\epsilon_{0}}(Z_{1})]I_{\overline{\Omega}\times[0,\ T]}\|_{\infty}=0.\]\end{lemma}
{\it Proof of lemma 4.1}. (1)From\ $\|Z_{1,\ i}\|_{-m_{1}}\leq \|Z_{1,\ i}\|_{L^{2}}\leq \|Z_{1,\ i}\|_{\infty}\sqrt{m(\overline{\Omega})T},\ 1\leq i\leq 33$, where\ $m(\overline{\Omega})$\ is the Lebesgue measure of\ $\overline{\Omega}$. Hence we get\ $\|Z_{1}\|_{-m_{1}}\leq \|Z_{1}\|_{\infty}\sqrt{m(\overline{\Omega})T}$. We may let\ $C=\sqrt{m(\overline{\Omega})T}$.\\
There will not exist\ $C>0$, such that\[ \|Z_{1}\|_{\infty}\leq C \|Z_{1}\|_{-m_{1}},\ \forall Z_{1}\in C(\overline{\Omega}\times[0,\ T]).\]
We may select\ $ Z_{1k}\in C(\overline{\Omega}\times[0,\ T]),\ k\geq1,\ \|Z_{1k}\|_{\infty}\equiv1$, but\ $Z_{1k}\rightarrow 0,\ \mbox{a.e.}$, here a.e. means that almost everywhere. Then from the Lebesgue dominated convergence theorem,\ $\|Z_{1k}\|_{-m_{1}}\rightarrow 0$. Hence (4.19) will not stand.\ $\|Z_{1}\|_{-m_{1}}$\ and\ $\|Z_{1}\|_{\infty}$\ are not equivalent.\\
(2)At first, we prove\ $(1+|\xi|^{2})^{-c-3}B_{1}^{-1}\in L^{1}(R^{4})$, where\ $c=\max\{\partial(b_{1}B^{-1}_{1}),\ \partial(b_{1}B^{-1}_{1}B_{2})\}$, here\ $\partial(\cdot)$\ means the highest degree.\\ We assume\ $B_{1}^{-1}=(b_{j,\ k,\ 1})_{33\times 33}$. From the previous section, we can obtain that\ $b_{1}B_{1}^{-1}$\ is a polynomial matrix. Moreover, the least degree of\ $b_{1}B_{1}^{-1}$\ is\ $2$. This means that there exists a constant\ $C_{4,\ 0}>0$, such that\ $|b_{1}b_{j,\ k,\ 1}|\leq C_{4,\ 0}|\xi|^{2},\ 1\leq j,\ k\leq 33$, if\ $|\xi|\leq 1$.\\
We can work out the following,\begin{eqnarray*}|a|^{2}&=&\xi_{0}^{2}+(\xi_{1}^{2}+\xi_{2}^{2}+\xi_{3}^{2})^{2}
\\ &=&|\xi|^{2}\cos^{2}\theta+|\xi|^{4}\sin^{4}\theta
\\&=&(|\xi|^{2}-|\xi|^{4})\cos^{2}\theta+|\xi|^{4}[(\sin^{2}\theta-1/2)^{2}+3/4]
\\ &=&(|\xi|^{4}-|\xi|^{2})\sin^{4}\theta+|\xi|^{2}[(\sin^{2}\theta-1/2)^{2}+3/4],
\end{eqnarray*}where\ $\cos\theta=\xi_{0}/|\xi|$. So we obtain that\ $|a|\geq \sqrt{3}|\xi|^{2}/2$,\ if\ $|\xi|\leq 1$,\ $|a|\geq \sqrt{3}|\xi|/2$,\ if\ $|\xi|> 1$. \\
Hence, there exists a constant\ $C_{4,\ 1}>0$, such that\ $|b_{1}b_{j,\ k,\ 1}/a|\leq C_{4,\ 1},\ 1\leq j,\ k\leq 33$, if\ $|\xi|\leq 1$.\\By spherical coordinates transformation on\ $R^{4}$, $$\begin{cases} \xi_{1}=\rho \sin \theta \sin \theta_{1} \sin \theta_{2},\\
\xi_{2}=\rho \sin \theta \sin \theta_{1} \cos \theta_{2},\\
\xi_{3}=\rho \sin \theta \cos \theta_{1},\\
\xi_{0}=\rho \cos \theta,\ \theta,\ \theta_{1}\in[0,\ \pi],\ \theta_{2}\in[0,\ 2\pi],
\end{cases},$$ and\ $d\xi_{1}d\xi_{2}d\xi_{3}d\xi_{0}=\rho^{3} \sin^{2}\theta \sin \theta_{1}d\rho d\theta d\theta_{1} d\theta_{2}$, we can obtain the following,
\begin{eqnarray*}\int_{|\xi|\leq 1}(1+|\xi|^{2})^{-c-3}|b_{j,\ k,\ 1}|d\xi&=&\int_{|\xi|\leq 1}(1+|\xi|^{2})^{-c-3}|(b_{1}b_{j,\ k,\ 1}/a)|/|(aa_{1})|d\xi
\\(|aa_{1}|=\rho^{2} \sin^{2}\theta)&\leq&\int_{0}^{2\pi}d\theta_{2}\int_{0}^{\pi}d\theta_{1}\int_{0}^{\pi}d\theta\int_{0}^{1}C_{4,\ 1} \rho \sin \theta_{1}d\rho<+\infty,\end{eqnarray*}$1\leq j,\ k\leq 33$. \\ Because\ \ $c=\max\{\partial(b_{1}B^{-1}_{1}),\ \partial(b_{1}B^{-1}_{1}B_{2})\}$, there exists a constant\ $C_{4,\ 2}>0$, such that\ $|b_{1}b_{j,\ k,\ 1}|\leq C_{4,\ 2}(1+|\xi|^{2})^{c},\ 1\leq j,\ k\leq 33$. Hence, we have
\begin{eqnarray*}\int_{|\xi|> 1}(1+|\xi|^{2})^{-c-3}|b_{j,\ k,\ 1}|d\xi&=&\int_{|\xi|> 1}(1+|\xi|^{2})^{-c-3}|(b_{1}b_{j,\ k,\ 1}/a)|/|(aa_{1})|d\xi
\\(|aa_{1}|=\rho^{2} \sin^{2}\theta)&\leq&\int_{0}^{2\pi}d\theta_{2}\int_{0}^{\pi}d\theta_{1}\int_{0}^{\pi}d\theta\int_{1}^{+\infty}2C_{4,\ 2} (1+\rho^{2})^{-3} \sin \theta_{1}d\rho<+\infty,\end{eqnarray*}$1\leq j,\ k\leq 33$.\\So\ $(1+|\xi|^{2})^{-c-3}B_{1}^{-1}\in L^{1}(R^{4})$.\\ Next we will prove\ $(F[T_{0}(Z_{1})])F^{-1}(\varphi)\in L^{1}(R^{4})$. We can work out the following,
\begin{eqnarray*}<T_{0}(Z_{1}),\ \varphi>&=&<F[T_{0}(Z_{1})],\ F^{-1}(\varphi)>\\ &=&<(B^{-1}_{1}\beta_{1}+B^{-1}_{1}B_{2}FI(Z_{2})),\ F^{-1}(\varphi)>\\&=&<(B^{-1}_{1}\beta_{1}+B^{-1}_{1}B_{2}FI(Z_{2}))(1+|\xi|^{2})^{-c-3},\\&& (1+|\xi|^{2})^{c+3}F^{-1}(\varphi)>.
\end{eqnarray*}From\ $\psi(Z_{1})$\ is a continuous, we get that\ $\beta_{1},\ Z_{2}$\ are all bounded. \\ We have the following,\begin{eqnarray*}|(i\xi)^{\alpha}F^{-1}(\varphi)|/\|\varphi\|_{m_{1}}&\leq&\int_{\Omega_{1}}|\partial^{\alpha}\varphi|dX/\|\varphi\|_{m_{1}}\\
&\leq&(\int_{\Omega_{1}}|\partial^{\alpha}\varphi|^{2}dX)^{1/2}(\int_{\Omega_{1}}dX)^{1/2}/\|\varphi\|_{m_{1}}\\&\leq&\sqrt{m(\overline{\Omega})T},\end{eqnarray*}
where\ $|\alpha|\leq 6+2c,\ c=\max\{\partial(b_{1}B^{-1}_{1}),\ \partial(b_{1}B^{-1}_{1}B_{2})\},\ m_{1}=6+2c$.\\
Hence, we can obtain\ $(F[T_{0}(Z_{1})])F^{-1}(\varphi)\in L^{1}(R^{4})$, moreover \[\|T_{0}(Z_{1})\|_{-m_{1}}\leq C_{4,\ 3},\]
where\ $C_{4,\ 3}>0$,\ is a constant only related to\ $M$. It is not related with\ $Z_{1}$.\\
We can work out as follows,
\begin{eqnarray*}<T_{0\epsilon}(Z_{1})-T_{0}(Z_{1}),\ \varphi>&=&<F[T_{0\epsilon}(Z_{1})-T_{0}(Z_{1})],\ F^{-1}(\varphi)>\\ &=&<(\widetilde{\delta_{\epsilon}}-1)F[T_{0}(Z_{1})],\ F^{-1}(\varphi)>,\\ &=&<(\widetilde{\delta_{\epsilon}}-1)(B^{-1}_{1}\beta_{1}+B^{-1}_{1}B_{2}FI(Z_{2})),\ F^{-1}(\varphi)>\\&=&<(\widetilde{\delta_{\epsilon}}-1)(B^{-1}_{1}\beta_{1}+B^{-1}_{1}B_{2}FI(Z_{2}))(1+|\xi|^{2})^{-c-3},\\&& (1+|\xi|^{2})^{c+3}F^{-1}(\varphi)>.
\end{eqnarray*}
Hence, we can get\[\|T_{0\epsilon}(Z_{1})-T_{0}(Z_{1})\|_{-m_{1}}\leq
\int_{R^{4}}|\widetilde{\delta_{\epsilon}}-1||(B^{-1}_{1}\beta_{1}+B^{-1}_{1}B_{2}FI(Z_{2}))(1+|\xi|^{2})^{-c-3}|\sqrt{m(\overline{\Omega})T}d\xi.\]
We can also work out the following,
\begin{eqnarray*}\widetilde{\delta_{\epsilon}}&=&\int_{R^{4}}\cfrac{1}{(\sqrt{\pi\epsilon})^{4}}e^{-|X|^{2}/\epsilon}e^{-i\xi\cdot X}dX,\ (X=\sqrt{\epsilon}Y,\ dX=(\sqrt{\epsilon})^{4}dY)
\\&=&\int_{R^{4}}\cfrac{1}{\pi^{2}}e^{-|Y|^{2}}e^{-i\xi\cdot \sqrt{\epsilon}Y}dY,
\\ &=&\int_{R^{4}}\cfrac{1}{\pi^{2}}e^{-|Y+i\sqrt{\epsilon}\xi/2 |^{2}}e^{-\epsilon|\xi|^{2}/4}dY,\\ &=&\int_{R^{4}}\cfrac{1}{\pi^{2}}e^{-|Y |^{2}}e^{-\epsilon|\xi|^{2}/4}dY=e^{-\epsilon|\xi|^{2}/4}.\end{eqnarray*}Here we use Cauchy contour integral as follows,
\[\int_{-\infty}^{+\infty} e^{-(x+i\xi_{1})^{2}}dx=\int_{-\infty}^{+\infty} e^{-x^{2}}dx,\ \forall \xi_{1}\in R.\]
We can get the following,
$$ \int_{R^{4}}|\widetilde{\delta_{\epsilon}}-1||(B^{-1}_{1}\beta_{1}+B^{-1}_{1}B_{2}FI(Z_{2}))(1+|\xi|^{2})^{-c-3}|\sqrt{m(\overline{\Omega})T}d\xi=I_{3,\ 1}+I_{3,\ 2},$$ where\begin{eqnarray*} I_{3,\ 1}&=&\int_{|\xi|>M_{0}}|\widetilde{\delta_{\epsilon}}-1||(B^{-1}_{1}\beta_{1}+B^{-1}_{1}B_{2}FI(Z_{2}))(1+|\xi|^{2})^{-c-3}|\sqrt{m(\overline{\Omega})T}d\xi,\\ I_{3,\ 2}&=&\int_{|\xi|\leq M_{0}}|\widetilde{\delta_{\epsilon}}-1||(B^{-1}_{1}\beta_{1}+B^{-1}_{1}B_{2}FI(Z_{2}))(1+|\xi|^{2})^{-c-3}|\sqrt{m(\overline{\Omega})T}d\xi.\end{eqnarray*}
$\forall \epsilon^{\prime}>0$, there exists\ $M_{0}>0$, which is only related to\ $M$,
such that\ $|I_{3,\ 1}|\leq \epsilon^{\prime}/2$. And for such\ $M_{0}$, there exists\ $\delta_{0}>0$, such that\ $|I_{3,\ 2}|\leq \epsilon^{\prime}/2$, if\ $\epsilon\leq \delta_{0}$.\\
Hence, we obtain\ $\|T_{0\epsilon}(Z_{1})-T_{0}(Z_{1})\|_{-m_{1}}\rightarrow0,$\ uniformly.\\
(3)We will prove\ $T_{0}(Z_{1k}),\ k\geq 1,$\ is totally bounded. From (2),\ $\forall \epsilon_{1}>0,\ \exists \delta_{1}>0$, such that \[\|T_{0\epsilon}(Z_{1k})-T_{0}(Z_{1k})\|_{-m_{1}}\leq \epsilon_{1}/3,\ \forall \epsilon\in (0,\ \delta_{1}),\ \forall k,\ k\geq 1.\]  If we choose\ $\epsilon_{0}\in (0,\ \delta_{1})$, then\ $T_{0\epsilon_{0}}(Z_{1k}),\ k\geq 1$\ is sequentially compact. There exist finite\ $\epsilon_{1}/(3C)$\ net\ $T_{0\epsilon_{0}}(Z_{1k_{1}}),\ T_{0\epsilon_{0}}(Z_{1k_{2}}),\cdots,\ T_{0\epsilon_{0}}(Z_{1k_{s}})$, where\ $C$\ is defined in (1).\\ This means that\ $\forall k,\ \exists l,\ 1\leq l\leq s$, such that\ $\|T_{0\epsilon_{0}}(Z_{1k})-T_{0\epsilon_{0}}(Z_{1k_{l}})\|_{\infty}\leq \epsilon_{1}/(3C)$.\\ From\ $\|T_{0\epsilon_{0}}(Z_{1k})-T_{0}(Z_{1k})\|_{-m_{1}}\leq \epsilon_{1}/3,\ \|T_{0\epsilon_{0}}(Z_{1k})-T_{0\epsilon_{0}}(Z_{1k_{l}})\|_{-m_{1}}\leq \epsilon_{1}/3$,\\ $\|T_{0\epsilon_{0}}(Z_{1k_{l}})-T_{0}(Z_{1k_{l}})\|_{-m_{1}}\leq \epsilon_{1}/3$, we obtain\ $\|T_{0}(Z_{1k})-T_{0}(Z_{1k_{l}})\|_{-m_{1}}\leq \epsilon_{1}$. \\So\ $T_{0}(Z_{1k_{1}}),\ T_{0}(Z_{1k_{2}}),\cdots,\ T_{0}(Z_{1k_{s}})$\ is a finite\ $\epsilon_{1}$\ net for\ $T_{0}(Z_{1k}),\ k\geq 1$. \\This means that\ $T_{0}(Z_{1k}),\ k\geq 1,$\ is totally bounded. From Hausdorff theorem on page 14 in [7],\ $T_{0}(Z_{1k}),\ k\geq1,$\ is sequentially compact under norm\ $\|\cdot\|_{-m_{1}}$.\\
(4)From\ $T_{0\epsilon,\ i}(Z_{1})=\delta_{\epsilon}.\ast T_{0,\ i}(Z_{1}),\ 1\leq i\leq 33$, we obtain
\begin{eqnarray*}T_{0\epsilon,\ i}(Z_{1})-T_{0\epsilon_{0},\ i}(Z_{1})&=&(\delta_{\epsilon}-\delta_{\epsilon_{0}}).\ast T_{0,\ i}(Z_{1})
\\&=&F^{-1}(F((\delta_{\epsilon}-\delta_{\epsilon_{0}}).\ast T_{0,\ i}(Z_{1})))\\&=&\cfrac{1}{(2\pi)^{4}}\int_{R^{4}}(e^{-\epsilon|\xi|^{2}/4}-e^{-\epsilon_{0}|\xi|^{2}/4})F(T_{0,\ i}(Z_{1}))e^{i\xi\cdot X}d\xi
\\&=& I_{3,\ 3}+I_{3,\ 4},\end{eqnarray*}
where\begin{eqnarray*} I_{3,\ 3}&=&\cfrac{1}{(2\pi)^{4}}\int_{|\xi|>M_{0}}(e^{-\epsilon|\xi|^{2}/4}-e^{-\epsilon_{0}|\xi|^{2}/4})(1+|\xi|^{2})^{c+3}(F(T_{0,\ i}(Z_{1}))(1+|\xi|^{2})^{-c-3})e^{i\xi\cdot X}d\xi,\\ I_{3,\ 4}&=&\cfrac{1}{(2\pi)^{4}}\int_{|\xi|\leq M_{0}}(e^{-\epsilon|\xi|^{2}/4}-e^{-\epsilon_{0}|\xi|^{2}/4})(1+|\xi|^{2})^{c+3}(F(T_{0,\ i}(Z_{1}))(1+|\xi|^{2})^{-c-3})e^{i\xi\cdot X}d\xi.\end{eqnarray*}From (2), we know\ $F(T_{0,\ i}(Z_{1}))(1+|\xi|^{2})^{-c-3}\in L^{1}(R^{4}),\ 1\leq i\leq 33$.\\ If we let\ $\epsilon\in [\epsilon_{0}/2,\ 3\epsilon_{0}/2]$, then\ $\forall \epsilon^{\prime}>0$, there exists\ $M_{0}>0$, which is related to\ $\epsilon_{0}$,
such that\ $|I_{3,\ 3}|\leq \epsilon^{\prime}/2$. And for such\ $M_{0}$, there exists\ $\delta_{0}\in (0,\ \epsilon_{0}/2)$, such that\ $|I_{3,\ 4}|\leq \epsilon^{\prime}/2$, if\ $|\epsilon-\epsilon_{0}|\leq \delta_{0}$. So the statement holds.
\qed\\
We will use result (4) into Leray-Schauder degree. This is the reason why we choose\ $\delta_{\epsilon}$\ instead of\ $I_{\{|X|\leq \epsilon\}}/ |I_{\{|X|\leq \epsilon\}}|$, where \[|I_{\{|X|\leq \epsilon\}}|=\int_{|X|\leq \epsilon}dX.\] The latter will not satisfy (4). Maybe it works after being polished. We haven't tested it yet.\\ Now we introduce Leray-Schauder degree.
\begin{definition}\label{definition}$\Omega_{0}$\ is bounded open set of real Banach space\ $B,\ T:\ \overline{\Omega_{0}}\rightarrow B$\ is totally continuous,\ $f(x)=x-T(x),\ \forall x\in \overline{\Omega_{0}},\ p\in B\setminus f(\partial\Omega_{0})$, \[ \tau=\inf_{x\in\partial\Omega_{0}}\|f(x)-p\|>0.\] There exists\ $B^{(n)}$\ is subspace of\ $B$\ with finite dimensions,\ $p\in B^{(n)}$, and there exists bounded continuous operator\ $T_{n}:\ \overline{\Omega_{0}}\rightarrow B^{(n)}$, such that\[ \|T(x)-T_{n}(x)\|<\tau,\ \forall x\in \overline{\Omega_{0}}.\] Then Leray-Schauder degree of totally continuous field\ $f$\ is
\[ deg(f,\ \Omega_{0},\ p)=deg(f_{n},\ \Omega_{0,\ n},\ p),\]where\ $f_{n}=I-T_{n},\ \Omega_{0,\ n}=B^{(n)}\cap\Omega_{0}$.\end{definition}
Maybe you are not very familiar with Leray-Schauder degree or even you know nothing about it. That's not the problem. We only apply three primary theorems to\ $Z_{1}=T_{0\epsilon}(Z_{1}),\ \forall \epsilon>0$. We write them together into a lemma as follows.
\begin{lemma} \label{lemma1}(1)(Kronecker) If\ $deg(f,\ \Omega_{0},\ p)\neq 0$, then there exists solution for\ $f(x)=p$\ in\ $\Omega_{0}$.\\
(2)(Rothe) If\ $\Omega_{0}$\ is bounded and open convex set in Banach space\ $B$,\ $T: \overline{\Omega_{0}}\rightarrow B$\ is totally continuous,\ $T(\partial\Omega_{0})\subset \overline{\Omega_{0}},$\ and\ $T(x)\neq x,\ \forall x\in \partial\Omega_{0}$, then\ $deg(I-T,\ \Omega_{0},\ 0)\neq 0$.\\
(3) If\ $f=I-T$\ and\ $f_{1}=I-T_{1}$\ are all totally continuous fields mapping from\ $\overline{\Omega_{0}}$\ to Banach space\ $B$,\ $p$\ is not in\ $f(\partial\Omega_{0})\cup f_{1}(\partial\Omega_{0})$, moreover\[ \|T_{1}(x)-T(x)\|\leq \|x-T(x)-p\|,\ \forall x\in \partial\Omega_{0},\]
then\ $deg(f,\ \Omega_{0},\ p)=deg(f_{1},\ \Omega_{0},\ p)$.
\end{lemma}
(3) is called homotopic. You may read the explanation of definition and all the proofs of three primary theorems from page 135 to page 165 in [8]. We don't repeat them again.\\
We denote as follows, \[\tau(M,\ \epsilon,\ T)=\inf_{\|Z_{1}\|_{\infty}=M}\|Z_{1}-T_{0\epsilon}(Z_{1})\|_{\infty}.\]
If time\ $T$\ is fixed, then we denote\ $\tau(M,\ \epsilon,\ T)$\ into\ $\tau(M,\ \epsilon)$.\\
If\ $\tau(M,\ \epsilon,\ T)=0$, then there exists a sequence\ $Z_{1k},\ \|Z_{1k}\|_{\infty}=M,\ k\geq1$, such that
\[\lim_{k\rightarrow+\infty}\|Z_{1k}-T_{0\epsilon}(Z_{1k})\|_{\infty}=0.\]
Because\ $T_{0\epsilon}$\ is compact, there exist a sub-sequence\ $Z_{1n_{k}},\ k\geq1,\ \mbox{and}\ Z_{1\epsilon}\in C(\overline{\Omega}\times[0,\ T])$, such that\[\lim_{k\rightarrow+\infty}\|Z_{1\epsilon}-T_{0\epsilon}(Z_{1n_{k}})\|_{\infty}=0.\] From (4.30), we obtain\[\lim_{k\rightarrow+\infty}\|Z_{1n_{k}}-T_{0\epsilon}(Z_{1n_{k}})\|_{\infty}=0.\]
And we can get that\[\lim_{k\rightarrow+\infty}\|Z_{1\epsilon}-Z_{1n_{k}}\|_{\infty}=0.\] Because\ $T_{0\epsilon}$\ is continuous, we can obtain
\[\lim_{k\rightarrow+\infty}\|T_{0\epsilon}(Z_{1\epsilon})-T_{0\epsilon}(Z_{1n_{k}})\|_{\infty}=0.\] (4.31) and (4.34) mean that\ $Z_{1\epsilon}=T_{0\epsilon}(Z_{1\epsilon})$, where\ $Z_{1\epsilon}\in C(\overline{\Omega}\times[0,\ T]),\ \|Z_{1\epsilon}\|_{\infty}=M.$\\
Now we see the solution of\ $Z_{1\epsilon}=T_{0\epsilon}(Z_{1\epsilon})$\ as follows.
\begin{theorem} \label{Theorem4-1}(Local existence) If the following condition stands,
\[ \exists M>0,\ \exists \delta>0,\ \exists \delta^{\prime}>0,\ \forall \epsilon\in (0,\ \delta],\ \forall T\in (0,\ \delta^{\prime}],\ \mbox{we have}\ \tau(M,\ \epsilon,\ T)>0,\] then\ $\exists \delta^{\prime\prime}\in (0,\ \delta^{\prime}],\ \forall T\in (0,\ \delta^{\prime\prime}],\ \forall \epsilon\in (0,\ \delta]$, there exists\ $Z_{1\epsilon}\in C(\overline{\Omega}\times[0,\ T]),\ \|Z_{1\epsilon}\|_{\infty}<M$\, such that\ $Z_{1\epsilon}=T_{0\epsilon}(Z_{1\epsilon})$.\end{theorem}
{\it Proof of theorem 4.1}. If we denote\ $\Omega_{M}=\{Z_{1}\in C(\overline{\Omega}\times[0,\ T]): \|Z_{1}\|_{\infty}<M\}$, then we will see\ $\forall T\in (0,\ \delta^{\prime}],\  deg(Z_{1}-T_{0\epsilon}(Z_{1}),\ \Omega_{M},\ 0)$\ keeping constant,\ $\forall \epsilon\in (0,\ \delta]$.\\If we select\ $ \epsilon_{0}\in (0,\ \delta]$, then\ $\forall \epsilon_{1}\in [\epsilon_{0},\ \delta]$, according to (4) in lemma 4.1, we have\ $\exists \delta(\epsilon_{1})>0$, such that$$ \|[T_{0\epsilon_{2}}(Z_{1})-T_{0\epsilon_{1}}(Z_{1})]I_{\overline{\Omega}\times[0,\ T]}\|_{\infty}\leq \tau(M,\ \epsilon_{1},\ T),\ \forall\ Z_{1}\in \Omega_{M},\ \forall  \epsilon_{2}\in U(\epsilon_{1},\ \delta(\epsilon_{1})).  $$
 And from (3) in lemma 3.2, we obtain$$deg(Z_{1}-T_{0\epsilon_{2}}(Z_{1}),\ \Omega_{M},\ 0)=deg(Z_{1}-T_{0\epsilon_{1}}(Z_{1}),\ \Omega_{M},\ 0),\ \forall  \epsilon_{2}\in U(\epsilon_{1},\ \delta(\epsilon_{1})).  $$
Hence\ $ deg(Z_{1}-T_{0\epsilon_{2}}(Z_{1}),\ \Omega_{M},\ 0)$\
keep constant,\ $\forall \epsilon_{2}\in U(\epsilon_{1},\ \delta(\epsilon_{1}))$.\\ We can see that\ $U(\epsilon_{1},\ \delta(\epsilon_{1})),\ \forall \epsilon_{1}\in [\epsilon_{0},\ \delta]$\ is an open cover for\ $[\epsilon_{0},\ \delta]$. From Heine-Borel theorem, there exists finite sub-cover. This means that\ $ deg(Z_{1}-T_{0\epsilon_{1}}(Z_{1}),\ \Omega_{M},\ 0)$\ stays constant,\ $\forall \epsilon_{1}\in [\epsilon_{0},\ \delta]$. From the arbitrary nature of\ $ \epsilon_{0}$, we know\ $ deg(Z_{1}-T_{0\epsilon}(Z_{1}),\ \Omega_{M},\ 0)$\ keep constant,\ $\forall \epsilon\in (0,\ \delta]$.\\ Now we choose\ $T$\ is sufficiently small, such that\ $T_{0\delta}$\ is a contract mapping. \\There exists\ $\delta^{\prime\prime}\in (0,\ \delta^{\prime}],\ \forall T\in (0,\ \delta^{\prime\prime}]$,\ $T_{0\delta}(\partial\Omega_{M})\subset\overline{\Omega_{M}}$. From Rothe theorem, (2) in lemma 4.2, we get\ $deg(Z_{1}-T_{0\delta}(Z_{1}),\ \Omega_{M},\ 0)\neq 0$. \\Hence we obtain\ $deg(Z_{1}-T_{0\epsilon}(Z_{1}),\ \Omega_{M},\ 0)\neq 0,\ \forall \epsilon\in (0,\ \delta]$.\\
From Kronecker theorem, (1) in lemma 4.2, we know there exists\ $Z_{1\epsilon}\in \Omega_{M}$, such that\ $Z_{1\epsilon}=T_{0\epsilon}(Z_{1\epsilon}),\ \forall \epsilon\in (0,\ \delta]$.\qed\\
 If (4.35) is not true, then the following will stand,
\[ \forall M>0,\ \exists \epsilon_{k}>0,\ \exists T_{k}>0,\ \lim_{k\rightarrow+\infty}\epsilon_{k}=\lim_{k\rightarrow+\infty}T_{k}=0,\ \mbox{such that}\ \tau(M,\ \epsilon_{k},\ T_{k})\equiv 0,\ \forall k\geq1.\]
This means that\ $\forall M>0$, there exists\ $Z_{1}(M,\ \epsilon_{k},\ T_{k})\in C(\overline{\Omega}\times[0,\ T]),\ k\geq1$, such that
\[ Z_{1}(M,\ \epsilon_{k},\ T_{k})=T_{0\epsilon_{k}}(Z_{1}(M,\ \epsilon_{k},\ T_{k})),\ \|Z_{1}(M,\ \epsilon_{k},\ T_{k})\|_{\infty}=M,\ \forall k\geq1.\]
That's not easy. And if we take\[M_{2}-M_{1}\geq M,\ Z_{1}(M_{1},\ \epsilon_{1},\ T_{1})=Z_{1}(\epsilon_{1},\ T_{1}),\ Z_{1}(M_{2},\ \epsilon_{2},\ T_{2})=Z_{1}(\epsilon_{2},\ T_{2}),\] then we can get as follows,\ $\forall M>0,\ \forall \delta>0,\ \forall \delta^{\prime}>0,\ \exists \epsilon_{1},\ \epsilon_{2}\in (0,\ \delta],\ \exists T_{1},\ T_{2}\in (0,\ \delta^{\prime}]$, such that\[\|Z_{1}(\epsilon_{1},\ T_{1})-Z_{1}(\epsilon_{2},\ T_{2})\|_{\infty}\geq M,\] where\ $Z_{1}(\epsilon_{j},\ T_{j})=T_{0\epsilon_{j}}(Z_{1}(\epsilon_{j},\ T_{j})),\ j=1,\ 2$. From (4.39) it looks something related to blow-up is happening.
\begin{corollary}(Global existence) If the following condition stands,
\[ \exists M>0,\ \exists \delta>0,\ \forall \epsilon\in (0,\ \delta],\ \tau(M,\ \epsilon)>0,\ \mbox{and}
\ \exists \epsilon_{0}\in (0,\ \delta],\ deg(Z_{1}-T_{0\epsilon_{0}}(Z_{1}),\ \Omega_{M},\ 0)\neq0,\] then\ $\forall \epsilon\in (0,\ \delta]$, there exists\ $Z_{1\epsilon}\in \Omega_{M}$, such that\ $Z_{1\epsilon}=T_{0\epsilon}(Z_{1\epsilon})$.\end{corollary}
{\it Proof of corollary 4.1}. We can obtain the proof by Theorem 4.1.\qed\\Next we see the solution of\ $Z_{1}^{\ast}=T_{0}(Z_{1}^{\ast})$\ is the following.
\begin{theorem} \label{Theorem4-2} (1)(Strong solution) $\forall \epsilon_{k}>0,\ \epsilon_{k}\rightarrow0,\ k\rightarrow+\infty,\ Z_{1\epsilon_{k}}\in \Omega_{M}$,\ $Z_{1\epsilon_{k}}=T_{0\epsilon_{k}}(Z_{1\epsilon_{k}}),\ k\geq1$, there exist sub-series\ $n_{k},\ k\geq1$, and\ $Z_{1}^{\ast}\in H^{-m_{1}}(\Omega_{1})$, such that\[ \lim_{k\rightarrow+\infty}\|Z_{1\epsilon_{n_{k}}}-Z_{1}^{\ast}\|_{-m_{1}}=0,\ \mbox{moreover}\ \lim_{k\rightarrow+\infty}\|Z_{1\epsilon_{n_{k}}}-T_{0}(Z_{1\epsilon_{n_{k}}})\|_{-m_{1}}=0.\]If\ $Z_{1}^{\ast}$\ is locally integrable, then\ $\|Z_{1}^{\ast}\|_{L^{\infty}}\leq M$, where\[\|Z_{1}^{\ast}\|_{L^{\infty}}=\max_{1\leq i\leq33}\|Z_{1,\ i}^{\ast}\|_{L^{\infty}}.\](2)($L^{\infty}$\ solution)If there exist\ $\epsilon_{k}>0,\ Z_{1\epsilon_{k}}\in \Omega_{M}$,\ $Z_{1\epsilon_{k}}=T_{0\epsilon_{k}}(Z_{1\epsilon_{k}}),\ k\geq1$, moreover\[\lim_{k\rightarrow+\infty}\epsilon_{k}=0,\ \lim_{k\rightarrow+\infty}Z_{1\epsilon_{k}}\mbox{exists almost everywhere on}\ \overline{\Omega}\times[0,\ T],\] then there exists\ $Z_{1}^{\ast}\in L^{\infty}(\overline{\Omega}\times[0,\ T]),\ \|Z_{1}^{\ast}\|_{L^{\infty}}\leq M$, such that\ $Z_{1}^{\ast}=T_{0}(Z_{1}^{\ast})$.\end{theorem}
{\it Proof of theorem 4.2}. (1)From (3) in lemma 3.1, we can get that there exist sub-series\ $n_{k},\ k\geq1$, and\ $Z_{1}^{\ast}\in H^{-m_{1}}(\Omega_{1})$, such that\[\lim_{k\rightarrow +\infty}\epsilon_{n_{k}}=0,\ \lim_{k\rightarrow +\infty}\|T_{0}(Z_{1\epsilon_{n_{k}}})-Z_{1}^{\ast}\|_{-m_{1}}=0.\] From\ $Z_{1\epsilon_{n_{k}}}=T_{0\epsilon_{n_{k}}}(Z_{1\epsilon_{n_{k}}}),\ k\geq1$, and (2) in lemma 3.1, we can obtain (4.41) holds.\\
If\ $\|Z_{1}^{\ast}\|_{L^{\infty}}> M$, then there exists\ $Z_{1,\ i}^{\ast}$, such that\ $\|Z_{1,\ i}^{\ast}\|_{L^{\infty}}>M$.\\
Hence there exists\ $\epsilon_{0}>0$, such that\ $m(\Omega_{i}(\epsilon_{0}))>0$, where\[\Omega_{i}(\epsilon_{0})=\{ X\in \overline{\Omega}\times[0,\ T]:\ |Z_{1,\ i}^{\ast}(X)|\geq M+\epsilon_{0}\}.\] Because\ $m(\Omega_{i}(\epsilon_{0}))=m(\Omega_{i}^{+}(\epsilon_{0}))+m(\Omega_{i}^{-}(\epsilon_{0}))$, where
\begin{eqnarray*}\Omega_{i}^{+}(\epsilon_{0})&=&\{ X\in \overline{\Omega}\times[0,\ T]:\ Z_{1,\ i}^{\ast}(X)\geq M+\epsilon_{0}\},\\ \Omega_{i}^{-}(\epsilon_{0})&=&\{ X\in \overline{\Omega}\times[0,\ T]:\ Z_{1,\ i}^{\ast}(X)\leq -(M+\epsilon_{0})\},\end{eqnarray*} we know at least one of\ $m(\Omega_{i}^{+}(\epsilon_{0})),\ m(\Omega_{i}^{-}(\epsilon_{0}))$\ is bigger than\ $0$.\\ We assume as well\ $m(\Omega_{i}^{+}(\epsilon_{0}))>0$. If we choose\ $\varphi_{0}\in C_{0}^{\infty}(\Omega_{i}^{+}(\epsilon_{0}))$, and\ $\varphi_{0}\geq0,\ supp(\varphi_{0})\neq\emptyset$, then we can get that\[\|Z_{1\epsilon_{n_{k}}}-Z_{1}^{\ast}\|_{-m_{1}}\geq \cfrac{\int_{\Omega_{i}^{+}(\epsilon_{0})}\epsilon_{0}\varphi_{0}dX}{\|\varphi_{0}\|_{m_{1}}}>0.\] That's contradict with (4.41).\\
(2)From (4.43), we know there exist\ $Z_{1}^{\ast}$, such that\[ \lim_{k\rightarrow +\infty}Z_{1\epsilon_{k}}=Z_{1}^{\ast},\ \mbox{a.e.}\]
Because\ $Z_{1\epsilon_{k}}\in \Omega_{M}$, we get\ $Z_{1}^{\ast}\in L^{\infty}(\overline{\Omega}\times[0,\ T])$,
moreover\ $\|Z_{1}^{\ast}\|_{L^{\infty}}\leq M.$\\ From\ $\|Z_{1\epsilon_{k}}-Z_{1}^{\ast}\|_{-m_{1}}\leq \|Z_{1\epsilon_{k}}-Z_{1}^{\ast}\|_{L^{2}}$,
where\[\|Z_{1\epsilon_{k}}-Z_{1}^{\ast}\|_{L^{2}}=\max_{1\leq i\leq 33}\|Z_{1,\ i,\ \epsilon_{k}}-Z_{1,\ i}^{\ast}\|_{L^{2}},\]and Lebesgue dominated convergence theorem, we obtain\[\lim_{k\rightarrow +\infty}\|Z_{1\epsilon_{k}}-Z_{1}^{\ast}\|_{-m_{1}}=0.\]
At last we see\ $\|T_{0}(Z_{1\epsilon_{k}})-T_{0}(Z_{1}^{\ast})\|_{-m_{1}}$. From
\[<T_{0}(Z_{1\epsilon_{k}})-T_{0}(Z_{1}^{\ast}),\ \varphi>=<F[T_{0}(Z_{1\epsilon_{k}})-T_{0}(Z_{1}^{\ast})],\ F^{-1}(\varphi)>,\]
we only need to discuss\ $\psi(Z_{1\epsilon_{k}})-\psi(Z_{1}^{\ast})$\ in\ $F[T_{0}(Z_{1\epsilon_{n_k}})-T_{0}(Z_{1}^{\ast})]$. Because\ $\psi$\ is continuous, we can get there exists a constant\ $C_{M}>0$, such that
\[ \|\psi(Z_{1\epsilon_{k}})-\psi(Z_{1}^{\ast})\|_{L^{\infty}}\leq C_{M}.\]
Again from Lebesgue dominated convergence theorem, we obtain
\[\lim_{k\rightarrow +\infty}\|T_{0}(Z_{1\epsilon_{k}})-T_{0}(Z_{1}^{\ast})\|_{-m_{1}}=0.\]Together with\ $Z_{1\epsilon_{k}}=T_{0\epsilon_{k}}(Z_{1\epsilon_{k}})$, we obtain\ $\|Z_{1}^{\ast}-T_{0}(Z_{1}^{\ast})\|_{-m_{1}}=0$. This means\ $Z_{1}^{\ast}=T_{0}(Z_{1}^{\ast}),\ \mbox{a.e.}$\ which completed the statement.\qed\\
From the previous theorem and corollary, we know that the strong solution will exist locally under the condition (4.35) and exist globally under the condition (4.40).\\
Finally, we discuss a little more for\ $L^{\infty}$\ solution as follows.
\begin{theorem} \label{Theorem4-3} A necessary and sufficient condition for (4.43) holding is that there exist\ $\eta_{k}>0,\ Z_{1\eta_{k}}\in \Omega_{M}$,\ $Z_{1\eta_{k}}=T_{0\eta_{k}}(Z_{1\eta_{k}}),\ k\geq1$, moreover\[\lim_{k\rightarrow+\infty}\eta_{k}=0,\ \lim_{k,\ l\rightarrow+\infty}\|Z_{1\eta_{k}}-Z_{1\eta_{l}}\|_{L^{2}}=0.\]\end{theorem}
{\it Proof of theorem 4.3}. Necessity. If (4.43) holds, then there exists\ $Z_{1}^{\ast}\in L^{\infty}(\overline{\Omega}\times[0,\ T])$, such that
\[\lim_{k\rightarrow +\infty}\|Z_{1\epsilon_{k}}-Z_{1}^{\ast}\|_{L^{2}}=0.\] If we let\ $\eta_{k}=\epsilon_{k},\ k\geq1$, then (4.53) stands.\\
Sufficiency. If (4.53) holds, then\ $Z_{1\eta_{k}},\ k\geq 1$, are convergent by the Lebesgue measure as follows,
\[ \forall \epsilon>0,\ \lim_{k,\ l\rightarrow +\infty}m(\{X\in \overline{\Omega}\times[0,\ T]: \|Z_{1\eta_{k}}-Z_{1\eta_{l}}\|_{\infty}\geq \epsilon\})=0.\]
From the Riesz theorem on page 142 in [14], we will see that there exist sub-series\ $n_{k},\ k\geq1$, such that
\[ \lim_{k\rightarrow+\infty}\eta_{n_{k}}=0,\ \lim_{k\rightarrow+\infty}Z_{1\eta_{n_{k}}}\mbox{exists almost everywhere on}\ \overline{\Omega}\times[0,\ T].\]
If we let\ $\epsilon_{k}=\eta_{n_{k}}$, then (4.43) stands.\qed\\
We know\ $L^{2}(\Omega_{1})$\ is not completed under norm\ $\|\cdot\|_{-m_{1}}$. After being completed, it will be\ $H^{-m_{1}}(\Omega_{1})$. So there exists\ $f_{k}\in L^{2}(\Omega_{1}),\ k\geq1$, such that\[ \lim_{k,\ l\rightarrow +\infty}\|f_{k}-f_{l}\|_{-m_{1}}=0,\ \mbox{and}\ \|f_{k}-f_{l}\|_{L^{2}}\geq c>0,\ \forall\ k\neq l.\]
From (4.53), we know that (4.43) will not always stand. But we can see that it holds in many cases such as the following.
\begin{theorem} \label{Theorem4-4}A sufficient condition for (4.43) holding is that there exist\ $\epsilon_{k}>0,\ \epsilon_{k}\rightarrow0,\ k\rightarrow+\infty,\ Z_{1\epsilon_{k}}\in \Omega_{M}$,\ $Z_{1\epsilon_{k}}=T_{0\epsilon_{k}}(Z_{1\epsilon_{k}}),\ k\geq1$, such that\ $\forall c>0,\ \exists d>0,\ \forall k,\ l,\ i,\ 1\leq i\leq 33$, we have\[m(\Omega_{k,\ l,\ i}^{+}(d)\cup \Omega_{k,\ l,\ i}^{-}(d))\leq c,\]
where\begin{eqnarray*}\Omega_{k,\ l,\ i}^{+}&=&\{X\in \overline{\Omega}\times[0,\ T]:Z_{1,\ i,\ \epsilon_{k}}-Z_{1,\ i,\ \epsilon_{l}}\geq 0\},\\
\Omega_{k,\ l,\ i}^{-}&=&\{X\in \overline{\Omega}\times[0,\ T]:Z_{1,\ i,\ \epsilon_{k}}-Z_{1,\ i,\ \epsilon_{l}}< 0\},\\ \Omega_{k,\ l,\ i}^{+}(d)&=&\{X\in \overline{\Omega}\times[0,\ T]:\ dist(X,\ \partial \Omega_{k,\ l,\ i}^{+})\leq d\},\\ \Omega_{k,\ l,\ i}^{-}(d)&=&\{X\in \overline{\Omega}\times[0,\ T]:\ dist(X,\ \partial \Omega_{k,\ l,\ i}^{-})\leq d\}.\end{eqnarray*}\end{theorem}
{\it Proof of theorem 4.4}. From Theorem 3.2, we know there exist sub-series\ $n_{k},\ k\geq1$, such that
\[\lim_{k,\ l\rightarrow +\infty}\|Z_{1\epsilon_{n_{k}}}-Z_{1\epsilon_{n_{l}}}\|_{-m_{1}}=0.\]
We will prove that\ $\forall i,\ 1\leq i\leq 33,\ Z_{1,\ i,\ \epsilon_{n_{k}}},\ k\geq1$, is convergent by the Lebesgue measure.\\
If that is not true, then there exists\ $i_{0},\ 1\leq i_{0}\leq 33,\ \exists a_{2}>0,\ \exists b_{2}>0$, and sub-series\ $k_{j},\ l_{j},\ j\geq1$, such that
\[ m(\{X\in \overline{\Omega}\times[0,\ T]:\ |Z_{1,\ i_{0},\ \epsilon_{n_{k_{j}}}}-Z_{1,\ i_{0},\ \epsilon_{n_{l_{j}}}}|\geq a_{2}\})\geq b_{2},\ j\geq1.\]
We denote this in an easy way as follows,\ $\forall j,\ j\geq1$,\begin{eqnarray*}f_{j}&=&Z_{1,\ i_{0},\ \epsilon_{n_{k_{j}}}}-Z_{1,\ i_{0},\ \epsilon_{n_{l_{j}}}},\\
\Omega_{j}^{+}&=&\{X\in \overline{\Omega}\times[0,\ T]:f_{j}\geq 0\},\\
\Omega_{j}^{-}&=&\{X\in \overline{\Omega}\times[0,\ T]:f_{j}< 0\},\\ \Omega_{j}^{+}(d)&=&\{X\in \overline{\Omega}\times[0,\ T]:\ dist(X,\ \partial \Omega_{j}^{+})\leq d\},\\ \Omega_{j}^{-}(d)&=&\{X\in \overline{\Omega}\times[0,\ T]:\ dist(X,\ \partial \Omega_{j}^{-})\leq d\}.\end{eqnarray*}
From (4.58), we can obtain that there exists\ $d_{0}>0$, such that
\[ m(\Omega_{j}^{+}(d_{0})\cup\Omega_{j}^{-}(d_{0}))\leq \cfrac{a_{2}b_{2}}{4M}.\]
If we let\ $\varphi_{j}\in C_{0}^{\infty}(\Omega_{1}),\ j\geq1$, as follows,
\[ \varphi_{j}(X)=\int_{\Omega_{j}^{+}\setminus \Omega_{j}^{+}(d_{0})} \alpha_{d_{0}}(X-Y)dY-\int_{\Omega_{j}^{-}\setminus \Omega_{j}^{-}(d_{0})} \alpha_{d_{0}}(X-Y)dY,\ j\geq1,\]
where$$\alpha_{d_{0}}(X)=\cfrac{1}{d_{0}^{4}}\alpha(\cfrac{X}{d_{0}}),\ \alpha(X)= \begin{cases}
                                                                                                                              Ce^{1/(|X|^{2}-1)},\ |X|<1, \\
                                                                                                                              0,\ |X|\geq 1.
                                                                                                                   \end{cases} ,\ C=(\int_{|X|<1}e^{1/(|X|^{2}-1)}dX)^{-1},$$
then\ $\varphi_{j}=1$\ on\ $\Omega_{j}^{+}\setminus \Omega_{j}^{+}(d_{0})$,\ $\varphi_{j}\in [0,\ 1]$\ on\ $\Omega_{j}^{+}(d_{0})$, and\ $\varphi_{j}=-1$\ on\ $\Omega_{j}^{-}\setminus \Omega_{j}^{-}(d_{0})$,\ $\varphi_{j}\in [-1,\ 0]$\ on\ $\Omega_{j}^{-}(d_{0})$.
\\Moreover, we can get as follows, \[ |\partial^{\gamma}\varphi_{j}|\leq 2\int_{\Omega_{1}}|\partial^{\gamma}\alpha_{d_{0}}(X-Y)|dY,\ 0\leq |\gamma|\leq m_{1}.\]
Hence, there exists\ $d_{1}>0$, such that\[ \|\varphi_{j}\|_{m_{1}}\leq d_{1},\ \forall j,\ j\geq1.\]
From $$ |<f_{j},\ \varphi_{j}>|\geq|<f_{j},\ sign f_{j}>|-|<f_{j},\ sign f_{j}-\varphi_{j}>|\geq a_{2}b_{2}-2M\cfrac{a_{2}b_{2}}{4M}=\cfrac{a_{2}b_{2}}{2},$$where
$$ sign f_{j}=\begin{cases}
1,\ f_{j}\geq0,\\
-1,\ f_{j}<0. \end{cases},$$ we can obtain that\[ \|f_{j}\|_{-m_{1}}\geq \cfrac{|<f_{j},\ \varphi_{j}>|}{\|\varphi_{j}\|_{m_{1}}}\geq \cfrac{a_{2}b_{2}}{2d_{1}}.\]
This contradicts (4.59).\\
So\ $\forall i,\ 1\leq i\leq 33,\ Z_{1,\ i,\ \epsilon_{n_{k}}},\ k\geq1$, is convergent by the Lebesgue measure. From the Riesz theorem on page 142 in [7], we know that (4.43) holds. \qed\\ If (4.58) is not true, then\ $\forall \epsilon_{k}>0,\ \epsilon_{k}\rightarrow0,\ k\rightarrow+\infty,\ Z_{1\epsilon_{k}}\in \Omega_{M}$,\ $Z_{1\epsilon_{k}}=T_{0\epsilon_{k}}(Z_{1\epsilon_{k}}),\ k\geq1$,\ $\exists c>0,\ \forall d_{j}>0,\ \exists k_{j},\ l_{j},\ i_{j},\ 1\leq i_{j}\leq 33$, such that\[m(\Omega_{k_{j},\ l_{j},\ i_{j}}^{+}(d_{j})\cup \Omega_{k_{j},\ l_{j},\ i_{j}}^{-}(d_{j}))\geq c.\]
If we let\ $ d_{j}\rightarrow 0,\ j\rightarrow +\infty$, then we can get the following,\[ \lim_{j\rightarrow +\infty} S(\partial\Omega_{k_{j},\ l_{j},\ i_{j}}^{+})=+\infty,\] where $$S(\partial\Omega_{k_{j},\ l_{j},\ i_{j}}^{+})=\int_{\{X\in \overline{\Omega}\times[0,\ T]:\ Z_{1,\ i_{j},\ \epsilon_{k_{j}}}-Z_{1,\ i_{j},\ \epsilon_{l_{j}}}= 0\}}dS. $$ And the following is another sufficient condition for (4.43) holding,
\[ \exists\ \epsilon_{k}>0,\ Z_{1\epsilon_{k}}\in \Omega_{M},\ Z_{1\epsilon_{k}}=T_{0\epsilon_{k}}(Z_{1\epsilon_{k}}),\ k\geq1,\ \lim_{k\rightarrow +\infty}\epsilon_{k}=0,\ \mbox{and}\ \sup_{k,\ l,\ i}S(\partial\Omega_{k,\ l,\ i}^{+})< +\infty,\]where $$S(\partial\Omega_{k,\ l,\ i}^{+})=\int_{\{X\in \overline{\Omega}\times[0,\ T]:\ Z_{1,\ i,\ \epsilon_{k}}-Z_{1,\ i,\ \epsilon_{l}}= 0\}}dS,\ k\neq l,\ 1\leq i\leq 33. $$
We may define blow-up if one of the following happens,\\(1)(4.35) is not true,
\[ \forall M>0,\ \exists \epsilon_{k}>0,\ \exists T_{k}>0,\ \lim_{k\rightarrow+\infty}\epsilon_{k}=\lim_{k\rightarrow+\infty}T_{k}=0,\ \mbox{such that}\ \tau(M,\ \epsilon_{k},\ T_{k})\equiv 0,\ \forall k\geq1.\]
(2)(4.40) is not true,
\[ \forall M>0,\ \forall \delta>0,\ \exists \epsilon\in (0,\ \delta],\ \tau(M,\ \epsilon)=0,\ \mbox{or}
\ \forall \epsilon_{0}\in (0,\ \delta],\ deg(Z_{1}-T_{0\epsilon_{0}}(Z_{1}),\ \Omega_{M},\ 0)\equiv 0,\]
(3)The strong solution is not unique. If there exist at least two strong solutions\ $Z_{1}^{\ast},\ Z_{1}^{\ast\prime}$, then we can get that
\begin{eqnarray}&& \lim_{k\rightarrow+\infty}\epsilon_{k}=0,\ \lim_{k\rightarrow+\infty}\|Z_{1\epsilon_{k}}-Z_{1}^{\ast}\|_{-m_{1}}=0,\ \lim_{k\rightarrow+\infty}\|Z_{1\epsilon_{k}}-T_{0}(Z_{1\epsilon_{k}})\|_{-m_{1}}=0,\\&&  \lim_{k\rightarrow+\infty}\eta_{k}=0,\ \lim_{k\rightarrow+\infty}\|Z_{1\eta_{k}}-Z_{1}^{\ast\prime}\|_{-m_{1}}=0,\ \lim_{k\rightarrow+\infty}\|Z_{1\eta_{k}}-T_{0}(Z_{1\eta_{k}})\|_{-m_{1}}=0.\end{eqnarray}
If we assuming that\[ \epsilon_{0}=\|Z_{1}^{\ast}- Z_{1}^{\ast\prime}\|_{-m_{1}},\] then there exists\ $N>0$,\ $\forall\ k\geq N$, we have
\[ \|Z_{1\epsilon_{k}}-Z_{1}^{\ast}\|_{-m_{1}}\leq \cfrac{\epsilon_{0}}{4},\ \|Z_{1\eta_{k}}-Z_{1}^{\ast\prime}\|_{-m_{1}}\leq\cfrac{\epsilon_{0}}{4}.\]
But we can obtain the following,
\[ \|Z_{1\epsilon_{k}}-Z_{1\eta_{k}}\|_{-m_{1}}\geq \|Z_{1}^{\ast}- Z_{1}^{\ast\prime}\|_{-m_{1}}-\|Z_{1\epsilon_{k}}-Z_{1}^{\ast}\|_{-m_{1}}-\|Z_{1\eta_{k}}-Z_{1}^{\ast\prime}\|_{-m_{1}}\geq\cfrac{\epsilon_{0}}{2}.\]
This means that\ $ \exists \epsilon_{0}>0,\ \forall \epsilon>0,\ \exists \epsilon_{1}>0,\ \epsilon_{2}>0,\ |\epsilon_{1}-\epsilon_{2}|\leq \epsilon,\ Z_{1\epsilon_{1}}=T_{0\epsilon_{1}}(Z_{1\epsilon_{1}}),\ Z_{1\epsilon_{2}}=T_{0\epsilon_{2}}(Z_{1\epsilon_{2}})$,\ such that
\[\|Z_{1\epsilon_{1}}-Z_{1\epsilon_{2}}\|_{-m_{1}}\geq\cfrac{\epsilon_{0}}{2}.\]
This also means that the solution of the equation\ $Z_{1\epsilon}=T_{0\epsilon}(Z_{1\epsilon})$\ is not stable on\ $\epsilon$.\\
(4)(4.58) is not true,\ $\forall \epsilon_{k}>0,\ \epsilon_{k}\rightarrow0,\ k\rightarrow+\infty,\ Z_{1\epsilon_{k}}\in \Omega_{M}$,\ $Z_{1\epsilon_{k}}=T_{0\epsilon_{k}}(Z_{1\epsilon_{k}}),\ k\geq1$,\ $\exists c>0,\ \forall d_{j}>0,\ \exists k_{j},\ l_{j},\ i_{j},\ 1\leq i_{j}\leq 33$, such that\[m(\Omega_{k_{j},\ l_{j},\ i_{j}}^{+}(d_{j})\cup \Omega_{k_{j},\ l_{j},\ i_{j}}^{-}(d_{j}))\geq c.\]
If we let\ $ d_{j}\rightarrow 0,\ j\rightarrow +\infty$, then we can get the following,\[ \lim_{j\rightarrow +\infty} S(\partial\Omega_{k_{j},\ l_{j},\ i_{j}}^{+})=+\infty,\] where $$S(\partial\Omega_{k_{j},\ l_{j},\ i_{j}}^{+})=\int_{\{X\in \overline{\Omega}\times[0,\ T]:\ Z_{1,\ i_{j},\ \epsilon_{k_{j}}}-Z_{1,\ i_{j},\ \epsilon_{l_{j}}}= 0\}}dS. $$
Each of them deserves to be discussed more carefully.\\
If there is no blow-up, then we obtain\ $u^{\ast}\in W^{2,\ +\infty}(\overline{\Omega}),\ p^{\ast}\in W^{1,\ +\infty}(\overline{\Omega})$, if\ $Z_{1}^{\ast}=T_{0}(Z_{1}^{\ast})$,\ $Z_{1}^{\ast}\in L^{\infty}(\overline{\Omega}\times[0,\ T])$, where\ $Z_{1}^{\ast}=(u^{\ast},\ p^{\ast},\ \partial u^{\ast}\setminus\ u^{\ast}_{1x},\ \partial^{2} u^{\ast}, grad p^{\ast})^{T}$. Here\ $W^{1,\ +\infty}(\overline{\Omega}),\ W^{2,\ +\infty}(\overline{\Omega})$\ are Sobolev spaces defined on page 153 in [2]. From the condition that domain\ $\Omega$\ satisfies a uniform exterior and interior cone, if\ $\Omega$\ is bounded,\ $\partial\Omega\in C^{1,\ \alpha},\ 0<\alpha\leq 1$, we can get that\ $u^{\ast}\in C^{1,\ 1}(\overline{\Omega}),\ p^{\ast}\in C^{0,\ 1}(\overline{\Omega})$\ by imbedding. By using Morrey's inequality defined on page 163 in [2], we get\ $u^{\ast}$\ is twice classically differentiable and\ $p^{\ast}$\ is classically differentiable almost everywhere in\ $\overline{\Omega}$.\\
If we looked back, we should have defined the classical solution as the\ $L^{\infty}$\ solution of\ $Z_{1}^{\ast}=T_{0}(Z_{1}^{\ast})$. It would always exist and be unique except the blow-up.\\If\ $F(T_{0}(Z_{1}^{\ast}))$\ is analytical, then\ $u^{\ast}$\ and\ $p^{\ast}$\ satisfy Eqs(1.1) and (1.2) almost everywhere in\ $\overline{\Omega}\times[0,\ T]$, where\ $F(T_{0}(Z_{1}^{\ast}))$\ is the Fourier transform of\ $T_{0}(Z_{1}^{\ast})$. That is near our goal.\\
\section{Leray-Schauder degree}\setcounter{equation}{0}
In this section, we will discuss the
the Leray-Schauder degree of nonlinear integral equation of Hammerstein type as follows,
\[ f(X)=g(X)+T(f(X)),\ \forall X\in \overline{\Omega_{1}},\]
where\[T(f(X))=\int_{\Omega_{1}}k(X,\ Y)\psi(Y,\ f(Y))dY,\]
$X=(x,\ y,\ z,\ t)^{T},\ Y=(x_{1},\ y_{1},\ z_{1},\ t_{1})^{T},\ \Omega_{1}=\Omega\times(0,\ T)$,\ $f(X)$\ is an unknown continuous function on\ $\overline{\Omega_{1}}$,\ $g(X)$\ is a known continuous function on\ $\overline{\Omega_{1}}$,\ $k(X,\ Y)$\ is a known continuous function on\ $\overline{\Omega_{1}}\times\overline{\Omega_{1}}$,\ $\psi$\ is a known continuous function on\ $\overline{\Omega_{1}}\times[-M,\ M]$.\\
From the definition of the Leray-Schauder degree on page 138 to page 139 in [8], we can work out the Leray-Schauder degree of Eq(5.1) directly as follows.\\
By the Weirstrass theorem, we know there exist two polynomials\[k_{N}(X,\ Y)=\sum_{|\alpha|=0}^{N}C_{\alpha}(Y)X^{\alpha},\ g_{N}(X)=\sum_{|\alpha|=0}^{N}g_{\alpha}X^{\alpha},\] where\ $\alpha=(\alpha_{1},\ \alpha_{2},\ \alpha_{3}, \alpha_{4})^{T},\ |\alpha|=\alpha_{1}+\alpha_{2}+\alpha_{3}+\alpha_{4},\ X^{\alpha}=x^{\alpha_{1}}y^{\alpha_{2}}z^{\alpha_{3}}t^{\alpha_{4}},\ C_{\alpha}(Y),\ 0\leq|\alpha|\leq N$, are all polynomials of\ $Y$, moreover\ $\alpha_{1},\ \alpha_{2},\ \alpha_{3},\ \alpha_{4}$\ are all nonnegative whole numbers,\ $g_{\alpha},\ 0\leq|\alpha|\leq N$,
are all real numbers, such that\ $\forall\ f(X)\in\Omega_{M}=\{f(X):\|f(X)\|_{\infty}< M\}$, we have
\[ \|\int_{\Omega_{1}}(k(X,\ Y)-k_{N}(X,\ Y))\psi(Y,\ f(Y))dY\|_{\infty}\leq \cfrac{\tau}{3},\ \|g(X)-g_{N}(X)\|_{\infty}\leq \cfrac{\tau}{3},\]
$\tau$\ is defined as follows,\[ \tau=\inf_{\|f(X)\|_{\infty}= M}\|f(X)-T(f(X))-g(X)\|_{\infty}>0.\]
From the Combination theory, the number of the solutions of\ $|\alpha|=\alpha_{1}+\alpha_{2}+\alpha_{3}+\alpha_{4}$\ is that\ $C_{3+|\alpha|}^{3}$. The number of all the items\ $X^{\alpha},\ 0\leq|\alpha|\leq N$, is that\ $L_{N}=C_{3}^{3}+C_{4}^{3}+\cdots+C_{3+N}^{3}=C_{4+N}^{4}$.\\By the homotopic, we have the following,
\[ deg(f(X)-T(f(X)),\ \Omega_{M},\ g(X))=deg(f(X)-T(f(X)),\ \Omega_{M},\ g_{N}(X)).\]
If we assume\[T_{N}(f(X))=\int_{\Omega_{1}}k_{N}(X,\ Y)\psi(Y,\ f(Y))dY,\] then from (5.4) we have\[ \|T(f(X))-T_{N}(f(X))\|_{\infty}\leq \cfrac{\tau}{3},\ \forall\ f(X)\in\Omega_{M}.\]
Moreover,\[T_{N}(f(X))=\sum_{|\alpha|=0}^{N}(\int_{\Omega_{1}}C_{\alpha}(Y)\psi(Y,\ f(Y))dY)X^{\alpha}\in E_{N},\]
where\ $E_{N}$\ is the sub-space with finite dimensions generated by\ $X^{\alpha},\ 0\leq|\alpha|\leq N$.\\From the definition of the Leray-Schauder degree, we can obtain the following,
\[deg(f(X)-T(f(X)),\ \Omega_{M},\ g_{N}(X))=deg(f(X)-T_{N}(f(X)),\ \Omega_{M,\ 1},\ g_{N}(X)),\]where\ $\Omega_{M,\ 1}=\Omega_{M}\cap E_{N}$.\\
If we denote\begin{eqnarray} f(X)&=&\tilde{X}^{T}D_{N},\ \tilde{X}=(X^{\alpha},\ 0\leq|\alpha|\leq N)^{T},\\
D_{N}&=&(D_{\alpha},\ 0\leq|\alpha|\leq N)^{T}\in R^{L_{N}},\end{eqnarray} \begin{eqnarray}
T_{N}(f(X))&=&\int_{\Omega_{1}}k_{N}(X,\ Y)\psi(Y,\ f(Y))dY\\&=&\sum_{|\alpha|=0}^{N}(\int_{\Omega_{1}}C_{\alpha}(Y)\psi(Y,\ \tilde{Y}^{T}D_{N})dY)X^{\alpha},\\&=&\sum_{|\alpha|=0}^{N}\phi_{\alpha}(D_{N})X^{\alpha}=\tilde{X}^{T}\phi(D_{N}),\end{eqnarray} \begin{eqnarray}
\phi_{\alpha}(D_{N})&=&\int_{\Omega_{1}}C_{\alpha}(Y)\psi(Y,\ \tilde{Y}^{T}D_{N})dY,\\
\phi(D_{N})&=&(\phi_{\alpha}(D_{N}),\ 0\leq|\alpha|\leq N)^{T},\\ g_{N}&=&(g_{\alpha},\ 0\leq|\alpha|\leq N)^{T},\\
\Omega_{M,\ 2}&=&\{D_{N}\in R^{L_{N}}: \|\tilde{X}^{T}D_{N}\|_{\infty}< M\},
\end{eqnarray} then we obtain\ $ f(X)=T_{N}(f(X))+g_{N}(X),\ f(X)\in \Omega_{M,\ 1}$, is equivalent to
\[ D_{N}=\phi(D_{N})+g_{N},\ D_{N}\in \Omega_{M,\ 2}.\]Hence, we obtain
\[ deg(f(X)-T(f(X)),\ \Omega_{M},\ g_{N}(X))=deg(D_{N}-\phi(D_{N}),\ \Omega_{M,\ 2},\ g_{N}).\]
From (5.16), we can see that\ $\phi(D_{N})$\ is explicit. At this time, we can work out the finite dimensions Brouwer degree\ $deg(D_{N}-\phi(D_{N}),\ \Omega_{M,\ 2},\ g_{N})$\ directly. From the definition of the Brouwer degree on page 89 in [8], we only need to calculate the double integral on\ $\Omega_{M,\ 2}\times\Omega_{1}$. This is just what we want, the Leray-Schauder degree of Eq(5.1).\\As\ $\tau$\ is fixed, we don't need that\ $N$\ is sufficient big. This is the reason why we haven't discuss the priori estimation yet.\\By the same way, we can also work out the Leray-Schauder degree of Eq(5.1) even if
\[ T(f(X))=\int_{\Omega_{1}}G(X,\ Y,\ f(Y))dY,\] where\ $G$\ is a known continuous function on\ $\overline{\Omega_{1}}\times\overline{\Omega_{1}}\times[-M,\ M]$.\\
We only need to change\ $k_{N}(X,\ Y)$\ into the following,
\[ k_{N}(X,\ Y,\ f(Y))=\sum_{|\alpha|=0}^{N}C_{\alpha,\ 1}(Y,\ f(Y))X^{\alpha},\]where\ $C_{\alpha,\ 1}(Y,\ f(Y)),\ 0\leq|\alpha|\leq N$, are all polynomials of\ $Y,\ f(Y)$.\\
\section{Declarations}
Availability of supporting data\\
This paper is available to all supporting data.\\
Competing interests\\
We declare that we have no competing interests.\\
Funding\\
The author was supported by mathematics research fund of Hohai university, No. 1014-414126.\\
Authors' contributions\\
We have transformed Navier-Stokes equations into the equivalent generalized integral equations. Moreover, we discuss the existence for the classical solution by Leray-Schauder degree and Sobolev space.\\
Acknowledgements\\
We would like to give our best thanks to Prof. Mark Edelman in Yeshiva University, Prof. Caisheng Chen in Hohai University, Prof. Junxiang Xu in Southeast University, Prof. Zuodong Yang in Nanjing Normal University, Prof. Jishan Fan in Nanjing Forestry University for their guidance and other helps.\\
Authors' information\\ Jianfeng Wang, Prof. in the Mathematics department of Hohai University.\\Email: wjf19702014@163.com, Acadmic email: 20020001@hhu.edu.cn.\\
Orcid: 0000-0002-3129-756x.\\

\end{document}